\documentclass[12pt]{article}
\usepackage{subfig}

\usepackage{graphicx,color}

\usepackage{amsmath, amsfonts}
\usepackage{amsthm}
\usepackage{comment}
\usepackage{xcolor}

\usepackage[colorlinks,pdfdisplaydoctitle]{hyperref}
\newcommand{\true}{u^{\dag}}
\newcommand{\tr}{{\rm tr}}
\newcommand{\N}{\mathcal{N}}
\newcommand{\tri}[3]{{#1}_{#2}^{#3}}

\newcommand{\rr}[1]{r_{n,\alpha}(#1)}

\newcommand{\D}[1]{\mathcal{D}\left(#1\right)}

\def\E{\mathbb{E}}
\newtheorem{thm}{Theorem}{\bf}{\rm}
\newtheorem{assumption}{Assumption}{\bf}{\rm}
\newtheorem{lem}{Lemma}{\bf}{\rm}
\newtheorem{algorithm}{Algorithm}{\bf}{\rm}
\newtheorem{proposition}{Proposition}[section]
\newtheorem{remark}{Remark}[section]

\begin{document}
\title{Filter Based Methods For\\Statistical Linear Inverse Problems}

\author{Marco A. Iglesias\thanks{School of Mathematical Sciences, University of Nottingham, UK,\newline Marco.Iglesias@nottingham.ac.uk} \and Kui Lin  \thanks{%
              School of Mathematical Sciences, Fudan University, China, link10@fudan.edu.cn}
\and  Shuai Lu\thanks{%
              School of Mathematical Sciences, Fudan University, China, slu@fudan.edu.cn}
\and Andrew M. Stuart \thanks{ Mathematics Institute, University of Warwick, UK, A.M.Stuart@warwick.ac.uk}
}

\maketitle



\maketitle
\begin{abstract}
Ill-posed inverse problems are ubiquitous in applications. Understanding
of algorithms for their solution has been greatly enhanced by a deep
understanding of the linear inverse problem. In the applied communities
ensemble-based filtering methods have recently been used to solve
inverse problems by introducing an artificial dynamical system. This
opens up the possibility of using a range of other filtering methods,
such as 3DVAR and Kalman based methods, to solve inverse problems,
again by introducing an artificial dynamical system. The aim of this paper
is to analyze such methods in the context of the ill-posed linear inverse
problem.

Statistical linear inverse problems are studied in the sense that the
observational noise is assumed to be derived via realization of
a Gaussian random variable.  We investigate the asymptotic behavior of filter
based methods for these inverse problems. Rigorous
convergence rates are established for 3DVAR and for the Kalman filters,
including minimax rates in some instances. Blowup of 3DVAR and a variant
of its basic form is also presented,
and optimality of the Kalman filter is discussed.
These analyses reveal a close connection between (iterative) regularization
schemes in deterministic inverse problems and filter based methods
in data assimilation. Numerical experiments are presented to illustrate
the theory.
\end{abstract}

\section{Introduction}\label{se_Intro}
In many geophysical applications, in particular in the petroleum industry
and in hydrology, distributed parameter estimation problems are often
solved by means of iterative ensemble Kalman filters \cite{orl08}.
The basic methodology is to introduce an artificial dynamical system,
to supplement this with observations, and to apply the ensemble Kalman
filter.  The methodology is described in a basic, abstract
form, applicable to a general, possibly nonlinear,
inverse problem in \cite{ILS2013}.
In this basic form of the algorithm
regularization is present due to dynamical preservation of a
subspace spanned by the ensemble during the iteration.
The paper \cite{reg_enkf} gives further insight into the development
of regularization for these ensemble Kalman inversion methods,
drawing on links with the Levenberg-Marquardt scheme \cite{h97}.
In this paper our aim is to further the study of filters
for the solution of inverse problems, going beyond the ensemble Kalman
filter to encompass the study of other filters such as
3DVAR and the Kalman filter itself -- see \cite{LSZ2015} for an overview
of these filtering methods. A key issue will be the implementation
of regularization with the aim of deriving optimal error estimates.

We focus on the linear inverse problem
\begin{equation}
\label{eq:base}
y=A\true+\eta,
\end{equation}
where $A$ is a compact operator acting between Hilbert spaces $X$ and $Y$. The exact solution is denoted by $\true\in X$ and $\eta$ is a noise polluting
the observations. We will consider two situations: {\bf Data Model $1$}
where multiple observations are made in the form \ref{eq:base}; and
{\bf Data Model $2$} where a single observation is made.
For modelling purposes we will assume that the noise $\eta$
is generated by the Gaussian $\N(0,\gamma^2I)$,
independently in the case of multiple observations.
In each case we create a sequence $\{y_n\}_{n \ge 0}$; for Data
Model $1$ the elements of this sequence are i.i.d. $\N(A\true,\gamma^2I)$
whilst for Data Model $2$ they are $y_n \equiv y$, with $y$ a single draw
from $\N(A\true,\gamma^2I)$.
The case where multiple independent observations are made is not
uncommon in applications (for example in electrical impedance tomography
(EIT, \cite{EIT}) and, although we do not pursue it here,
our methodology also opens up the possibility of
considering multiple instances with correlated observational noise,
by means of similar filtering-based techniques.

The artificial, partially observed linear dynamical system
that underlies our methodology is as follows:
\begin{equation}\label{eq:dynamic}
    \begin{aligned}
      u_{n}&=u_{n-1},\\
      y_{n}&=Au_{n}+\eta_{n}.
    \end{aligned}
\end{equation}
In {\em deriving} the filters we apply to this dynamical system,
it is assumed that
the $\{\eta_n\}_{n \ge 0}$ are i.i.d. from $\N(0,\gamma^2I)$.
Note, however, that whilst the data sequence $\{y_n\}_{n \ge 0}$
we use in Data Model $1$ is of this form, the assumption
is not compatible with Data Model $2$; thus
for Data Model $2$ we have a form of {\em model error} or {\em model
mis-specification} \cite{LSZ2015}.

By studying the application of filtering methods to the solution of
the linear inverse problem our aim is to open up the
possibility of employing the filtering methodology to (static) inverse
problems of the form \ref{eq:base}, and nonlinear generalizations.
We confine our analysis to the linear setting as experience has shown
that a deep understanding of this case is helpful both because there
are many linear inverse problems which arise in applications,
and because knowledge of the linear case guide methodologies
for the more general nonlinear problem \cite{EHN1996}.
The last few decades have seen a comprehensive development of the
theory of linear inverse problems, both classically
and statistically -- see \cite{B07,EHN1996} and the references
therein. Consider the Tikhonov-Phillips regularization method
 \begin{equation*}
  {\rm argmin}_{u}\Bigl(\frac{1}{2\gamma^2}\|y-Au\|_Y^2+\frac{\alpha}{2}\|u-u_0\|_E^2\Bigr).
\end{equation*}
This can be reformulated from a probabilistic perspective as the MAP
estimator for Bayesian inversion given a Gaussian smoothness prior, with mean
$u_0$ and Cameron-Martin space $E$ compactly embedded into $X$,
and a Gaussian noise model as defined
above; this connection is eludicated in \cite{KS2005,Stuart2010}.
We note that
from the point of view of Tikhonov-Phillips regularization
only the parameter $\alpha\gamma^2$ is relevant, but that each of $\alpha$
and $\gamma$ have
separate interpretations in the overarching Bayesian picture, the
first as a scaling of the prior precision and the second as observational
noise variance.  In this paper
we deepen the connection between the Bayesian methodology and classical
methods.

The recent paper \cite{ILS2013} opens up the prospect for a
statistical explanation of iterative regularization methods in the form of
\begin{equation*}
\begin{aligned}
u_{n}&=u_{n-1}+K_{n}(y-Au_{n-1})
\end{aligned}
\end{equation*}
with a general Kalman gain operator $K_{n}$.
In this paper, we establish the connection between iterative regularization methods (c.f. \cite{EHN1996,HG1998}) and filter based methods \cite{LSZ2015} with respect to an artificial dynamic system. More precisely, for a linear inverse problem, we verify that the iterative Tikhonov regularization method
\begin{equation}
u_{n}=u_{n-1}+A^*(AA^*+\alpha I)^{-1}(y-Au_{n-1})
\label{eq:sit}
\end{equation}
is closely related to filtering methods such as 3DVAR and the
Kalman filter when applied to the partially observed linear dynamical
system \ref{eq:dynamic}. The similarity between both schemes provides a
probabilistic interpretation of iterative regularization methods, and
allows the possibility of quantifying uncertainty via the variance.
On the other hand, we will employ techniques from the convergence
analysis arising in regularization theories \cite{EHN1996}
to shed light on the convergence of filter based methods,
especially when the linear observation operator is ill-posed.

The paper is organized as follows. We first introduce filter
based methods for the artificial dynamics \ref{eq:dynamic}
in Section \ref{se_ArtiDyna}. Section \ref{se_Assumptions} describes
some general useful formulae which are relevant to all the filters
we study, and lists our  main assumptions on the inverse problem
of interest. In Sections \ref{se_Kalman} and \ref{se_3DVAR}
respectively, detailed asymptotic analyses are given for
the Kalman filter method and 3DVAR, for both data models.
The final Section \ref{se_num} presents numerical illustrations
confirming the theoretical predictions.

\section{Filters For The Artificial Dynamics}\label{se_ArtiDyna}

\subsection{Filter Definitions}

Recall the artificial dynamics \eqref{eq:dynamic},
where the observation operator $A$ also defines the
inverse problem \eqref{eq:base}, and $\{\eta_n\}_{n \ge 0}$ is
an i.i.d. sequence with $\eta_1 \sim \N(0,\gamma^2I)$.
The aim of filters is to estimate $u_n$ given the data $\{y_j\}_{j=1}^n.$
In particular, probabilistic filtering aims to estimate the probability
distribution of the conditional random variable $u_n|\{y_j\}_{j=1}^n.$

If we assume that $u_0 \sim {\mathcal N}(m_0,C_0)$ then
the desired conditional random variable is Gaussian, because of the
linearity inherent in \eqref{eq:dynamic}, together with the assumed
Gaussian structure of the noise sequence $\{\eta_n\}_{n \ge 0}$.
Furthermore the independence of the elements of the noise sequence
means that the Gaussian can be updated sequentially in a Markovian
fashion. If we denote by $m_n$ the mean, and by $C_n$ the covariance,
then we obtain the Kalman filter updates for these two
quantities:
\begin{subequations}
\begin{eqnarray}
 \label{eq:Kgain}K_{n} &=&\tri{C}{n-1}{}A^*\left(A\tri{C}{n-1}{}A^*+\gamma^2I\right)^{-1}\\
 \label{eq:mean}\tri{m}{n}{}&=&\tri{m}{n-1}{}+K_{n}(y_{n}-A\tri{m}{n-1}{})  \\
 \label{eq:covariance}\tri{C}{n}{}&=& (I-K_{n}A)\tri{C}{n-1}{}.
\end{eqnarray}
\end{subequations}
The operator $K_n$ is known as the {\em Kalman gain} and
the inverse of the covariance, the {\em precision operator} $C_{n}^{-1}$,
may be shown to satisfy
\begin{eqnarray}
\label{eq:preci}
  C_{n}^{-1} = C_{n-1}^{-1} + \frac{1}{\gamma^2}A^*A.
\end{eqnarray}
All of these facts concerning the Kalman filter may be found in
Chapter 4 of \cite{LSZ2015}.
Expression \eqref{eq:preci} requires careful justification in
infinite dimensions, and this is provided in \cite{ALS2013}
in certain settings. However we will only use \eqref{eq:preci} as
a quick method for deriving useful formulae, not expressed in terms of
precision operators, which can be justified directly under the
assumptions we make.

A simplification of the Kalman filter method is the
3DVAR algorithm \cite{LSZ2015} which is not, strictly speaking,
a probabilistic filter because it
does not attempt to accurately track covariance information. Instead
the covariance is fixed in time at
\begin{equation}
C_{n-1}=\frac{\gamma^2}{\alpha}\Sigma_0
\label{eq:3DVARcovariance2}
\end{equation}
for some fixed positive
and self-adjoint operator $\Sigma_0.$
The parameter $\alpha$ is a scaling constant the inverse of which measures
the relative size of the fixed covariance of the filter relative
to that of the data.  Imposing this simplification on
equations \eqref{eq:Kgain}, \eqref{eq:mean} gives
\begin{subequations}
\label{eq:3DVARN}
\begin{eqnarray}
\label{eq:3DVARKgain} K_{n}& \equiv&  \mathcal{K} :=\Sigma_{0}A^*\left(A\Sigma_0 A^*+\alpha I\right)^{-1} \\
\label{eq:3DVARmean} \tri{\zeta}{n}{}&=&\tri{\zeta}{n-1}{}+\mathcal{K}(y_n-A\tri{\zeta}{n-1}{}).
\end{eqnarray}
\end{subequations}
It is also helpful to define, from \eqref{eq:covariance},
\begin{equation}
\label{eq:3DVARcovariance} \mathcal{C}\equiv \frac{\gamma^2}{\alpha}(I-\mathcal{K}A)\Sigma_0.
\end{equation}
Notice \cite{EHN1996,LP2013} that
the iteration \eqref{eq:3DVARmean} looks like a stationary
iterative Tikhonov method \eqref{eq:sit}
with $A$ replaced by $A\Sigma_0^{\frac12}.$

Throughout the paper $(K_n,m_n,C_n)$ stands for Kalman gain, updated mean
and updated covariance for the Kalman filter method and
$(\mathcal{K}, \zeta_n, \mathcal{C})$ is the related sequence
of quantities for 3DVAR.

\subsection{Asymptotic Behaviour of Filters}

We will view the filters as methods for reconstructing the truth
$\true$; in particular we will study the proximity of $m_n$
(for the Kalman filter) and $\zeta_n$ (for 3DVAR) to $\true$ for
various large $n$ asymptotics. Although the assumption in the
{\em derivation} of the filters is that $y_n$ is an i.i.d. sequence
of the form ${\mathcal N}(A\true,\gamma^2 I)$, we will not
always assume that the data available is of this form; to be
precise Data Model $1$ is compatible with this assumption whilst
Data Model $2$ is not.

Recall that Data Model $1$ refers to the
situation where the data used in the Kalman and 3DVAR filters
has the form  $y_n = Au^{\dag} + \eta_n$, where the $\eta_n$ are
i.i.d. ${\mathcal N}(0,\gamma^2 I).$ Given such a data sequence
we can generate an auxiliary element
  \begin{align*}
  \bar{y} =\frac{1}{n} \sum_{j=1}^n y_j  = A u^{\dag} + \frac{1}{n}\sum_{j=1}^n \eta_j
  \end{align*}
  with $\bar{\eta} = \frac{1}{n}\sum_{j=1}^n \eta_j$ and $\bar{\eta} \sim\mathcal{N}(0, \frac{\gamma^2}{\sqrt{n}}I)$. The law of large numbers and
central limit theorem thus
allows us to consider an inverse problem of the form \eqref{eq:base}
with noise level reduced by a factor of $\sqrt{n}.$ We will study,
in the sequel, whether the Kalman or 3DVAR filters are able
to automatically exploit the
decreased uncertainty inherent in an i.i.d. data set of this form.

For Data Model $2$ we simply assume that the data used in the
filters is of the form $y_n \equiv y$ where $y$ is given by
\eqref{eq:base} with $\eta \sim {\mathcal N}(0,\gamma^2 I).$
From the discussion in the preceding paragraph, we clearly
expect less accurate reconstruction in this case.
For this data model we may view 3DVAR as a stationary iterative
Tikhonov regularization, whilst the Kalman filter is an alternative
iterative non-stationary regularization scheme, since $K_n$
is updated in each step. In addition, the statistical perspective
not only allows us to obtain an estimator (the mean (\ref{eq:mean}) or (\ref{eq:3DVARmean})), but also in the case of the Kalman filter method,
to quantify the uncertainty (via the covariance (\ref{eq:covariance})).
This uncertainty quantification perspective provides additional
motivation for the filtering approaches considered herein.

In this paper our primary focus is the asymptotic large $n$
behavior of the Kalman filter method and 3DVAR. More precisely, we
are interested in the accurate recovery of the true state
$u^{\dag}$ when the noise variance vanishes, i.e.
$\gamma^2 \rightarrow 0$ for Data Models $1$ and $2$,
or as $n\rightarrow \infty$ for Data Model $1$ (by the law of
large numbers/central limit theorem discussion above).

To highlight the difficulties inherent in ill-posed inverse
problems in this regard, we note the following which is
a straightforward consequence of Theorem 4.10 in \cite{LSZ2015}
when specialized to linear problems. Here, and in what follows,
$\|\cdot\|$ denotes both the norm on $X$ and the operator
norm from $X$ into itself.

\begin{proposition} Consider the 3DVAR algorithm with $\{y_n\}_{n \ge 1}$.
Assume that there exists a constant $L$ such that
$\|I-\mathcal{K}A\|\leq L<1$ and that $\|\mathcal{K}\| < \infty.$
Then, for Data Model $2$, it yields
$$\limsup_{n\to \infty}\|\zeta_n-\true\|\leq  \frac{\|\mathcal{K}\|\|\eta\|}{1-L}.$$
\end{proposition}

\noindent
Note however, that if the observation operator $A$ is compact or the
inversion is ill-posed, the assumption $L<1$ in the preceding proposition
cannot hold. More precisely, the operator  $I-\mathcal{K}A$ is no longer
contractive since the spectrum of the operator $\mathcal{K}A$ clusters
at the origin. Our focus in the remainder of the paper will be on such
ill-posed inverse problems.

\section{Main Assumptions and General Properties of Filters}
\label{se_Assumptions}

\subsection{Assumptions}

Recall that $\|\cdot\|$ denotes both the norm on $X$ and the
operator norm from $X$ into itself. Throughout $C$ will denote
a generic constant, independent of the key limiting quantities
$\gamma$ and $n$. Our main assumption is:

\begin{assumption}\label{assp_main}
For both the Kalman filter and the 3DVAR filter, we assume
\begin{description}
  \item[(i)] the variance $C_0=\frac{\gamma^2}{\alpha}\Sigma_0$
and $\mathcal{R}(\Sigma_0^{1/2}) \subset \D{A}$, where
$\alpha$ is a positive constant and $\Sigma_0$ is positive
self-adjoint, and $\Sigma_0^{-1}$ is a densely defined
unbounded self-adjoint strictly positive operator;
  \item[(ii)] the forward operator $A$ satisfies
   \begin{equation}
   \label{eq:CA}
   C^{-1} \|\Sigma_{0}^{\frac a 2}x\|\leq \|Ax\| \leq  C\|\Sigma_{0}^{\frac a 2}x\|
   \end{equation}
   on $X$ for some constants $a>0$ and $1 \leq C < \infty$;
  \item[(iii)] the initial mean satisfies $m_0-\true\in \D{\Sigma_{0}^{-\frac s 2}}$ (or $\zeta_0-\true\in \D{\Sigma_{0}^{-\frac s 2}}$) with $0 \leq s \leq a+2$;
  \item[(iv)] the operator $\Sigma_0$ in item (i) is trace class on $X$.
\end{description}
\end{assumption}
We briefly comment on these items. Item $(i)$ allows a well defined operator \begin{equation}
\label{eq:B0}
B_0:= A\Sigma_0^{1/2}
\end{equation}
 which is essential in carrying out our analysis.
Item $(ii)$ is often called the {\em link condition} and it connects both
operators $A$ and $\Sigma_0$ (or $C_0$). The third item $(iii)$ is the
{\em source condition} (regularity) of the true solution \cite{EHN1996}.
The final item $(iv)$ makes $C_0$ a well-defined covariance
operator on $X$ \cite{bog98}.

Item $(ii)$ in the preceding assumption
is automatically satisfied if $A^*A$ and $\Sigma_0$ have the same
eigenfunctions and certain decaying singular values.
Item $(iii)$ can then be expressed in this eigenbasis.
When studying the Kalman filter we will, in some instances, employ
the following specific form of items $(ii)$, $(iii)$.
Comparison of Assumptions \ref{assp_main} and \ref{assp_main2} we
see that they are identical if $a=\frac{2p}{1+2\epsilon}$ and
$s=\frac{2\beta}{1+2\epsilon}.$

\begin{assumption}\label{assp_main2}
Let the variance $C_0 = \frac{\gamma^2}{\alpha} \Sigma_0$.
The operators $\Sigma_0$ and $A^*A$ have the same eigenfunctions
$\{e_i\}$ with their eigenvalues $\{\lambda_i\}$ and $\{\kappa_i^2\}$
satisfying
   \begin{equation*}
    \lambda_i = i^{-1-2\epsilon}, \qquad C^{-1} i^{-p} \leq \kappa_i \leq C i^{-p}
   \end{equation*}
   for some $\epsilon>0$, $p>0$ and $C\geq 1$. Furthermore, by choosing the initial mean $m_0=0$, the true solution $\true$ with its coordinates $\{u^{\dag,i}\}$ in the basis $\{e_i\}$ obeys $\sum_{i=1}^{\infty} (u^{\dag,i})^2 i^{2\beta} < \infty$.
\end{assumption}

\subsection{Filter Properties}

We start by deriving properties of the Kalman filter method under
Data Model $1$. Other cases can be simply derived from minor variants
of this setting.  Recall from (\ref{eq:mean})
$$m_n=(I-K_nA)m_{n-1}+K_ny_n$$
and note that
$$\true=(I-K_nA)\true+K_nA\true.$$
Under Data Model $1$ we have $y_n=A\true+\eta_n$ and hence
the total error $e_n: = m_n - u^{\dag}$ satisfies
\begin{align}
e_n & =(I-K_n A) e_{n-1} + K_n \eta_n  \nonumber \\
& = \prod_{j=1}^n (I-K_j A) e_0 + \sum_{j=1}^{n-1} \left(\prod_{i=n-j}^{n-1}(I-K_{i+1}A)\right) K_{n-j} \eta_{n-j} + K_n \eta_n \label{eq:generaldecomposition}\\
& := J_1+J_2. \nonumber
\end{align}
Here
\begin{align*}
J_1 & = \prod_{j=1}^n (I-K_j A) e_0\quad{\rm {and}}\\
J_2 & = \sum_{j=1}^{n-1} \left(\prod_{i=n-j}^{n-1}(I-K_{i+1}A)\right) K_{n-j} \eta_{n-j} + K_n \eta_n.
\end{align*}

To establish a rigorous convergence analysis, the mean squared error (MSE)
$\mathbb{E}\|m_n - u^{\dag}\|^2$ is of particular interest.
Since $u^{\dag}$ is deterministic and each $\eta_n$ is i.i.d Gaussian
we obtain a bias-variance decomposition of the MSE:
\begin{align}\label{eq:bivadecomp}
\mathbb{E}\|m_n - u^{\dag}\|^2 =  \|J_1\|^2 + \mathbb{E}\|J_2\|^2.
\end{align}
To proceed further, both terms $J_1$ and $J_2$ need to
be calibrated more carefully.

We consider the operator $I-K_n A$ which appears in both terms.
By (\ref{eq:covariance}), we obtain
\begin{equation*}
  C_n=(I-K_nA)C_{n-1} = \prod_{j=1}^n(I-K_jA)C_0,
\end{equation*}
which is equivalent to
\begin{equation}
\label{eq:prod}
  \prod_{j=1}^n(I-K_jA) = C_nC_0^{-1}.
\end{equation}
Notice that (\ref{eq:preci}) yields
\begin{eqnarray}
\label{eq:recur_preci}
   C_n^{-1} = C_{n-1}^{-1} + \frac{1}{\gamma^2}A^*A = C_0^{-1} + \frac{n}{\gamma^2} A^*A.
\end{eqnarray}
By (\ref{eq:prod}) and (\ref{eq:recur_preci}) we obtain
\begin{align}
  \prod_{j=1}^{n} (I-K_j A) & =  C_nC_0^{-1} = (C_0^{-1} + n A^*A/\gamma^2)^{-1} C_0^{-1}\nonumber\\
  &  =  C_0^{\frac1 2}\gamma^2(\gamma^2 I+ n C_0^{\frac 1 2}A^*AC_0^{\frac 1 2})^{-1} C_0^{-\frac 1 2} \label{eq:kgexplicit}.
\end{align}
We will use this expression (which is well-defined in view of
Assumption \ref{assp_main} $(i)$)
and the labelled equations preceding it
in this subsection, frequently in what follows.

\section{Asymptotic Analysis of the Kalman Filter} \label{se_Kalman}

In this section we investigate the asymptotic behaviour of the Kalman
filter (\ref{eq:Kgain})-(\ref{eq:covariance}), under
Assumption \ref{assp_main}. In particular, we are interested in
whether we can reproduce the minimax convergence rate.
This minimax rate is achieved by adopting
Assumption \ref{assp_main} in the diagonal
form of Assumption \ref{assp_main2}.

\subsection{Kalman Filter and Data Model $1$}\label{subse_KFDM2}
We present the main results in current subsection.
\begin{thm}\label{thm_KFDM1}
  Let Assumption \ref{assp_main} hold. Then the Kalman filter method (\ref{eq:Kgain})-(\ref{eq:covariance}) yields a bias-variance decomposition of the MSE
\begin{align*}
\E\|m_n - \true\|^2 \leq C \left(\frac{\alpha}{n}\right)^{\frac{s}{a+1}} + \frac{\gamma^2}{\alpha} \tr(\Sigma_0)
\end{align*}
for the Data Model $1$. Setting $\alpha=N^{\frac{s}{s+a+1}}$ and
stopping the iteration when $n=N$ then gives
\begin{eqnarray}\label{rate1_KF}
\E\|m_N - \true\|^2 \leq \left(C+\gamma^2 \tr(\Sigma_0)\right) N^{-\frac{s}{s+a+1}}.
\end{eqnarray}
\end{thm}
\begin{thm}\label{thm_KFDM1Ass2} Let Assumption \ref{assp_main2} hold. Then the Kalman filter method (\ref{eq:Kgain})-(\ref{eq:covariance}) yields a bias-variance decomposition of the MSE
\begin{align*}
\E\|m_n - \true\|^2 \leq C \left(\frac{\alpha}{n}\right)^{\frac{2\beta}{1+2\epsilon +2 p}} + \gamma^2 n^{-\frac{2\epsilon}{1+2\epsilon+2p}} \alpha^{-\frac{1+2p}{1+2\epsilon+2p}}
\end{align*}
for the Data Model $1$. Setting
$\alpha=N^{\frac{2(\beta-\epsilon)}{1+2\beta+2p}}$ and stopping
the iteration when $n=N$ then gives the following minimax
convergence rate:
\begin{align*}
\E\|m_N - \true\|^2 \leq  C N^{-\frac{2\beta}{1+2\beta+2p}}.
\end{align*}
\end{thm}

\begin{remark}

\begin{itemize}

\item (i)
Under the Assumption \ref{assp_main2} we prove unconditional convergence
of the Kalman filter method for any fixed $\alpha>0$ and $\gamma>0$,
noticing that both the bias and variance vanish when $n$ goes to infinity.
The key ingredient which leads to this unconditional convergence, in
comparison with Assumption \ref{assp_main}, is that the rate of
decay of the eigenvalues of the variance operator $\Sigma_0$ is made
explicit under Assumption \ref{assp_main2}; this is to be contrasted
with the weaker assumption $\tr(\Sigma_0)< \infty$ made in
Assumption \ref{assp_main} $(iv)$.

\item (ii)
By choosing $\alpha$ depending on the update step $N$, again with
fixed $\gamma$, both Theorems \ref{thm_KFDM1} and \ref{thm_KFDM1Ass2}
yield convergence rates in the MSE sense.
Indeed in the second case, where we use Assumption \ref{assp_main2},
the minimax rate of $N^{-\frac{2\beta}{1+2\beta+2p}}$ is achieved.
This minimax rate may also be achieved from the Bayesian posterior
distribution with appropriate tuning of the prior in terms of the
(effective) noise size $\sqrt{N}$ \cite{KVZ2011}; the tuning of
the prior is identical to the tuning of the initial condition for
the covariance $C_0$, via choice of $\alpha.$ \hfill$\diamondsuit$
\end{itemize}
\end{remark}

Proof of Theorems \ref{thm_KFDM1} and \ref{thm_KFDM1Ass2} is
straightforward by means of a bias-variance decomposition.
Let Assumption \ref{assp_main} $(i)$ hold, noting that then
$B_0:= A\Sigma_0^{1/2}$ is well-defined. We thus obtain,
by (\ref{eq:kgexplicit}),
\begin{eqnarray}
\label{eq:Rapproxi_km2}
      \prod_{i=1}^{n} (I-K_i A) =\Sigma_0^{\frac1 2}\alpha(\alpha I + n B_0^*B_0)^{-1}\Sigma_0^{-\frac1 2} =\Sigma_0^{\frac1 2}r_{1,\frac{\alpha}{n}}(B_0^*B_0) \Sigma_0^{-\frac 1 2},
\end{eqnarray}
where
\begin{equation}\label{eq_rn}
  r_{1,\frac{\alpha}{n}}(\lambda) : = \frac{\frac{\alpha}{n}}{\frac{\alpha}{n}+\lambda} = \frac{\alpha}{\alpha+n\lambda}.
\end{equation}
The operator-valued function $r_{1,\frac{\alpha}{n}}$ (\ref{eq_rn})
has been verified to be powerful in the convergence analysis of
deterministic regularization schemes -- see \cite[Ch.2]{LP2013}.
In that context the following inequality is useful:
\begin{eqnarray}
  \label{eq:pos_rr}
   |\lambda^t r_{1,\frac{\alpha}{n}}(\lambda)|&\leq&  \left(\frac{\alpha}{n}\right)^{t}, \quad \quad \lambda\in (0,\|B_0^*B_0\|], \quad 0\leq t \leq 1.
\end{eqnarray}

Following these ideas we
obtain the next two lemmas describing
the bias and variance error bounds.
We leave the proofs of both lemmas to the Appendix.
Theorems \ref{thm_KFDM1} and \ref{thm_KFDM1Ass2}
are consequences, by choosing the parameter $\alpha$ appropriately.

\begin{lem}[Bias for Kalman filter]\label{lemma_Kalmanbias}
Let Assumption \ref{assp_main} $(i)$-$(iii)$ hold. Then the Kalman filter method (\ref{eq:Kgain})-(\ref{eq:covariance}) yields
   \begin{align*}
   \|J_1\|^2 \leq C\left(\frac{\alpha}{n}\right)^{\frac{s}{a+1}}.
   \end{align*}
 Furthermore, if Assumption \ref{assp_main2} is valid, the bias obeys
\begin{align*}
\|J_1\|^2 \leq C\left(\frac{\alpha}{n}\right)^{\frac{2\beta}{1+2\epsilon+2p}}.
\end{align*}

\end{lem}

\begin{lem}[Variance for Kalman filter -- Data Model $1$]\label{lemma_Kalmanvariance}
Let Assumption \ref{assp_main} $(i)$, $(iv)$ hold and
$\{\eta_n\}$ in (\ref{eq:dynamic}) be an i.i.d sequence with
$\eta_1 \sim \mathcal{N}(0,\gamma^2 I)$.
Then the Kalman filter method (\ref{eq:Kgain})-(\ref{eq:covariance}) yields
\begin{align*}
\E\|J_2\|^2 \leq \frac{\gamma^2}{\alpha} \tr(\Sigma_0).
\end{align*}
Furthermore, if Assumption \ref{assp_main2} is valid, the variance obeys
\begin{align*}
\E\|J_2\|^2 \leq \gamma^2 n^{-\frac{2\epsilon}{1+2\epsilon+2p}} \alpha^{-\frac{1+2p}{1+2\epsilon+2p}}.
\end{align*}
\end{lem}

\subsection{Kalman Filter and Data Model $2$}\label{subse_KFDM1}

The key difference between Data Model $2$ and Data Model $1$ is
that, in the case $2$, the noises $\eta_n$ appearing in the
expression for the term $J_2$ are identical, rather than
i.i.d. mean zero as in case $1$. This results in a reduced
rate of convergence in case $2$ over case $1$, as seen
in the following two theorems:

\begin{thm}\label{thm_KFDM2}
Let Assumption \ref{assp_main} hold. Then the Kalman filter method (\ref{eq:Kgain})-(\ref{eq:covariance}) yields a bias-variance decomposition of the MSE
\begin{align*}
\E\|m_n - \true\|^2 \leq C \left(\frac{\alpha}{n}\right)^{\frac{s}{a+1}} + \frac{ n \gamma^2}{\alpha} \tr(\Sigma_0)
\end{align*}
for the Data Model $2$. Fix $\alpha=1$ and assume that the noise variance $\gamma^2=N^{-\frac{a+s+1}{a+1}}$. If the iteration is stopped at
$n=N$ then the following convergence rate is valid:
\begin{eqnarray}\label{rate2_KF}
\E\|m_N - \true\|^2 \leq \left(C+ \tr(\Sigma_0)\right) N^{-\frac{s}{a+1}}.
\end{eqnarray}
\end{thm}
\begin{thm}\label{thm_KFDM2Ass2}
Let Assumption \ref{assp_main2} hold. Then the Kalman filter method (\ref{eq:Kgain})-(\ref{eq:covariance}) yields a bias-variance decomposition of the MSE
\begin{align*}
\E\|m_n - \true\|^2 \leq C \left(\frac{\alpha}{n}\right)^{\frac{2\beta}{1+2\epsilon +2 p}} + \gamma^2 \left(\frac{n}{\alpha}\right)^{\frac{1+2p}{1+2\epsilon+2p}}
\end{align*}
for the Data Model $2$.  Fix $\alpha=1$ and assume that the noise variance
$\gamma^2=N^{-\frac{1+2\beta+2p}{1+2\epsilon+2p}}.$ If the iteration is
stopped at $n=N$ then the following convergence rate is valid:
\begin{align*}
\E\|m_N - \true\|^2 \leq  C N^{-\frac{2\beta}{1+2\epsilon+2p}}.
\end{align*}
\end{thm}

Both convergence rates in Theorems \ref{thm_KFDM2} and \ref{thm_KFDM2Ass2}
are of the same order since the variance has been tuned to scale
in the same way as the bias by choosing to stop the iteration
at $N$, depending on $\gamma \ll 1$, appropriately.
In comparison with the convergence rates in Theorems \ref{thm_KFDM1}
and \ref{thm_KFDM1Ass2}, the ones in this section under Data Model $2$
require small noise $\gamma$; those in the preceding subsection do not because
multiple observations, with additive independent noise, are made
of $A\true.$  Proof of the two preceding theorems is straightforward:
the $J_1$ terms is analyzed as in the preceding subsection
and the $J_2$ term must be carefully
analyzed under the assumptions of Data Model $2$. The key new
result is stated in the following lemma, whose proof may be found
in the Appendix.

\begin{lem}[Variance for Kalman filter method -- Data Model $2$]\label{lemma_KalmanvarianceDM2}
Let Assumption \ref{assp_main} hold and each observation $y_{n}\equiv y$ be fixed. Then the Kalman filter method (\ref{eq:Kgain})-(\ref{eq:covariance}) yields
\begin{align*}
\E\|J_2\|^2 \leq \frac{n \gamma^2}{\alpha} \tr(\Sigma_0).
\end{align*}
Furthermore, if Assumption \ref{assp_main2} is valid, the variance obeys
\begin{align*}
\E\|J_2\|^2 \leq \gamma^2 \left(\frac{n}{\alpha}\right)^{\frac{1+2p}{1+2\epsilon+2p}}.
\end{align*}
\end{lem}

\section{Asymptotic Analysis of 3DVAR}\label{se_3DVAR}

\subsection{Classical 3DVAR}

The mean of the 3DVAR algorithm is given by (\ref{eq:3DVARmean}) and
has the form
$$\tri{\zeta}{n}{}=(I-\mathcal{K}A)\tri{\zeta}{n-1}{}+\mathcal{K}y_n.$$
Furthermore
$$\true=(I-\mathcal{K}A)\true+\mathcal{K}A\true.$$
If we define $\varepsilon_n = \zeta_n -\true$ then we obtain, since
$y_n=A\true+\eta_n$,
\begin{align*}
\varepsilon_n & = (I-\mathcal{K}A) \varepsilon_{n-1} + \mathcal{K} \eta_n \\
& = (I-\mathcal{K}A)^n \varepsilon_0 + \sum_{j=0}^{n-1} (I-\mathcal{K}A)^j \mathcal{K} \eta_{n-j}
\end{align*}
with $\varepsilon_0 := \zeta_0 - \true$.
We further derive
\begin{eqnarray*}
(I-\mathcal{K}A)^n = \Sigma_0^{\frac1 2}(\alpha(\alpha I + B_0^*B_0)^{-1})^n\Sigma_0^{-\frac1 2} = \Sigma_0^{\frac1 2}r_{n,\alpha}(B_0^* B_0)\Sigma_0^{-\frac1 2},
\end{eqnarray*}
by inserting the definition of the
Kalman gain (\ref{eq:3DVARKgain}), and
by assuming Assumption \ref{assp_main} $(i)$.
The operator-valued function $r_{n,\alpha} (\cdot)$ is defined
\begin{align*}
r_{n,\alpha}(\lambda):= \left(\frac{\alpha}{\alpha+\lambda}\right)^n.
\end{align*}
Similarly to the analysis of the Kalman filter, we derive
\begin{align}
\varepsilon_n & =(I-\mathcal{K}A)^n \varepsilon_0 + \sum_{j=0}^{n-1} (I-\mathcal{K}A)^j \mathcal{K} \eta_{n-j} \label{eq:3dvarerror}\\
 & = \Sigma_0^{\frac1 2} r_{n,\alpha}(B_0^* B_0) \Sigma_0^{-\frac{1}{2}} \varepsilon_0 + \sum_{j=0}^{n-1} (I-\mathcal{K}A)^j \mathcal{K} \eta_{n-j} \nonumber \\
& : = I_1 + I_2. \nonumber
\end{align}
Thus, the MSE takes the bias-variance decomposition form
\begin{align*}
\E\|\zeta_n-\true\|^2 = \E\|\varepsilon_n\|^2 = \|I_1\|^2 + \E\|I_2\|^2.
\end{align*}
This leads to the following two theorems:

\begin{thm}\label{thm_3DVARDM1}
Let Assumption \ref{assp_main} hold. Then 3DVAR filter
(\ref{eq:3DVARKgain})-(\ref{eq:3DVARmean}) yields a bias-variance decomposition of the MSE
\begin{align*}
\E\|\zeta_n - \true\|^2 \leq C \left(\frac{\alpha}{n}\right)^{\frac{s}{a+1}} + C \frac{\gamma^2 \ln n}{\alpha} \tr(\Sigma_0)
\end{align*}
for the Data Model $1$. Setting $\alpha=N^{\frac{s}{s+a+1}}$ and stopping
the iteration when $n=N$ then gives
\begin{eqnarray}\label{rate1_3DVAR}
\E\|\zeta_n - \true\|^2 \leq C\left(1+\gamma^2 \tr(\Sigma_0)\right) N^{-\frac{s}{s+a+1}} \ln N.
\end{eqnarray}
 \end{thm}

\begin{thm}\label{thm_3DVARDM1Ass2}
Let Assumption \ref{assp_main2} hold. Then 3DVAR filter
(\ref{eq:3DVARKgain})-(\ref{eq:3DVARmean}) yields a bias-variance decomposition of the MSE
\begin{align*}
\E\|\zeta_n - \true\|^2 \leq C \left(\frac{\alpha}{n}\right)^{\frac{2\beta}{1+2\epsilon+2p}} + C \gamma^2 \alpha^{-\frac{1+2p}{1+2\epsilon+2p}}
\end{align*}
for the Data Model $1$. Setting $\alpha=N^{\frac{2\beta}{1+2\epsilon+2\beta+2p}}$ and stopping
the iteration when $n=N$ then gives
\begin{align*}
\E\|\zeta_N - \true\|^2 \leq C\left(1+\gamma^2\right) N^{-\frac{2\beta}{1+2\epsilon+2\beta+2p}} \ln N.
\end{align*}
\end{thm}

\begin{remark}
The decay rate at the end of the preceding Theorem
is the same as that in Theorem \ref{thm_3DVARDM1}
if $a=\frac{2p}{1+2\epsilon}$ and $s=\frac{2\beta}{1+2\epsilon}$.
This is the setting in which Assumptions \ref{assp_main} and
\ref{assp_main2} are identical.

The preceding two theorems show that, under Data Model $1$
and for any fixed $\alpha$, the (bound
on the) MSE of the 3DVAR filter blows up logarithmically
as $n\rightarrow \infty$ under Assumption \ref{assp_main}, and
is asymptotically bounded for Assumption \ref{assp_main2}.
In contrast, for the Kalman filter method the MSE is asymptotically
bounded or unconditionally converges in $n$ under the same assumptions -- see Theorems \ref{thm_KFDM1} and \ref{thm_KFDM1Ass2}.

With optimal choice of $\alpha$ in terms of the stopping time
of the iteration at $n=N$, comparison of the convergence rates in Theorems \ref{thm_KFDM1}
and \ref{thm_3DVARDM1} (or Theorems \ref{thm_KFDM1Ass2}
and \ref{thm_3DVARDM1Ass2}) shows that the Kalman filter outperforms 3DVAR,
but only by a logarithmic factor (or a H\"{o}lder factor).
For simplicity we only analyze and discuss Data Model $1$ for 3DVAR filter under the additional Assumption \ref{assp_main2}; as for Data Model $2$
one can derive consequences similar to those in the preceding section
in an analogous manner.
\hfill$\diamondsuit$
\end{remark}

We now study Data Model $2$ and 3DVAR. We consider only
Assumption \ref{assp_main}; however the reader may readily
extend the analysis to include Assumption \ref{assp_main2}.
In the case of Data Model $2$, both the Kalman and 3DVAR filters
have the same error bounds:

\begin{thm}\label{thm_3DVARDM2}
Let Assumption \ref{assp_main} hold. Then the 3DVAR
algorithm (\ref{eq:3DVARKgain})-(\ref{eq:3DVARmean}) yields a bias-variance decomposition of the MSE
\begin{align*}
\E\|\zeta_n - \true\|^2 \leq C \left(\frac{\alpha}{n}\right)^{\frac{s}{a+1}} + \frac{n\gamma^2}{\alpha} \tr(\Sigma_0)
\end{align*}
for the Data Model $2$. Fix $\alpha=1$ and assume that the noise variance $\gamma^2=N^{-\frac{a+s+1}{a+1}}$. If the iteration is stopped at $n=N$
then the following convergence rate is valid:
\begin{eqnarray}\label{rate2_3DVAR}
\E\|\zeta_n - \true\|^2 \leq \left(C+ \tr(\Sigma_0)\right) N^{-\frac{s}{a+1}}.
\end{eqnarray}
\end{thm}

The preceding three theorems can be proved by the
bias-variance decomposition and application of the following
three lemmas, whose proofs are left to the Appendix.
However the proof of Theorem \ref{thm_3DVARDM1Ass2}
is not as straightforward as the others and we present
details in the Appendix.

\begin{lem}[Bias for 3DVAR]\label{lemma_3DVARbias}
Let Assumption \ref{assp_main} $(i)$-$(iii)$ hold. Then
3DVAR (\ref{eq:3DVARN}) yields
   \begin{align*}
   \|I_1\|^2 \leq C\left(\frac{\alpha}{n}\right)^{\frac{s}{a+1}}.
   \end{align*}
Furthermore, if Assumption \ref{assp_main2} is valid, the bias obeys
\begin{align*}
\|I_1\|^2 \leq C\left(\frac{\alpha}{n}\right)^{\frac{2\beta}{1+2\epsilon+2p}}.
\end{align*}
\end{lem}

\begin{lem}[Variance for 3DVAR - Data Model $1$]\label{lemma_3DVARvarianceDM1}
Let Assumption \ref{assp_main} $(i)$, $(iv)$ hold and each noise $\eta_n$ in (\ref{eq:dynamic}) i.i.d. generated by $\mathcal{N}(0,\gamma^2 I)$. Then
3DVAR (\ref{eq:3DVARN}) yields
\begin{align*}
\E\|I_2\|^2 \leq C \frac{\gamma^2 \ln n}{\alpha} \tr(\Sigma_0)
\end{align*}
for Data Model $1$. Furthermore, if Assumption \ref{assp_main2} is valid, the variance obeys
\begin{align*}
\E\|I_2\|^2 \leq C \gamma^2 \alpha^{-\frac{1+2p}{1+2\epsilon+2p}}
\end{align*}
and simultaneously
\begin{align*}
\E\|I_2\|^2 \leq C \frac{\gamma^2 \ln n}{\alpha}.
\end{align*}
\end{lem}

\begin{lem}[Variance for 3DVAR - Data Model $2$]\label{lemma_3DVARvarianceDM2}
Let Assumption \ref{assp_main} hold and each observation $y_{n}\equiv y$ be fixed. Then 3DVAR (\ref{eq:3DVARN}) yields
\begin{align*}
\E\|I_2\|^2 \leq \frac{n \gamma^2}{\alpha} \tr(\Sigma_0)
\end{align*}
for Data Model $2$.
\end{lem}

\subsection{Variant of 3DVAR for Data Model $1$}
Recall that the 3DVAR iteration \eqref{eq:3DVARmean}
looks like a stationary
iterative Tikhonov method \eqref{eq:sit} \cite{EHN1996,LP2013},
with $A$ replaced by $A\Sigma_0^{\frac12}$
and with a fixed parameter $\alpha$. The non-stationary iterative
Tikhonov regularization, with varying $\alpha$, has been proven
to be powerful in deterministic inverse problems \cite{HG1998}. We
generalize this method to the 3DVAR setting. Furthermore
we demonstrate that, for Data Model $1$, it is possible to see
blow-up phenomena with this algorithm.

Starting from the classical form of the 3DVAR filter
as given in \eqref{eq:3DVARKgain}, \eqref{eq:3DVARmean},
\eqref{eq:3DVARcovariance}  we propose a variant method  in
which $\alpha$ varies with $n$. We obtain
\begin{subequations}
\begin{eqnarray}
\label{eq:variantKgain} \mathcal{K}_{n}&  := &\Sigma_{0}A^*\left(A\Sigma_0 A^*+\alpha_n I\right)^{-1} \\
\label{eq:variantmean} v_n & := &\tri{v}{n-1}{}+\mathcal{K}_n(y_n-A\tri{v}{n-1}{})  \\
\label{eq:variantRcovariance} \tri{\mathcal{C}}{n}{}& : = & \frac{\gamma^2}{\alpha_n}(I-\mathcal{K}_n A)\Sigma_0.
\end{eqnarray}
\end{subequations}

If we define $\epsilon_n = v_n -\true$ then we obtain,
analogously to the derivation of (\ref{eq:generaldecomposition})
and (\ref{eq:3dvarerror}), the bias-variance decomposition as follows:
\begin{align*}
\epsilon_n & = \left(I-\mathcal{K}_n A \right) \epsilon_{n-1} + \mathcal{K}_n \eta_n \\
& = \prod_{j=1}^n \left(I-\mathcal{K}_j A\right) \epsilon_0 + \sum_{j=1}^{n-1} \left(\prod_{i=n-j}^{n-1} \left(I-\mathcal{K}_{i+1} A\right)\right) \mathcal{K}_{n-j} \eta_{n-j} +\mathcal{K}_n \eta_n \\
& := \mathcal{I}_1 + \mathcal{I}_2
\end{align*}
with
\begin{align*}
\epsilon_0 & = v_0-\true; \\
\mathcal{I}_1 & = \prod_{j=1}^n \left(I-\mathcal{K}_j A\right) \epsilon_0 ; \\
\mathcal{I}_2 & =\sum_{j=1}^{n-1} \left(\prod_{i=n-j}^{n-1} \left(I-\mathcal{K}_{i+1} A\right)\right) \mathcal{K}_{n-j} \eta_{n-j} +\mathcal{K}_n \eta_n.
\end{align*}

By calibrating both terms $\mathcal{I}_1$, $\mathcal{I}_2$ carefully, we obtain the following blow-up result, proved in the Appendix.
Although the result only provides an upper bound, numerical
evidence does indeed show blow-up in this regime.

\begin{thm}\label{thm_3DVARvariant}
Let Assumption \ref{assp_main} hold and let $\alpha_n$ be a
sequence satisfying $\frac{1}{\alpha_n} \leq \tilde{c} \sigma_{n-1}$,
with constant $\tilde{c}$, for $\sigma_n:=\sum_{j=1}^n \frac{1}{\alpha_j}.$
Then the variant EnKF method
(\ref{eq:variantKgain})-(\ref{eq:variantRcovariance}) yields a
bias-variance decomposition of the MSE
\begin{align*}
\E \|v_n - \true\|^2 \leq C \left(\sigma_n^{-\frac{s}{a+1}} + \gamma^2 \tr(\Sigma_0) \sigma_{n}\right)
\end{align*}
for Data Model $1.$ In particular, the geometric sequence
$\alpha_n : = \alpha q^{n-1}$ with $\alpha>0$ and $0<q<1$ yields
\begin{align*}
\E \|v_n - \true\|^2 \leq C \left(q^{\frac{s}{a+1} n} + \gamma^2 \tr(\Sigma_0) q^{-n} \right).
\end{align*}
\end{thm}

\section{Numerical Illustrations}
\label{se_num}

In this section we provide numerical results which display the capabilities of
3DVAR and Kalman filter for solving linear inverse problems of the type described in Section \ref{se_Intro}. In addition, we verify numerically some of the theoretical results from Section \ref{se_Kalman} and Section \ref{se_3DVAR}.

\subsection{Set-Up}\label{se_num:setup}

We consider the two-dimensional domain $\Omega=[0,1]^2$ and define the operators
\begin{eqnarray}\label{num1}
A= (-\triangle)^{-1},\qquad \Sigma_{0}=A^2
\end{eqnarray}
with
\begin{eqnarray*}
D(\triangle)= \Big\{v\in H^{2}(\Omega) ~\vert~~~ \nabla v\cdot \mathbf{n}=0 ~~\textrm{on}~
\partial \Omega, ~~\int_{\Omega}v=0    \Big\}.
\end{eqnarray*}
With this domain $-\triangle$ is positive, self-adjoint and invertible.Note that (\ref{eq:CA}) from Assumption \ref{assp_main} is satisfied with $a=1$ and $C=1$. In the following subsections we produce synthetic data from a true function $u^{\dagger}$ that we generate as a two-dimensional random field drawn from the Gaussian measure $N(0,\Sigma)$ with covariance
\begin{eqnarray}\label{num3}
\Sigma=\Bigl(-\triangle+\frac{1}{10} I\Bigr)^{-(2s+1)}
\end{eqnarray}
with domain of definition $D(\Sigma)=D(\triangle)$ and where $s>0$ is
selected as described below. The shift of $-\triangle$ by a constant introduces
a correlation length into the true function. Furthermore, for simplicity we consider $m_{0}=0$ and note \cite{Stuart2010} that the given selection of $u^{\dagger}$ yields $u^{\dagger}\in  \mathcal{H}^{t}$ for all $0<t<2s$.
Therefore, for discussion of the present experiments we simply assume
that $u^{\dagger}-m_{0}\in \mathcal{D}(\Sigma_{0}^{-s/2})= \mathcal{H}^{2s}$.
Consequently,  in order to satisfy Assumption \ref{assp_main} (iii) we need $s$ such that $0<s\leq a+2=3$. Note that operator $\Sigma_0$ in (\ref{num1}) satisfies Assumption \ref{assp_main} (iv).

The numerical generation of $u^{\dagger}$ is carried out by means of the Karhunen-Loeve decomposition of $u^{\dagger}$ in terms of the eigenfunctions of $\Sigma$ which, from the definition of $D(\Sigma)$, are cosine functions. We recall that for the application of the Kalman Filter and 3DVAR with Data Model 1 we need to generate $N$ instances of synthetic data $\{y_{n}\}_{n=1}^{N}$ where $N$ is the maximum number of iterations of the scheme. Below we discuss the selection of such $N$. The aforementioned synthetic data are generated by means of $y_{n}\equiv  Au^{\dagger}+\eta_{n}$ with $\eta_{n} \sim N(0, \gamma^2 I)$ and $\gamma$ specified below.  For Data Model 2 we produce synthetic data simply by
setting $y_n \equiv y=Au^{\dagger}+\eta$ with $\eta \sim N(0, \gamma^2 I)$.  In order to avoid inverse crimes \cite{KS2005},
all synthetic data used in our experiments are generated on a finer grid ($120\times 120$ cells) than the one (of $60\times 60$ cells) used for the inversion. We use splines to interpolate synthetic data on the coarser grid that we use for the application of the filters.

For both 3DVAR and Kalman filter with Data Model 1, we fix the number of iterations $N$ for the scheme and consider the selection $\alpha=N^{\frac{s}{s+1+a}}$ stated in Theorem \ref{thm_KFDM1} and Theorem \ref{thm_3DVARDM1}. Provided that the filters are stopped according to $n=N$, these theorems ensure the convergence rates in (\ref{rate1_KF}) and (\ref{rate1_3DVAR}) that we verify numerically in the following subsection.  For Data Model 2, Theorem \ref{thm_KFDM2} and Theorem \ref{thm_3DVARDM2} suggest that convergence rates (\ref{rate2_KF}) and (\ref{rate2_3DVAR}) are satisfied for $\alpha=1$ and provided that $\gamma^2=N^{-\frac{a+s+1}{a+1}}$. The latter equality provides an expression for the  number of iterations $N=N(\gamma)$ that we may use as stopping criteria for these filtering algorithms applied to Data Model 2. In Algorithm \ref{Al1} we summarize the Kalman filter and 3DVAR schemes applied to both data models.

\begin{algorithm}{Kalman Filter/3DVAR (Data Model $1/$Data Model $2$)}\label{Al1}~~\\
Let
\begin{eqnarray*}
N=\left\{\begin{array}{cc}
 \textrm{fix integer selected a priori} &  \textrm{   for Data Model $1$}\\
\textrm{round}(\gamma^{-\frac{2(a+1)}{a+s+1}}) & \textrm{for Data Model $2$}\end{array}\right.
\end{eqnarray*}
and
\begin{eqnarray*}
\alpha=\left\{\begin{array}{cc}
N^{\frac{s}{s+1+a}} &  \textrm{for Data Model $1$}\\
1  & \textrm{for Data Model $2$}\end{array}\right.
\end{eqnarray*}
For $n=1,\dots,N$, update $m_{n-1}$ and $C_{n-1}$ as follows
\begin{eqnarray*}
  \tri{m}{n}{}&=&\tri{m}{n-1}{}+K_{n}(y_{n}-A\tri{m}{n-1}{})
\end{eqnarray*}
\begin{eqnarray*}
\tri{C}{n}{}=\left\{\begin{array}{cc}
 (I-K_{n}A)\tri{C}{n-1}{}, &  \textrm{for Kalman Filter}\\
C_{n-1}, & \textrm{for 3DVAR}\end{array}\right.
\end{eqnarray*}
where
\begin{eqnarray*}
K_{n} =\tri{C}{n-1}{}A^*\left(A\tri{C}{n-1}{}A^*+\gamma^2I\right)^{-1}.
\end{eqnarray*}
and where we recall that
\begin{eqnarray*}
y_{n}=\left\{\begin{array}{cc}
Au^{\dagger}+\eta_{n} &  \textrm{for Data Model $1$}\\
y  & \textrm{for Data Model $2$}\end{array}\right.
\end{eqnarray*}
Note that for Data Model $1$ we need at least $N$ independent instances of data.
\end{algorithm}

\subsection{Using Kalman Filter and 3DVAR for solving linear inverse problems}

\begin{figure}[htbp]
\begin{center}
\includegraphics[scale=0.23]{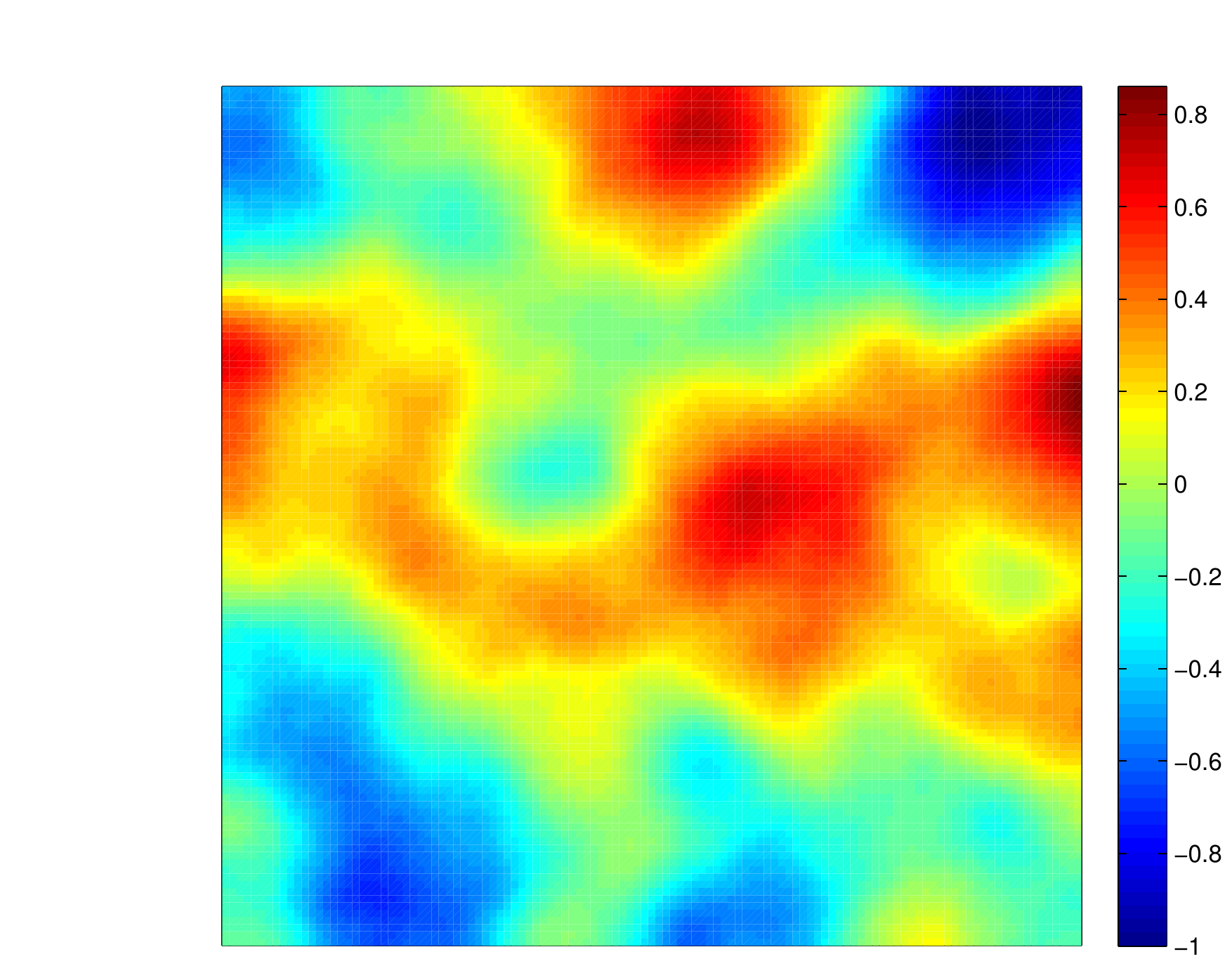}\\

\includegraphics[scale=0.23]{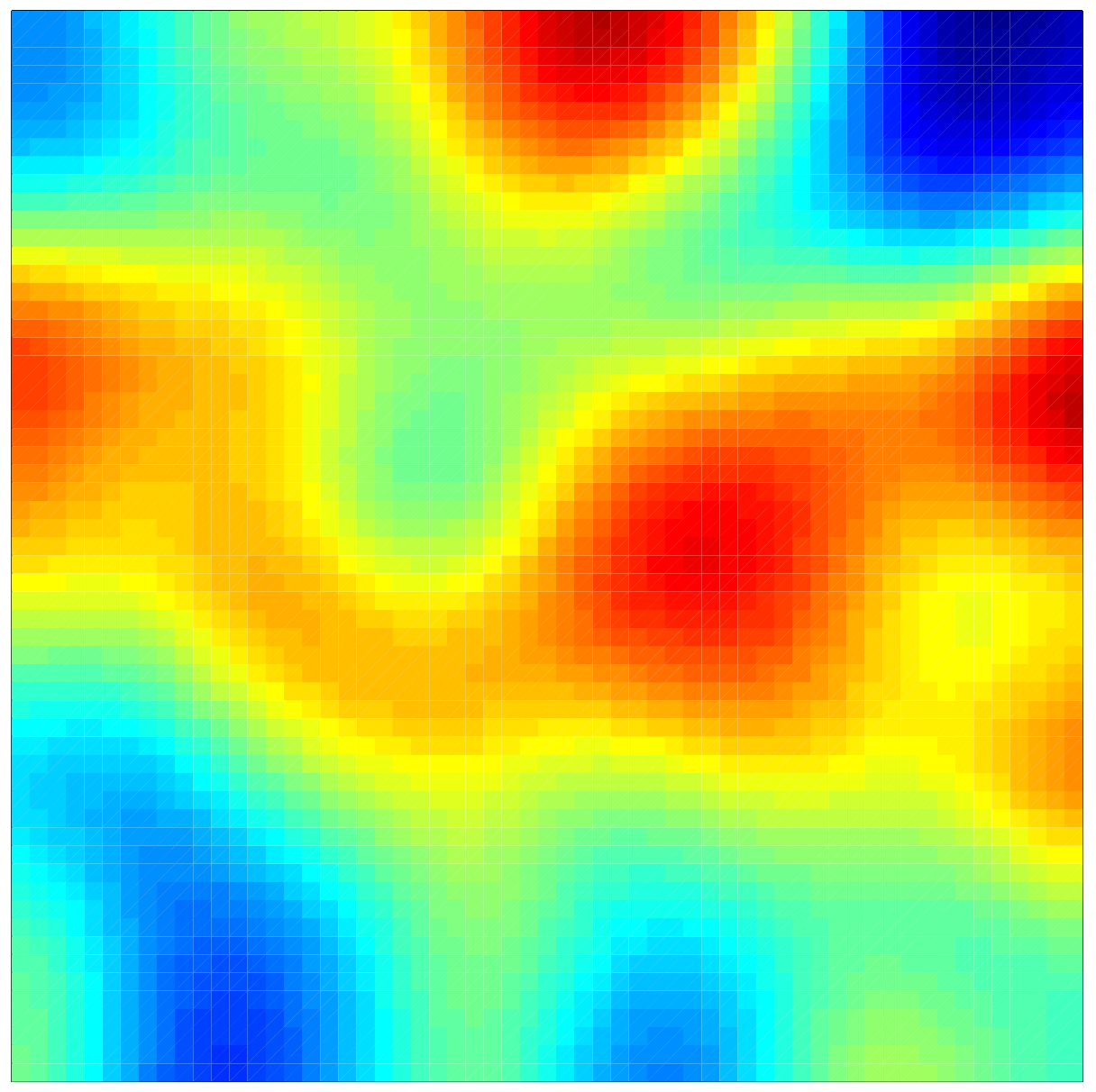}
\includegraphics[scale=0.23]{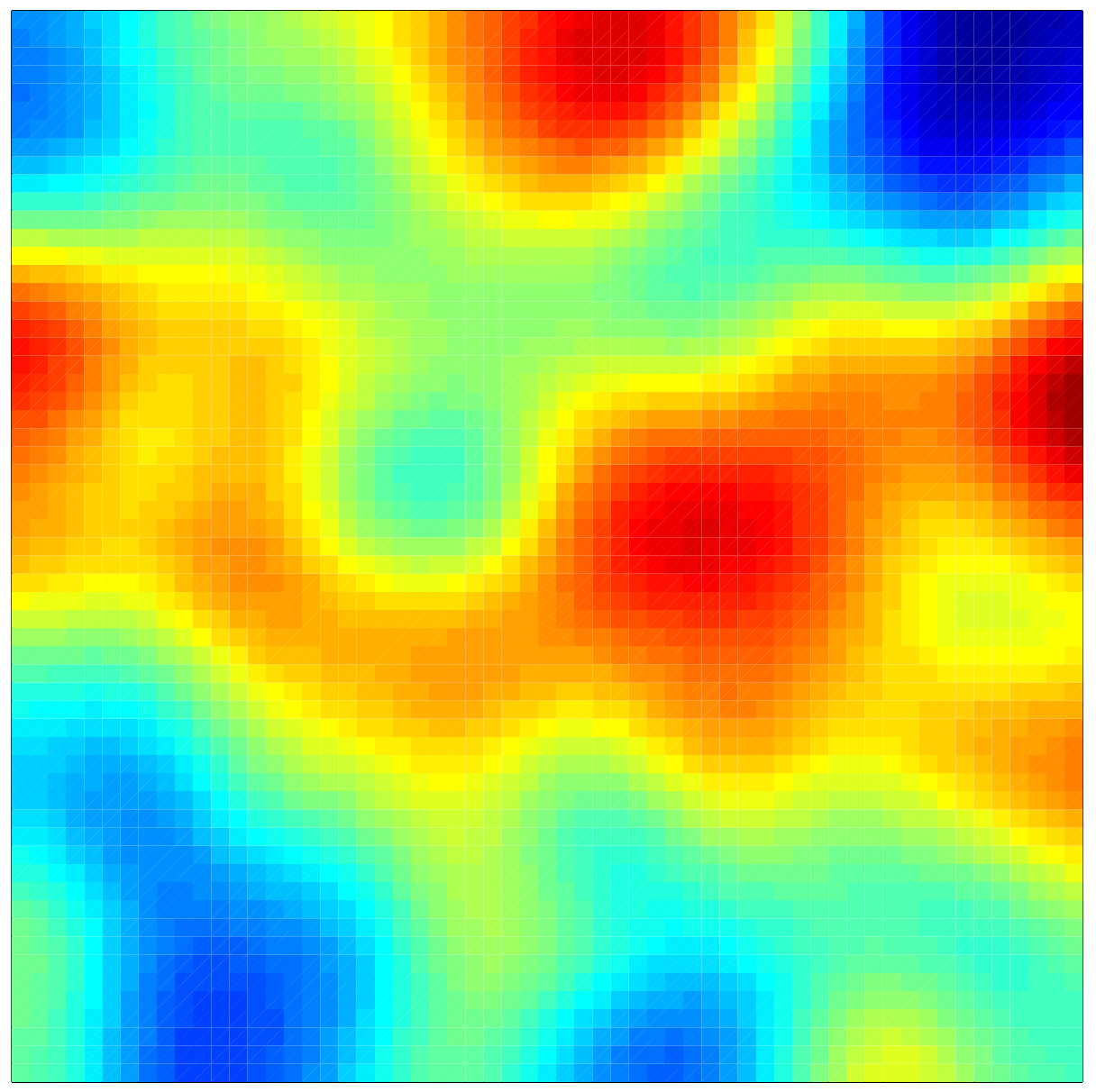}
\includegraphics[scale=0.23]{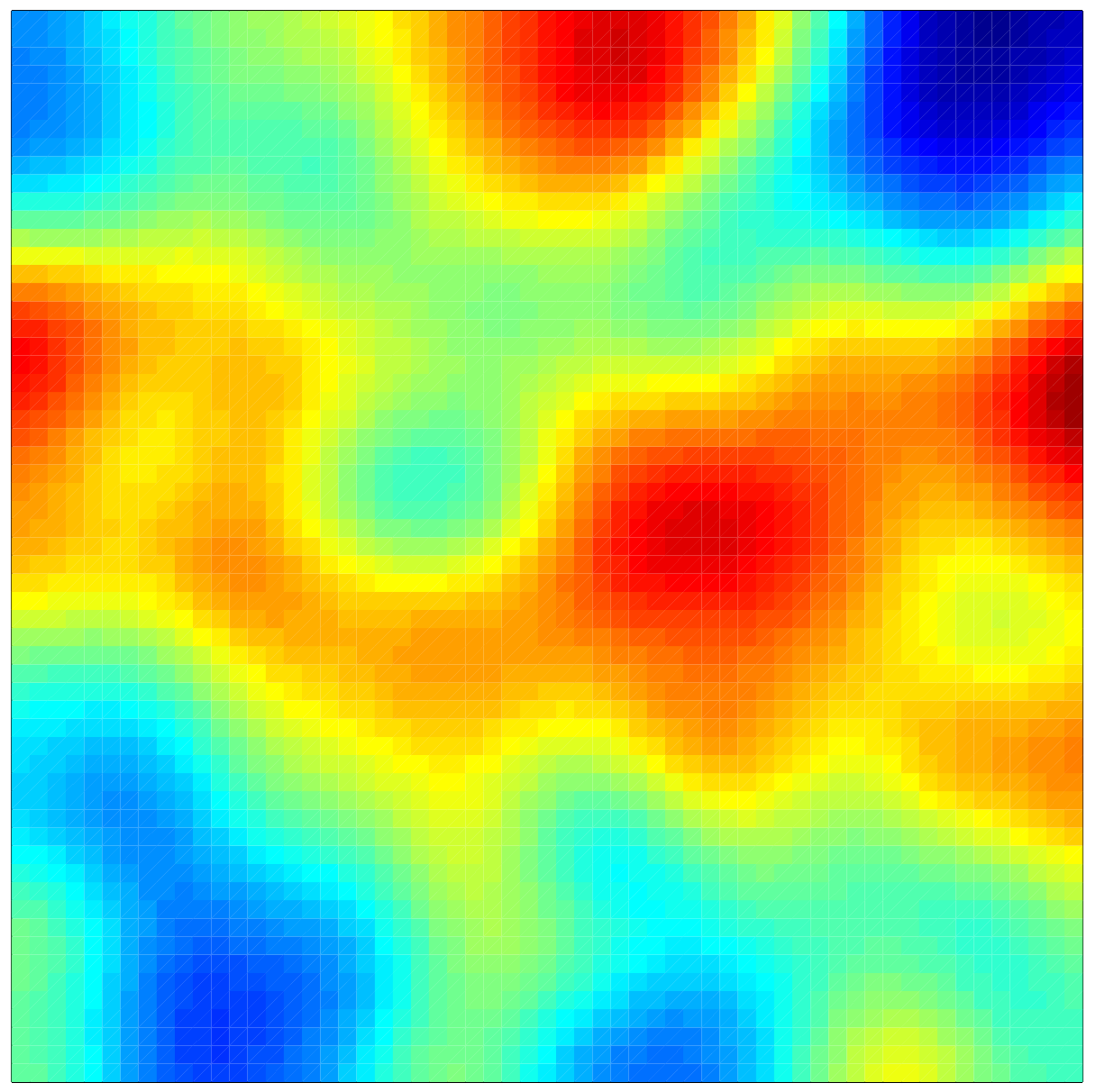}
\includegraphics[scale=0.23]{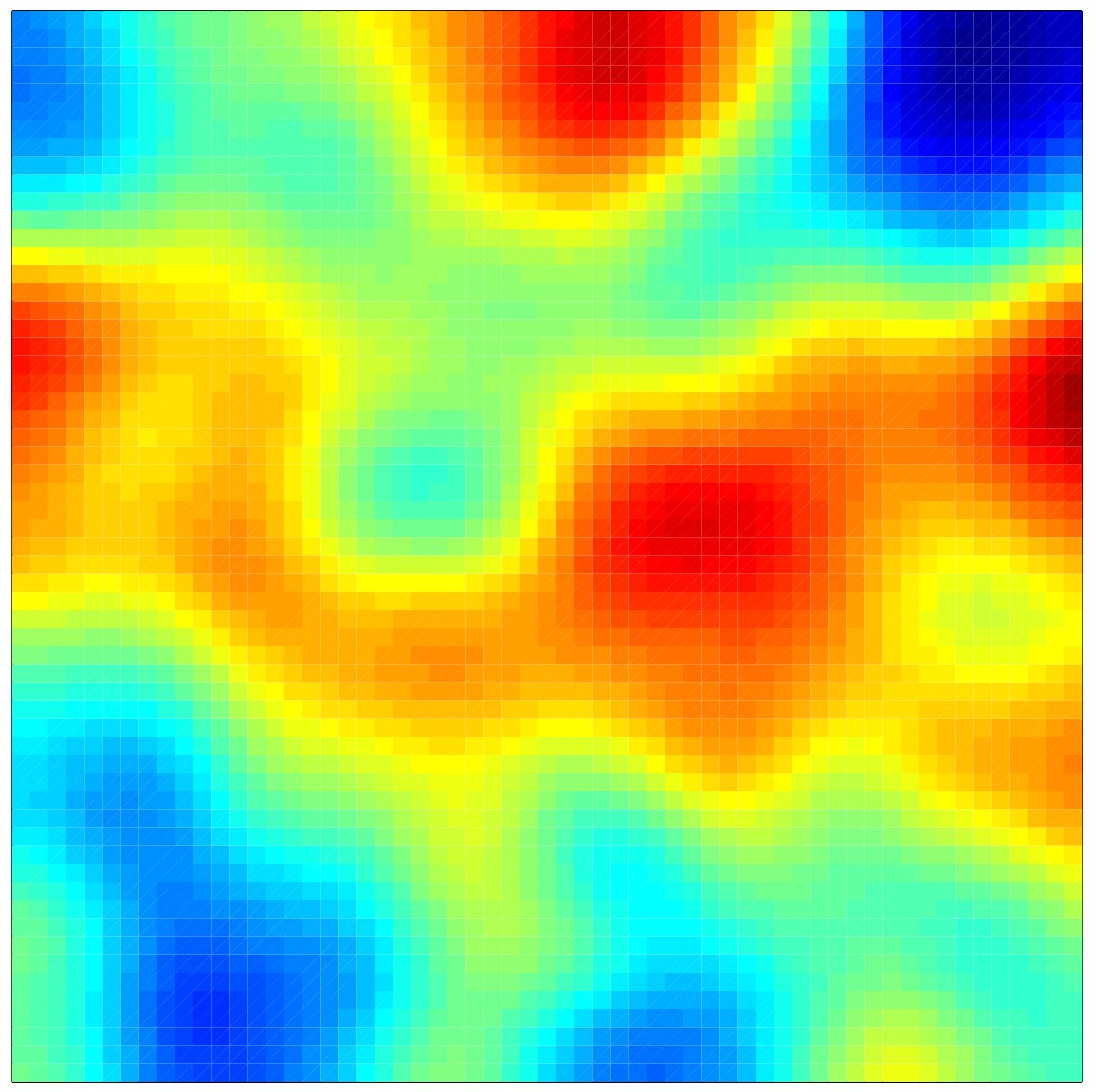}\\

\includegraphics[scale=0.23]{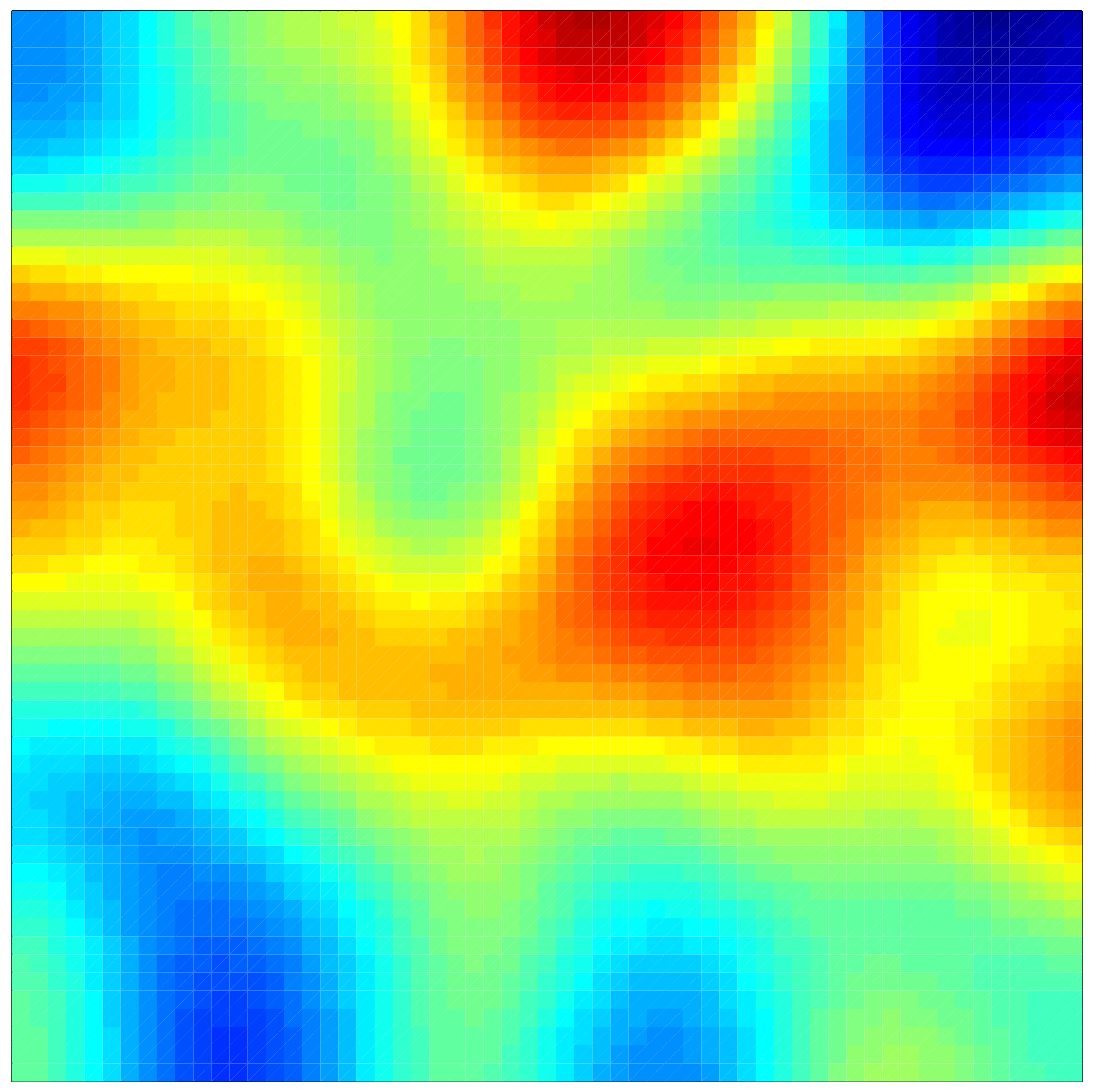}
\includegraphics[scale=0.23]{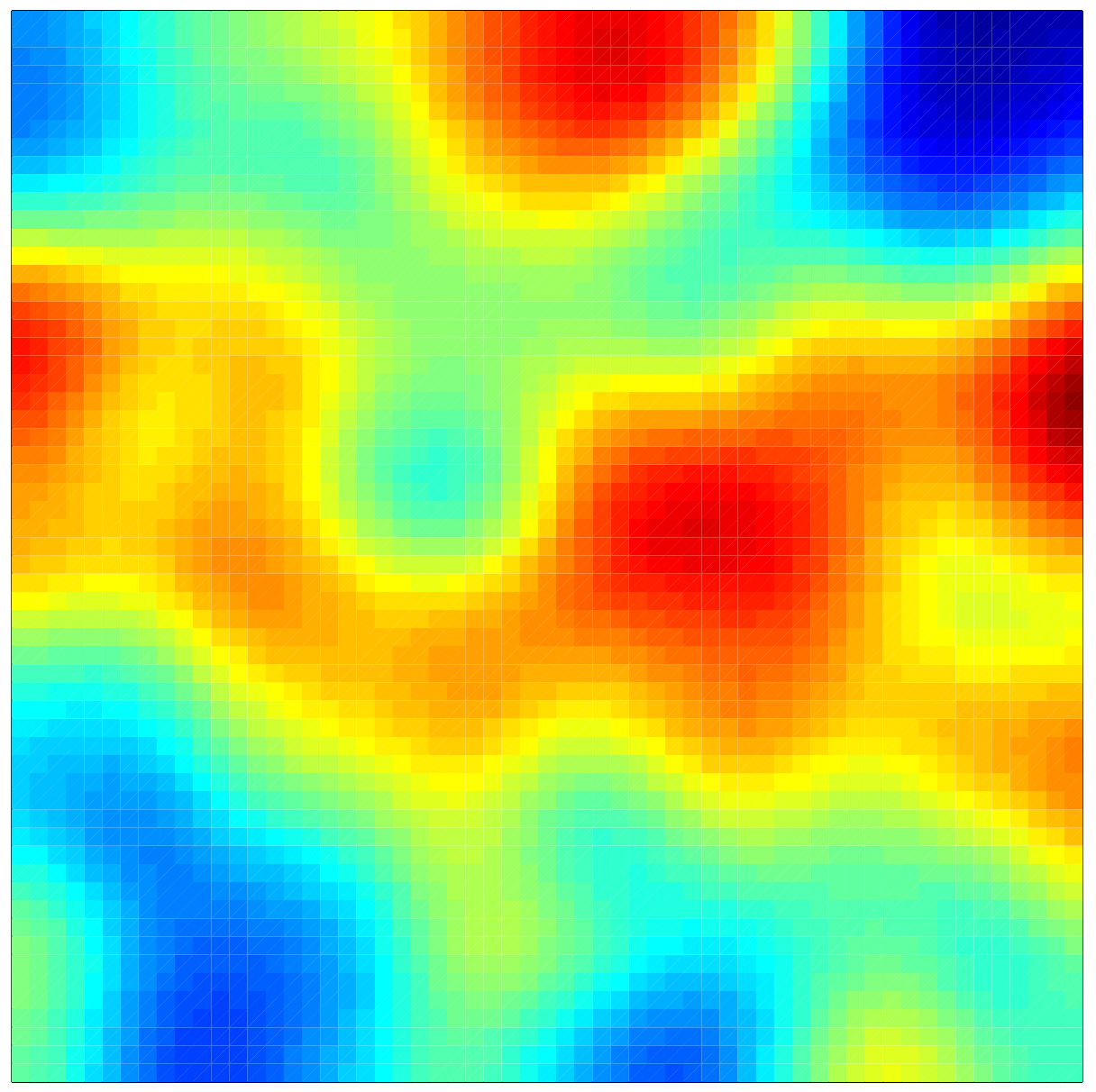}
\includegraphics[scale=0.23]{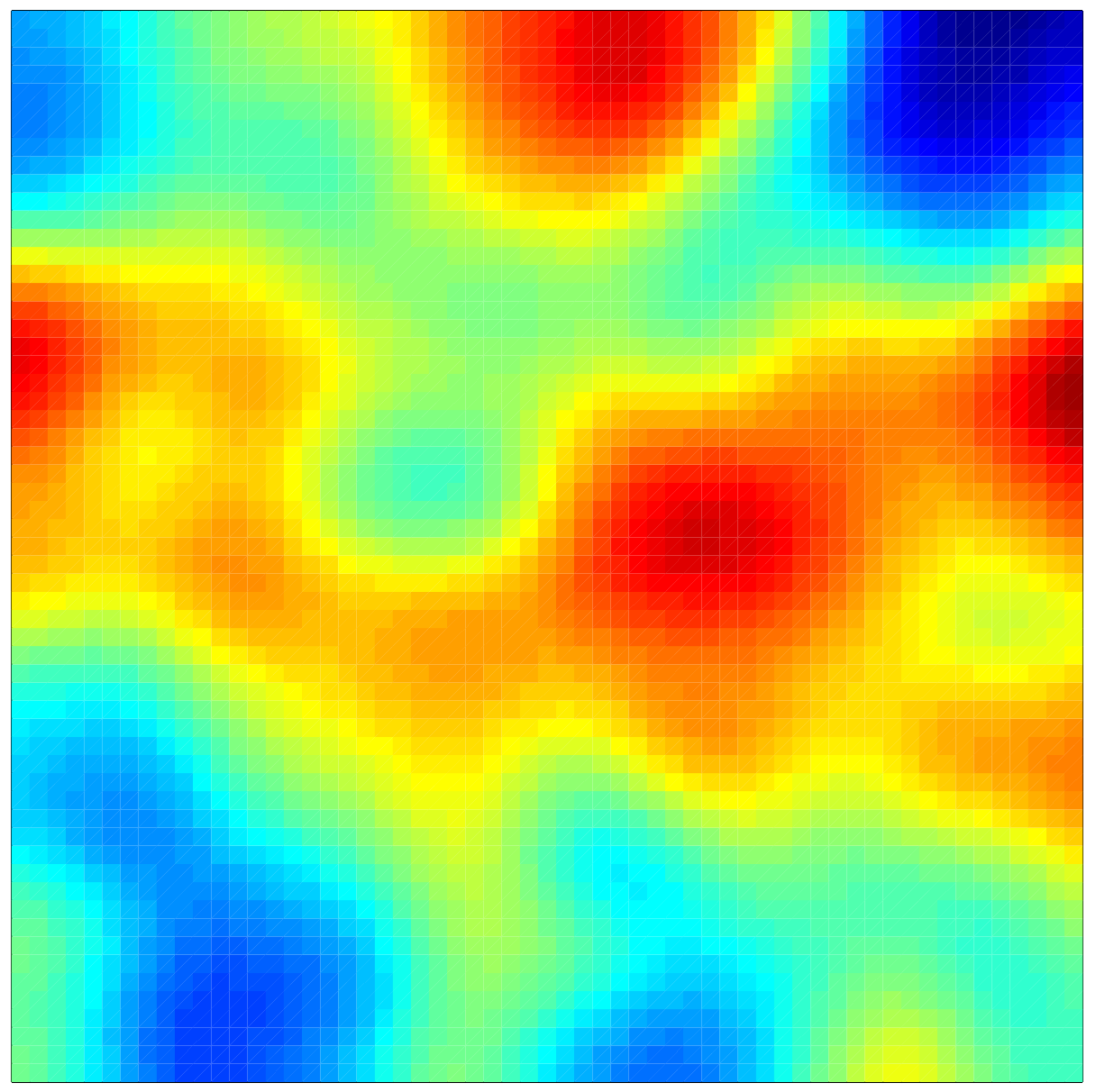}
\includegraphics[scale=0.23]{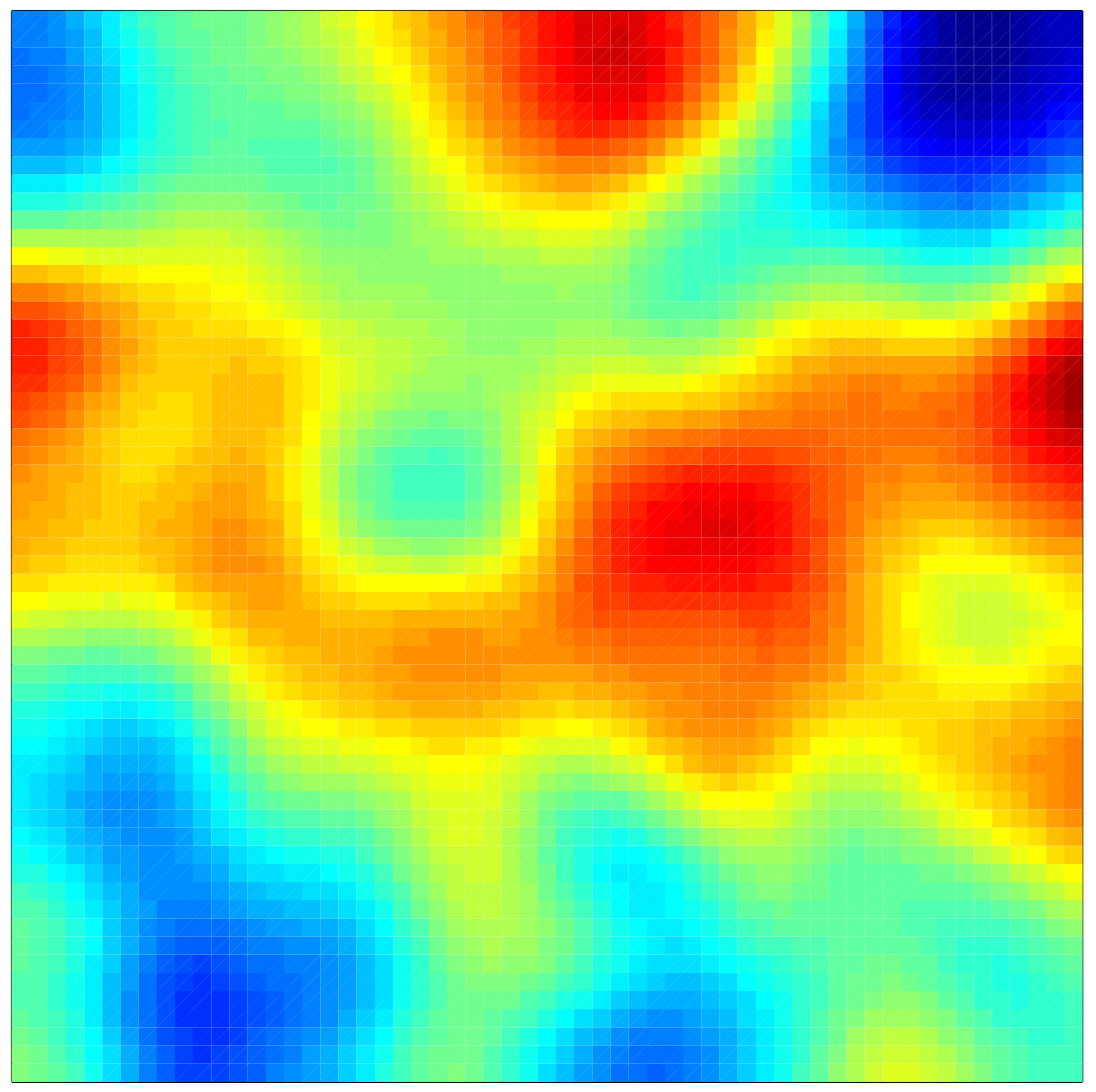}\\

\includegraphics[scale=0.23]{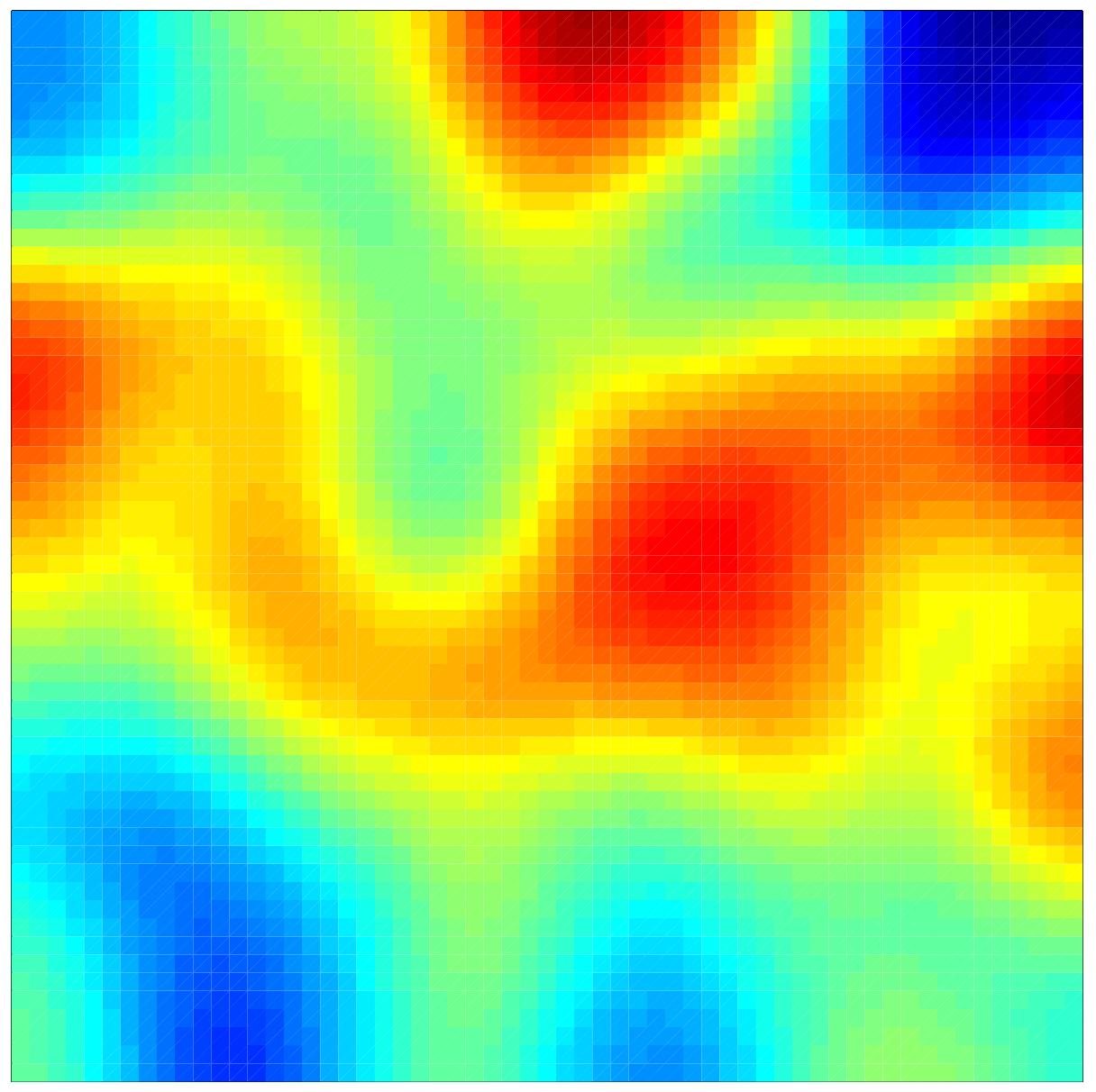}
\includegraphics[scale=0.23]{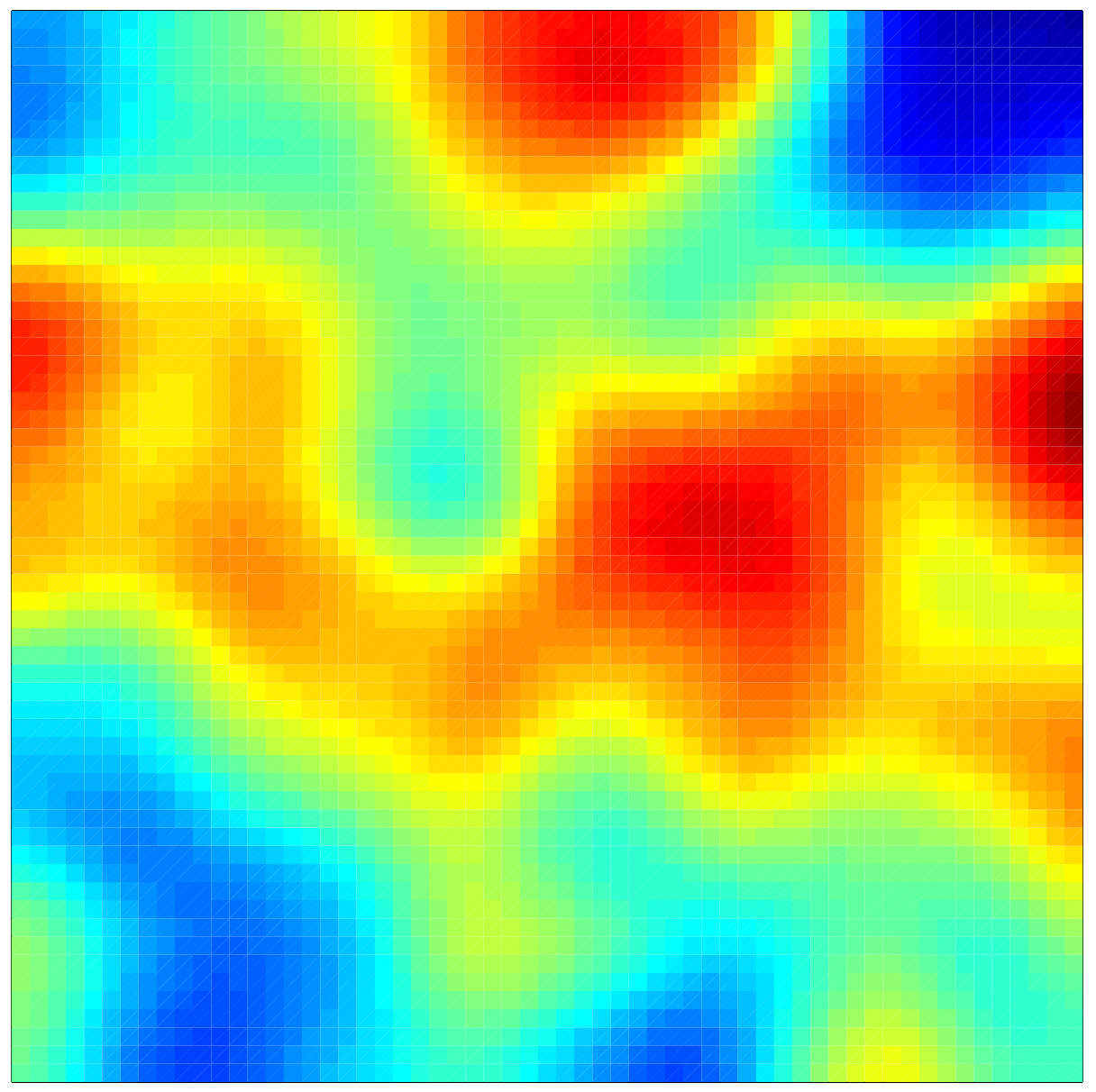}
\includegraphics[scale=0.23]{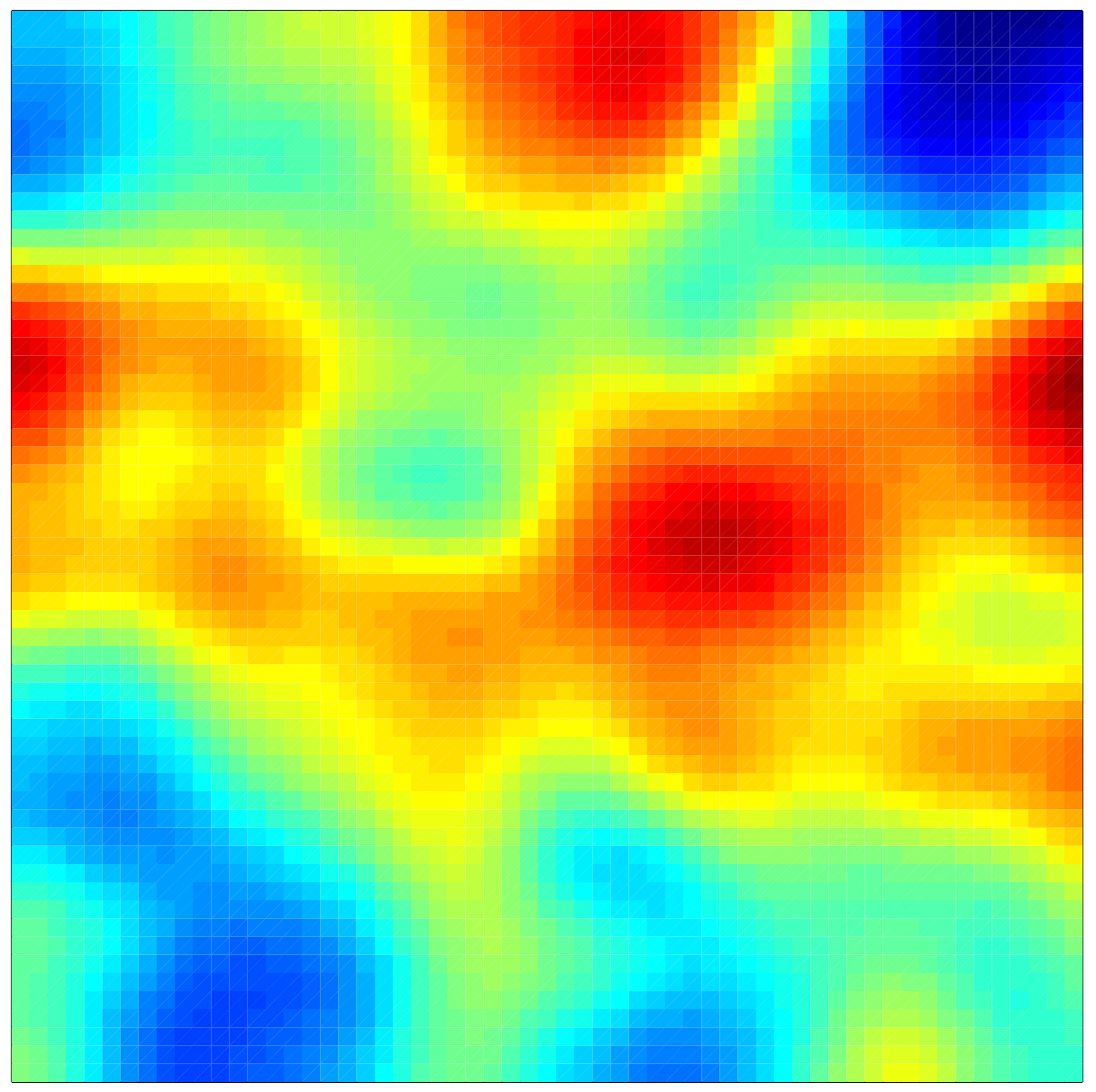}
\includegraphics[scale=0.23]{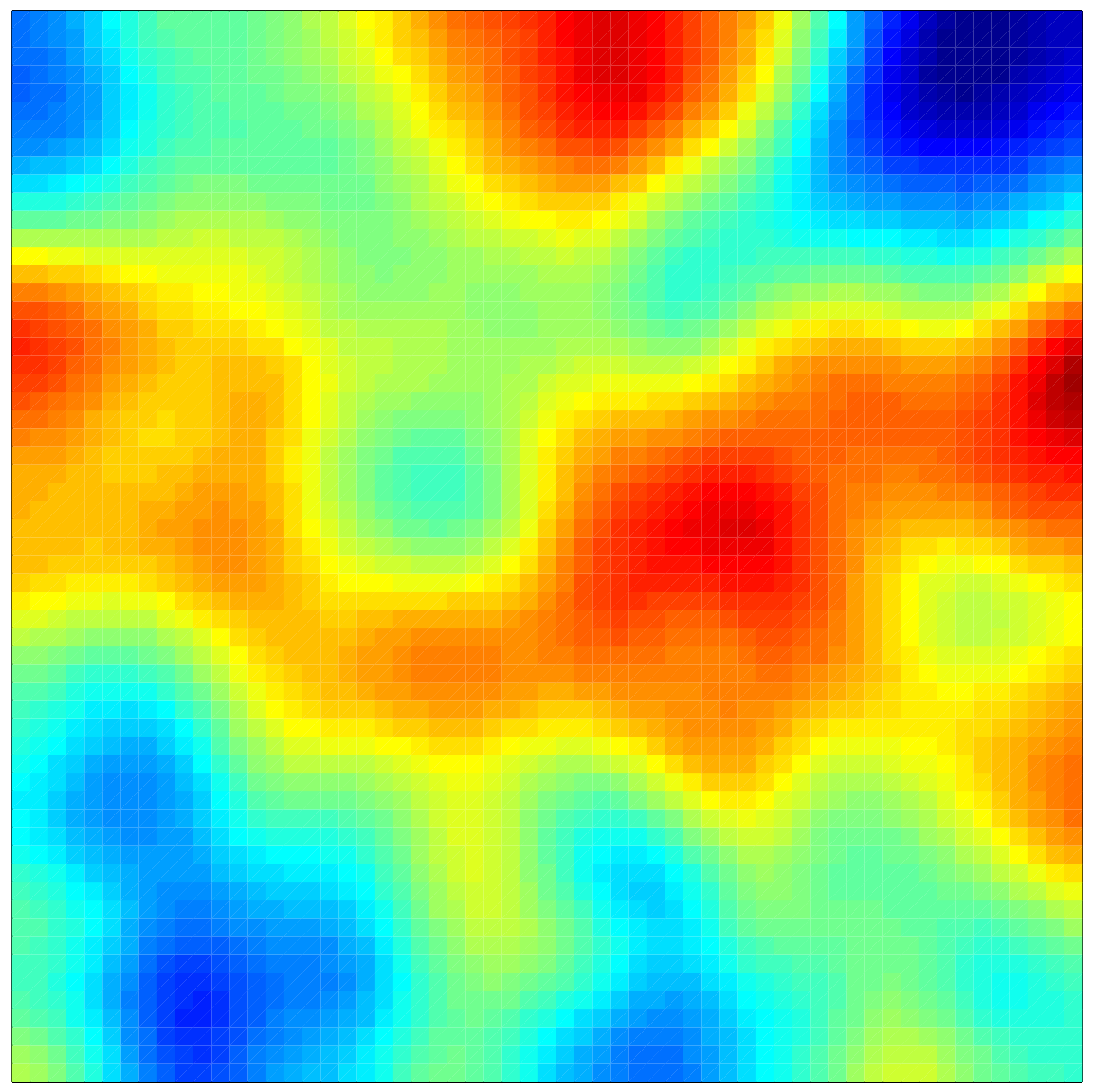}

 \caption{Top: truth $u^{\dagger}$. Top-middle, bottom-middle and bottom: Estimates obtained with 3DVAR and Data Model 1 at iterations (from left to right) 1, 10, 20, 30)  for noise levels of $1\%$ (top-middle), $2.5\%$ (bottom-middle) and $5\%$ (bottom) . } \label{Fig1}
\end{center}
\end{figure}

In this subsection we demonstrate how the filters under consideration in the iterative framework described in Algorithm \ref{Al1} can be used, with Data Model 1 and Data Model 2, to solve the linear inverse problem presented in Section \ref{se_Intro}. Let us consider the truth $u^{\dagger}$ displayed in Figure \ref{Fig1} (top) generated, as described in subsection \ref{se_num:setup}, from a Gaussian measure with covariance (\ref{num3}) and $s=1$. Synthetic data are generated as described above with three different choices of $\gamma$ that yield noise levels of approximately $1\%$,  $2.5\%$, and $5\%$ of the norm of the noise free data (i.e. $Au^{\dagger}$).

\begin{figure}[htbp]
\begin{center}

\includegraphics[scale=0.3]{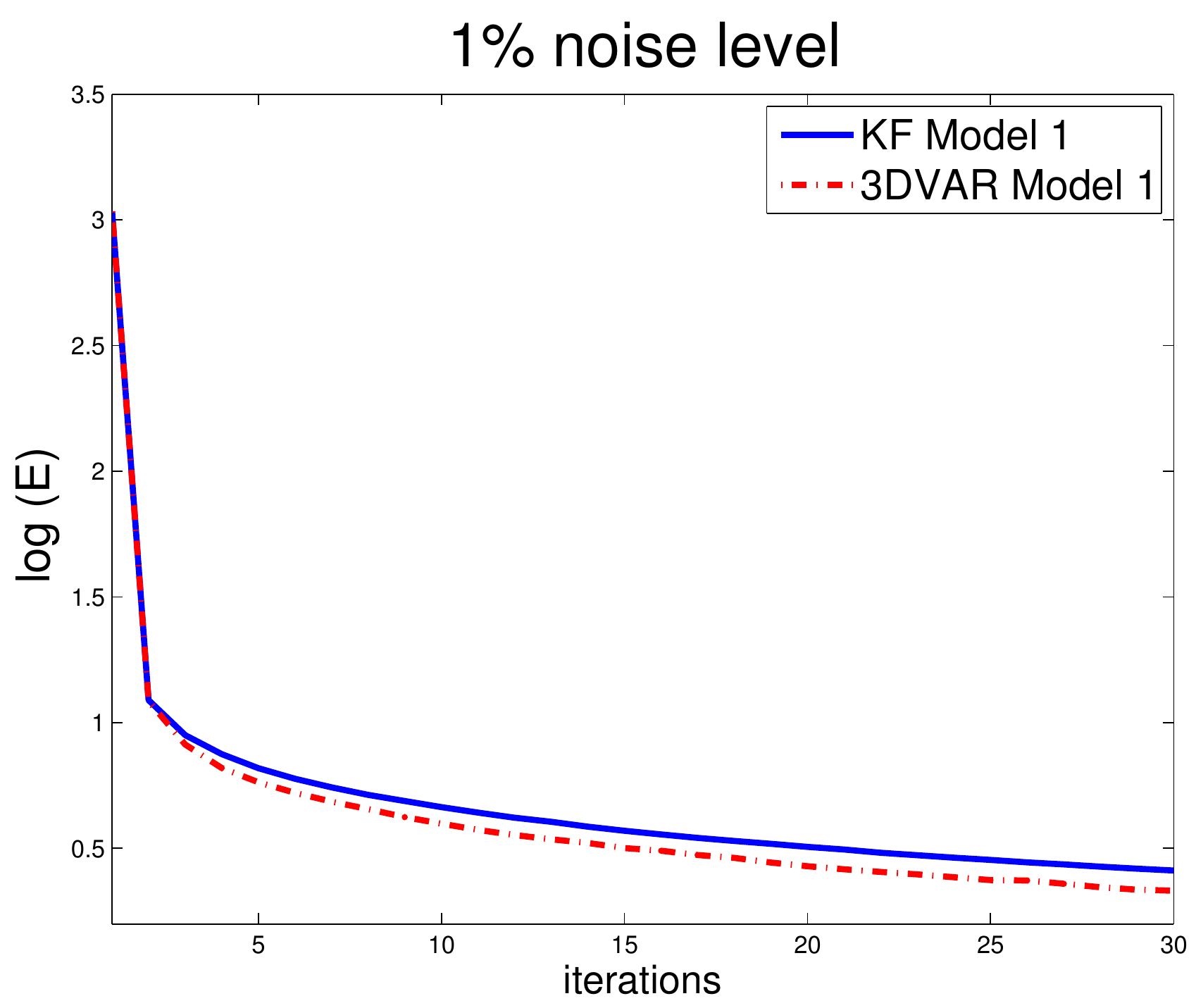}
\includegraphics[scale=0.3]{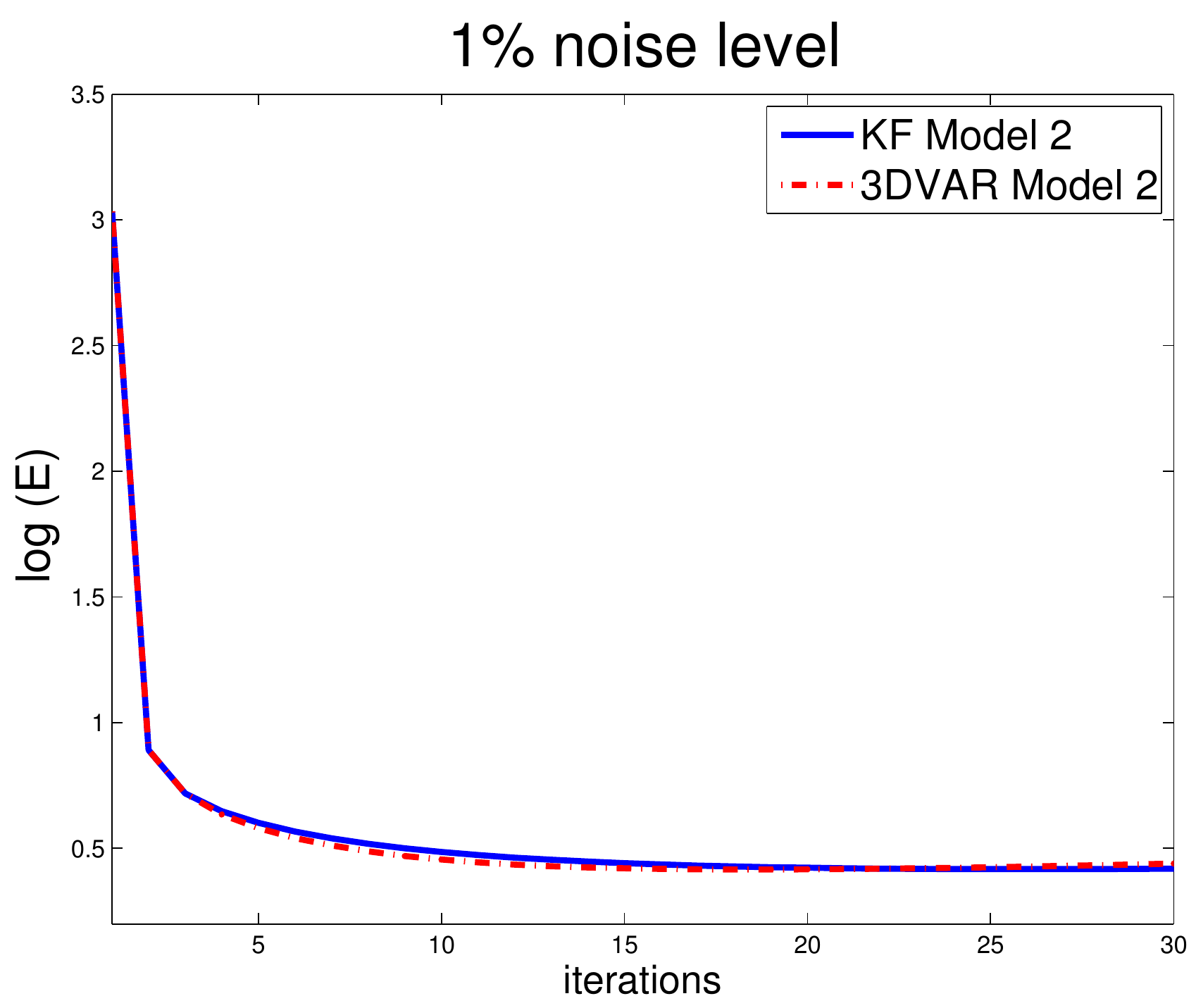}\\
\includegraphics[scale=0.3]{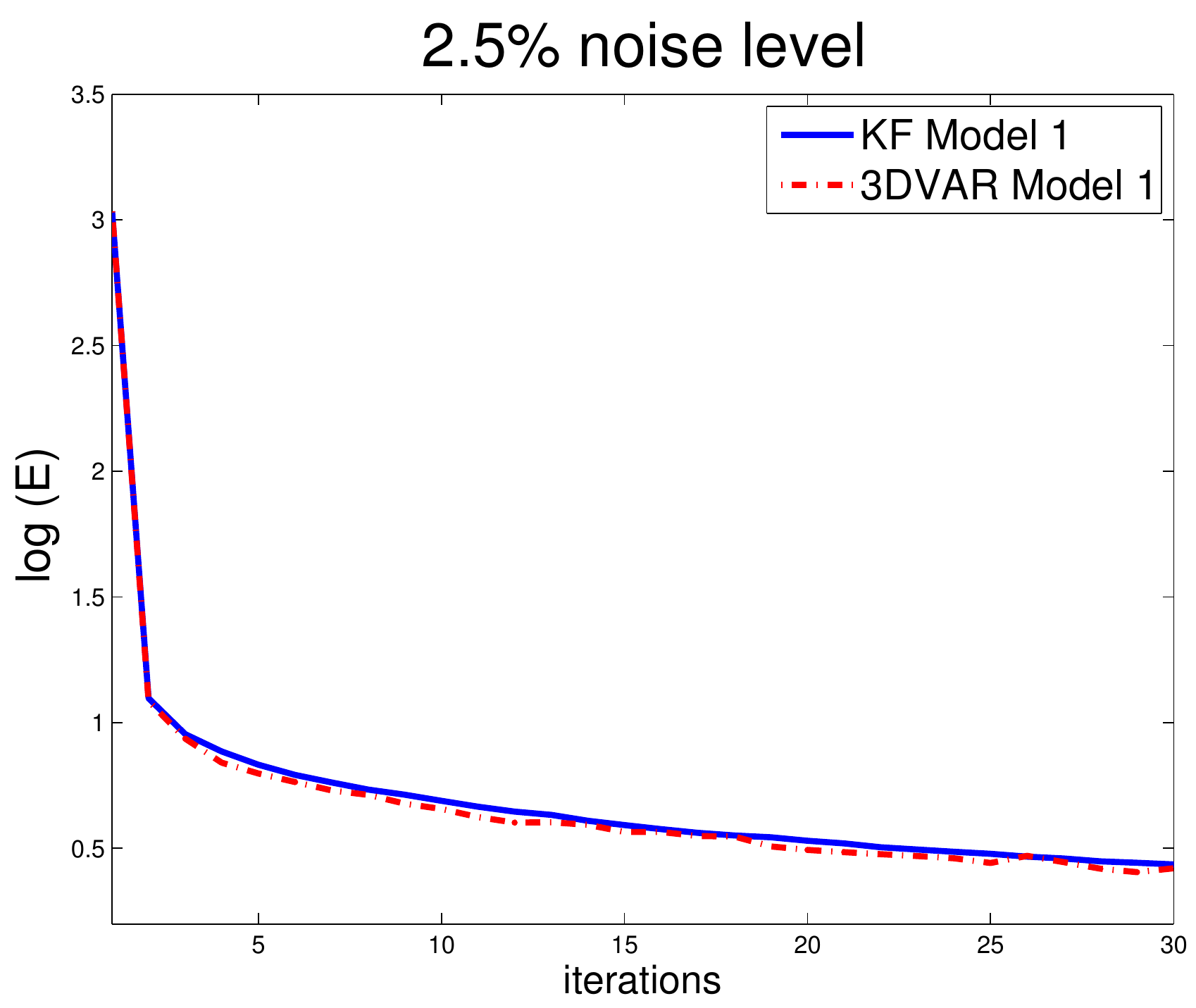}
\includegraphics[scale=0.3]{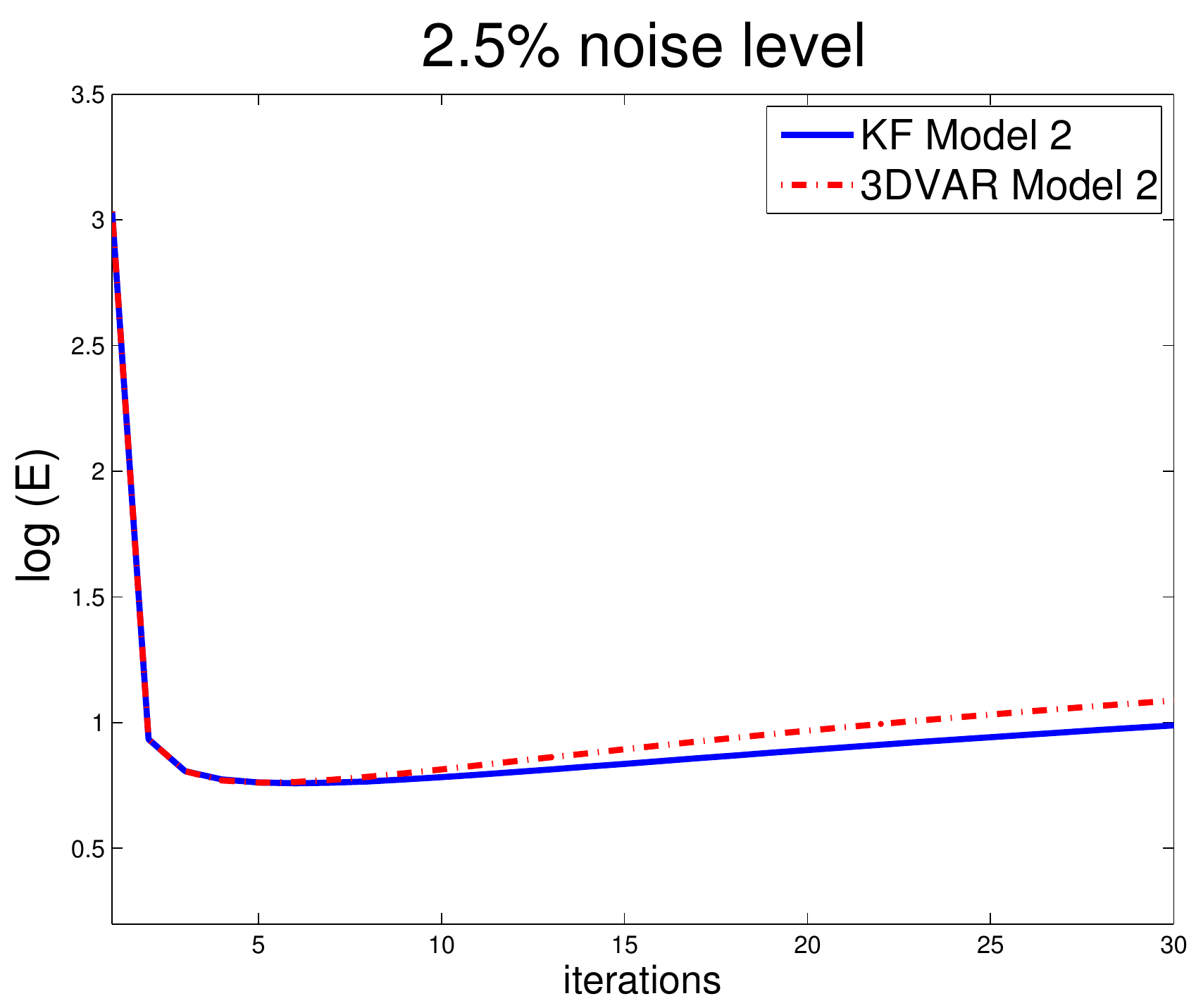}\\
\includegraphics[scale=0.3]{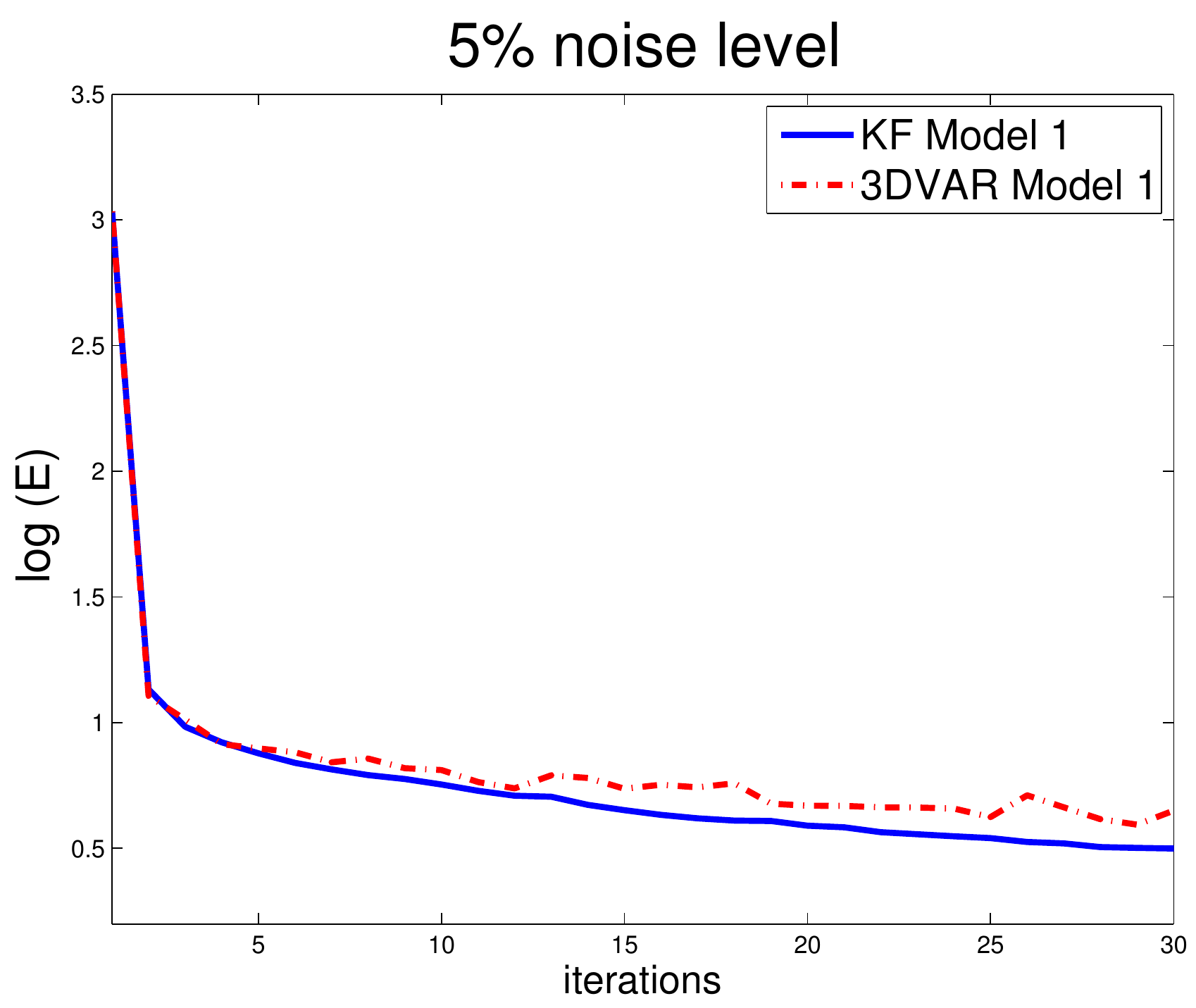}
\includegraphics[scale=0.3]{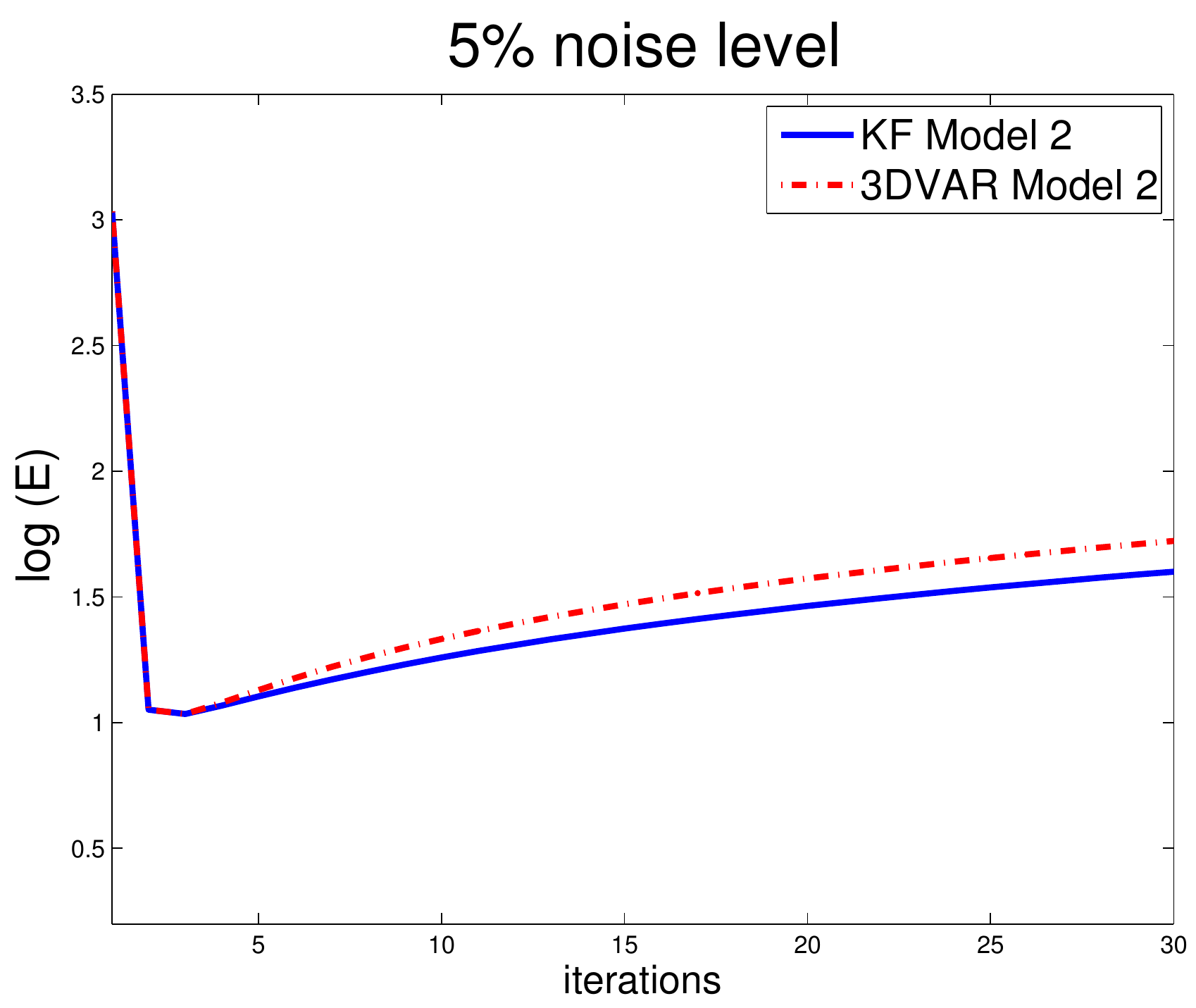}

 \caption{log Error with respect to the truth vs iterations of 3DVAR and Kalman filter applied to Data Model 1 (left) and Data Model 2 (right) for different noise levels} \label{Fig2}
\end{center}
\end{figure}

We apply Algorithm \ref{Al1} to this synthetic data generated both for the application of Data Model $1$ and Data Model $2$. For Data Model $1$ we consider a selection of $N=25$. Algorithm \ref{Al1} states that the these schemes should be stopped at iteration level $n=N$. However, in order to observe the performance of these schemes, in these experiments we allowed for a few more iterations. In the left column of Figure \ref{Fig2} we display the plots of the error w.r.t the truth of the estimator $m_{n}$ as a function of the iterations, i.e.
\begin{eqnarray*}
E(m_{n})\equiv \| m_n-\mathcal{P}u^{\dagger} \|
\end{eqnarray*}
where $\mathcal{P}u^{\dagger}$ denotes the interpolation of $u^{\dagger}$ on the coarse grid used for the inversion. Note that the error w.r.t. the truth of the estimates produced by both schemes decreases monotonically.  Interestingly, the value at the final iteration displayed in these figures is approximately the same for all these experiments independently of the noise level. Moreover, the stability of the scheme does not seem to depend on the early termination of the scheme. In addition, we note that both Kalman filter and 3DVAR exhibit very similar performance in terms of reducing the error w.r.t the truth. However, for larger noise levels we observe small fluctuations in the error obtained with 3DVAR.

In Figure \ref{Fig1} we display the estimates $m_n$ obtained with 3DVAR with Data Model $1$ at iterations $n=1, 10, 20, 30$ for noise levels (determined by the choice of $\gamma$) of $1\%$ (top-middle), $2.5\%$ (middle-bottom) and $5\%$ (bottom). We can visually appreciate from Figure \ref{Fig1} that the estimates obtained at the final iteration $n=30$ are indeed similar one to another even though the one in the bottom row was computed by inverting data five times noisier than the one in the first row. Similar estimates (not shown) were obtained with the Kalman filter for Data Model $1$.

\begin{figure}[htbp]
\begin{center}
\includegraphics[scale=0.23]{Truth.pdf}\\

\includegraphics[scale=0.23]{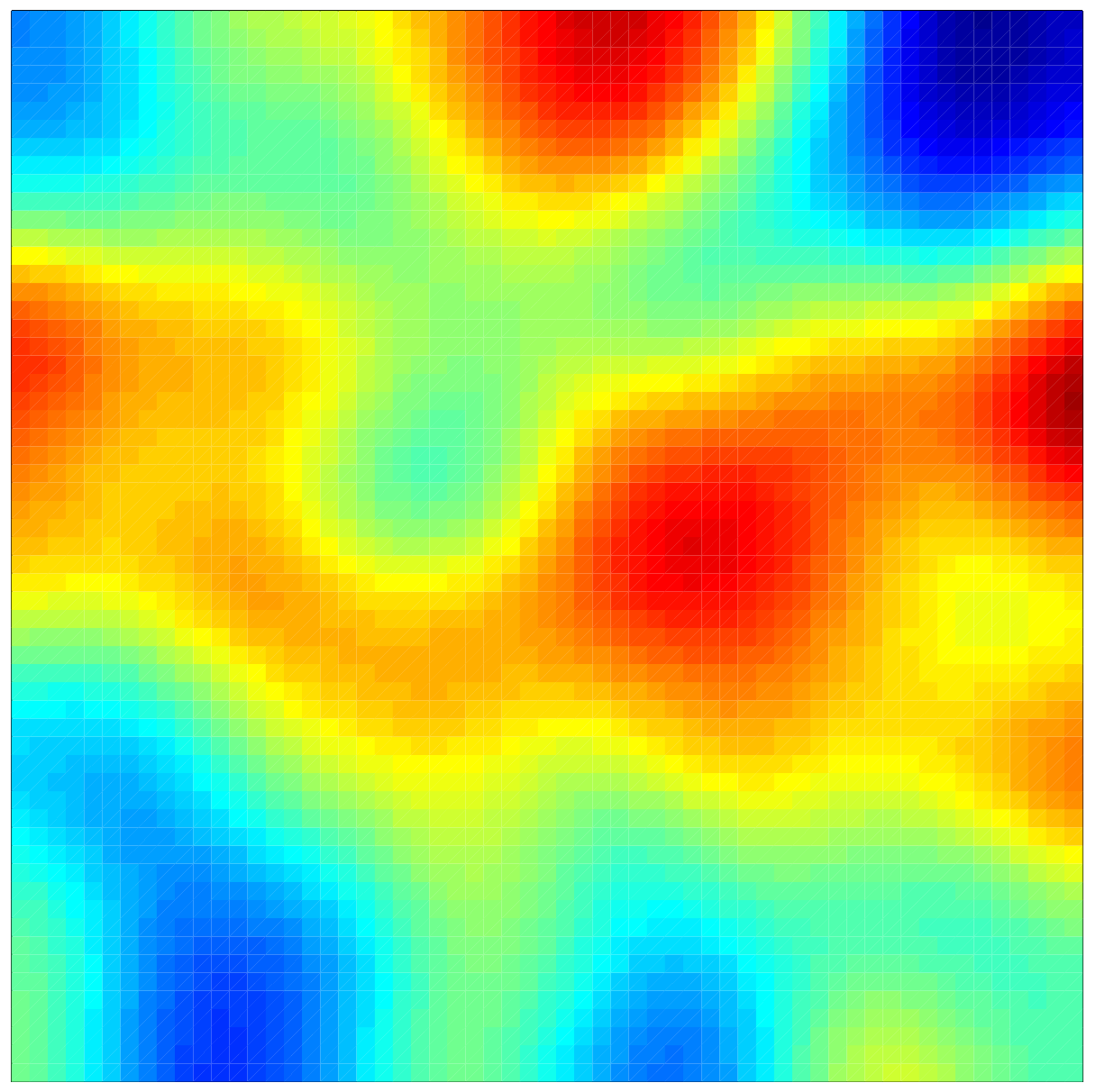}
\includegraphics[scale=0.23]{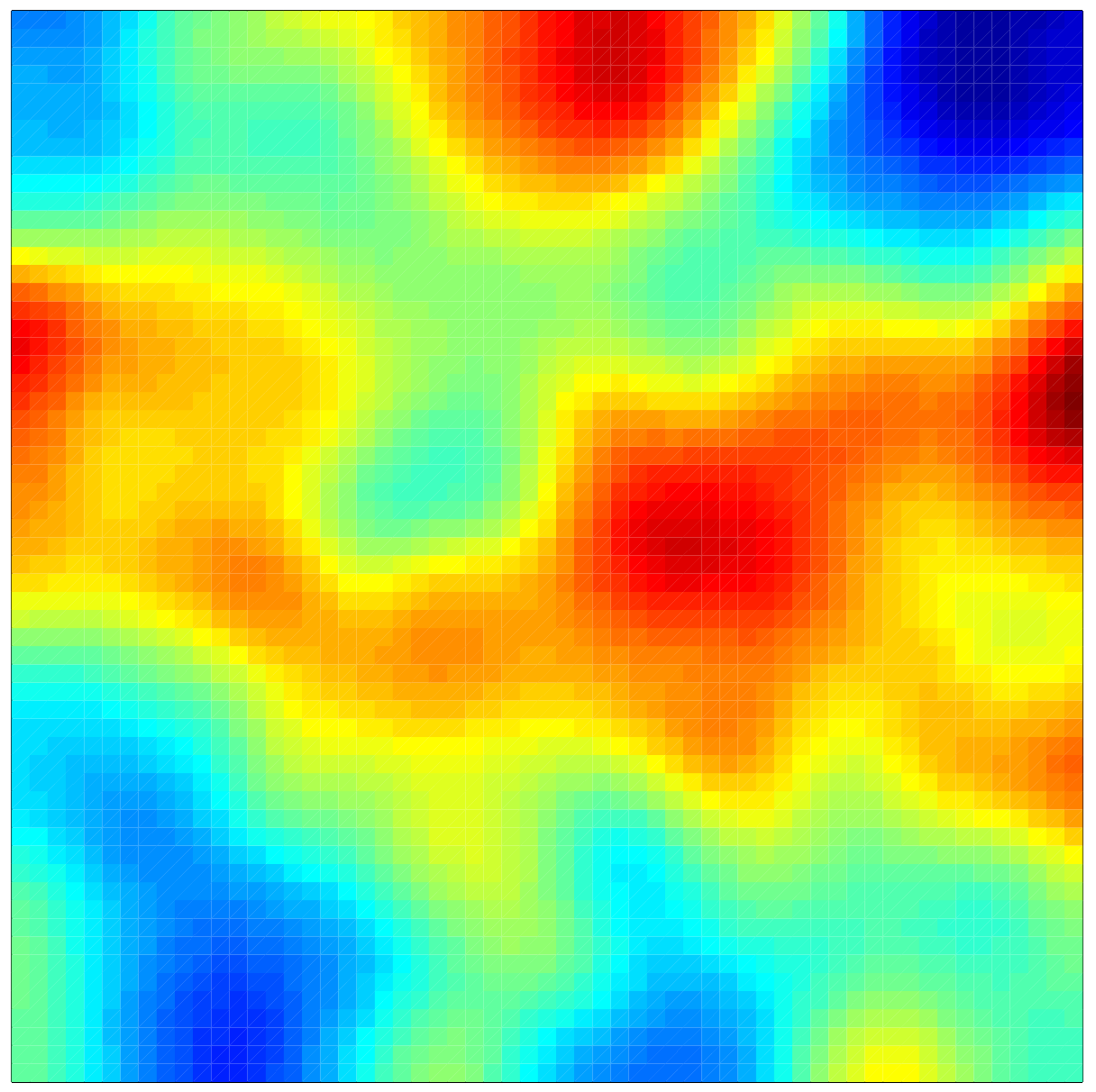}
\includegraphics[scale=0.23]{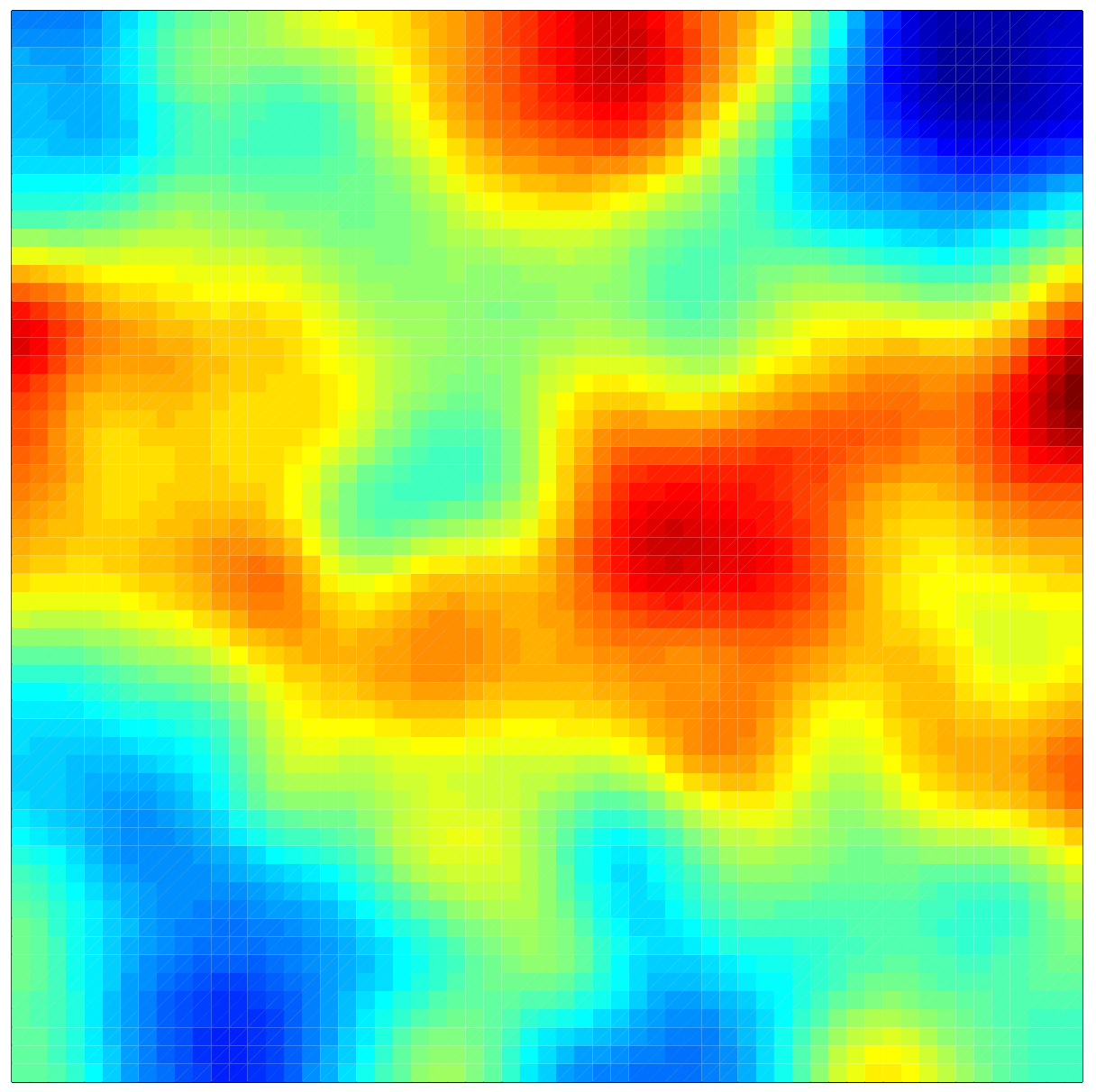}
\includegraphics[scale=0.23]{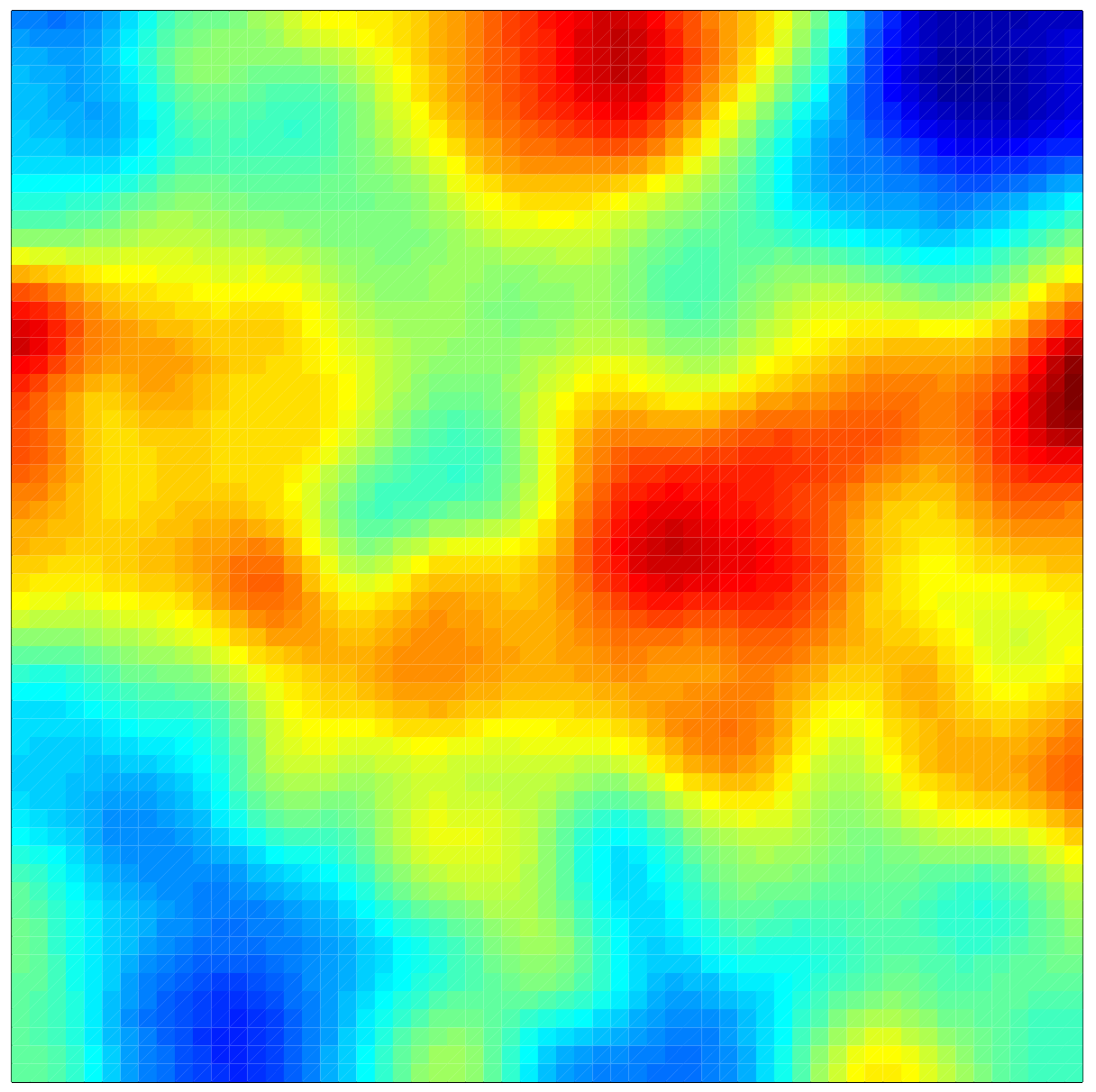}\\

\includegraphics[scale=0.23]{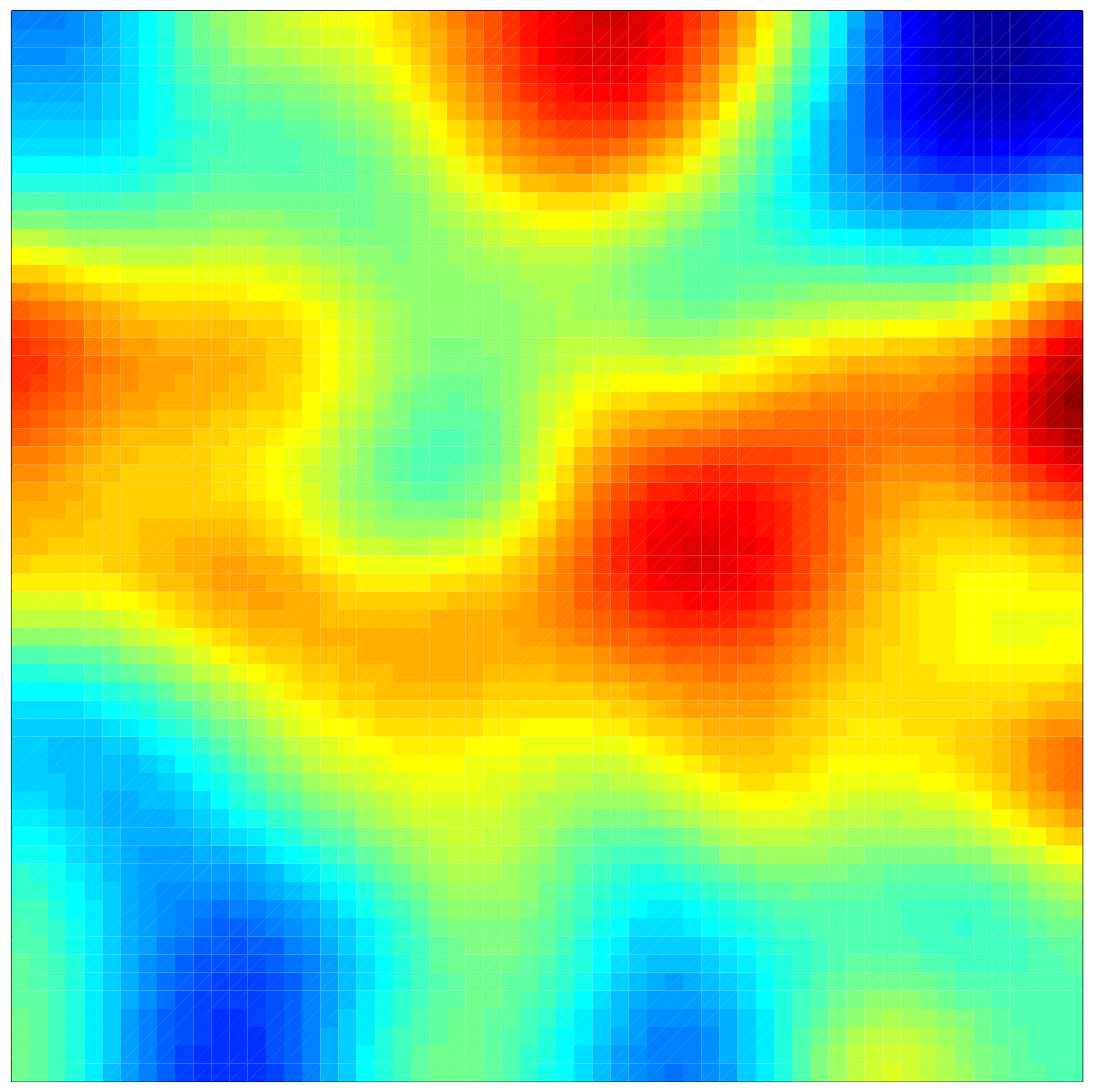}
\includegraphics[scale=0.23]{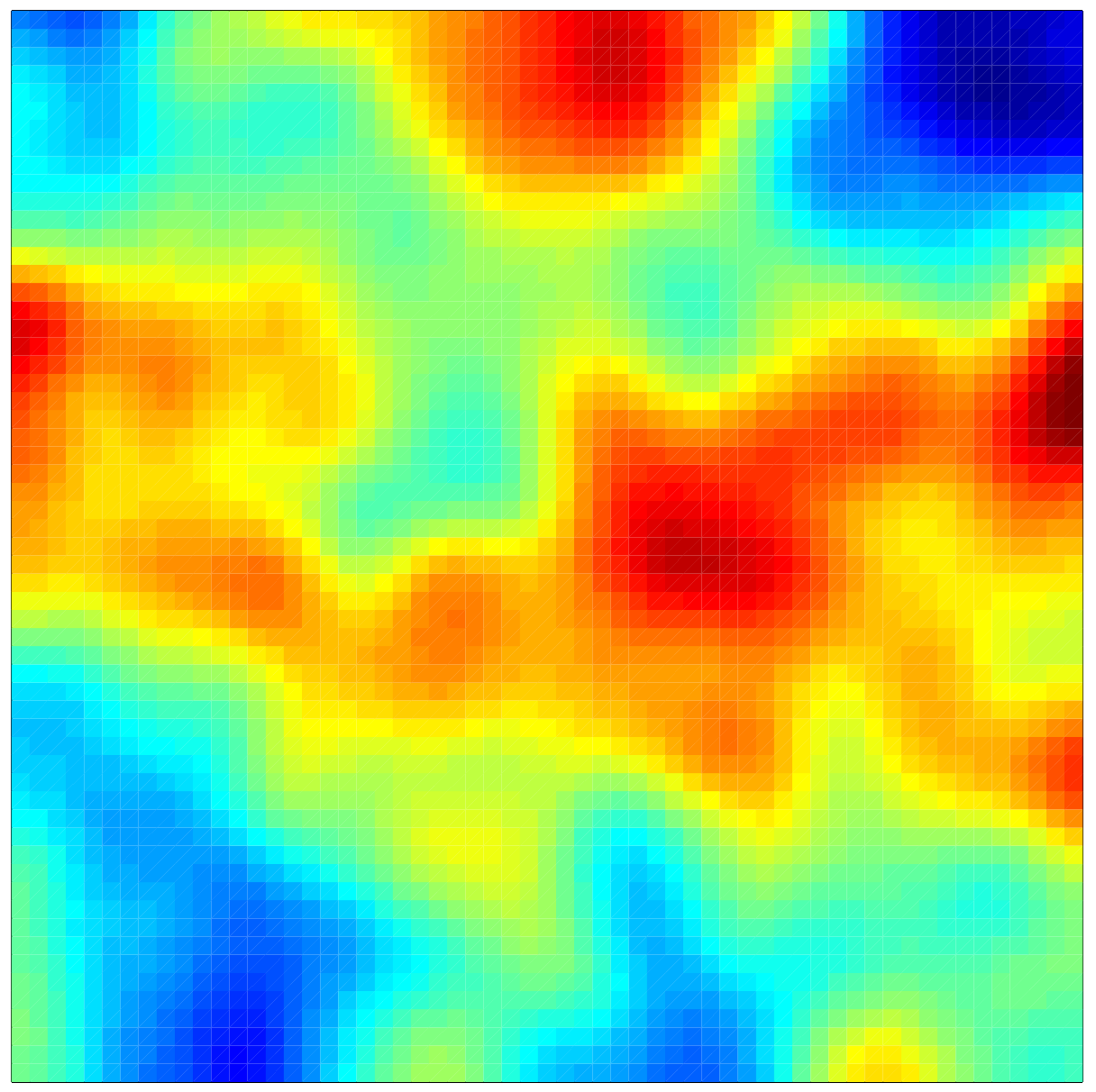}
\includegraphics[scale=0.23]{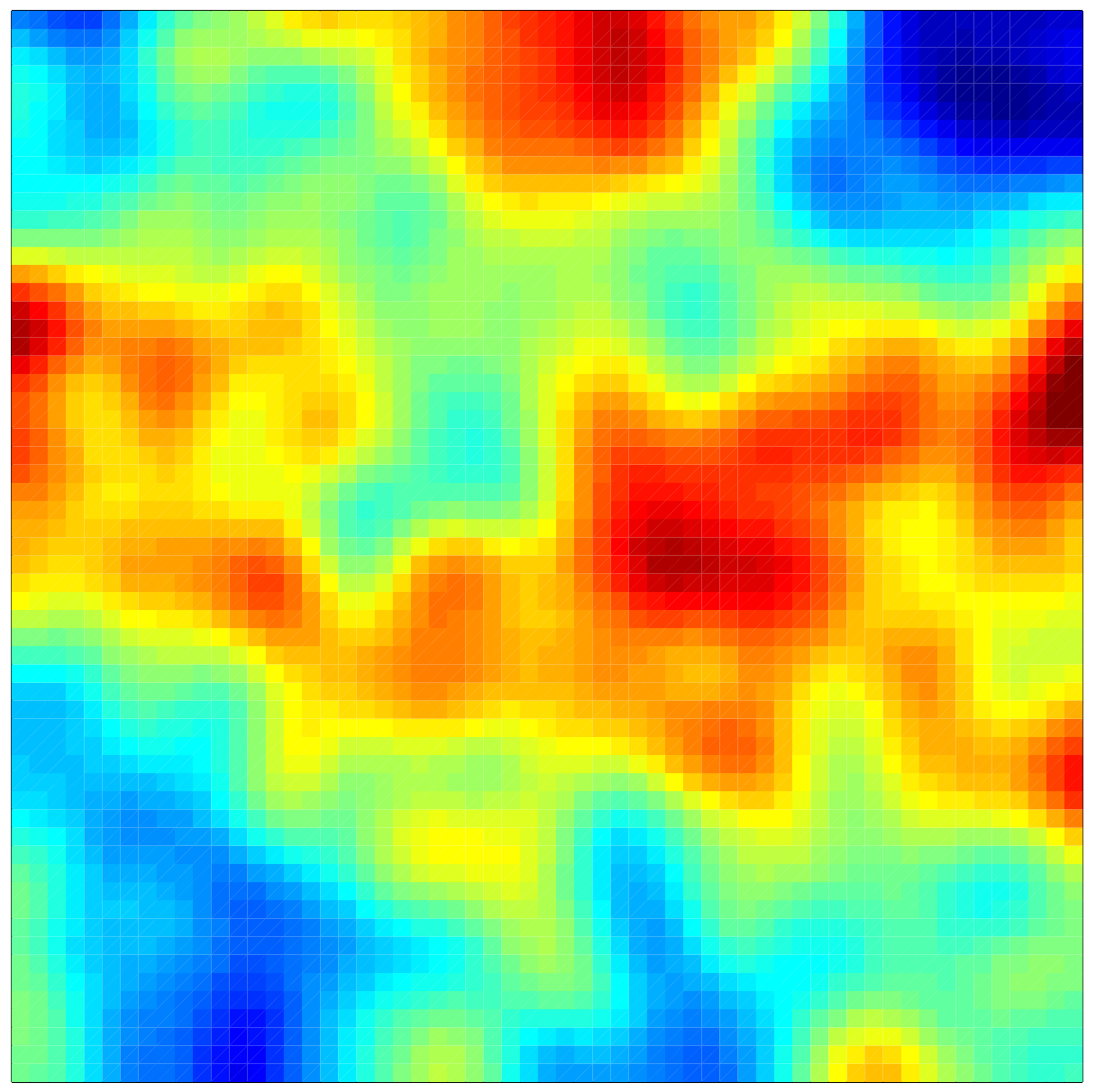}
\includegraphics[scale=0.23]{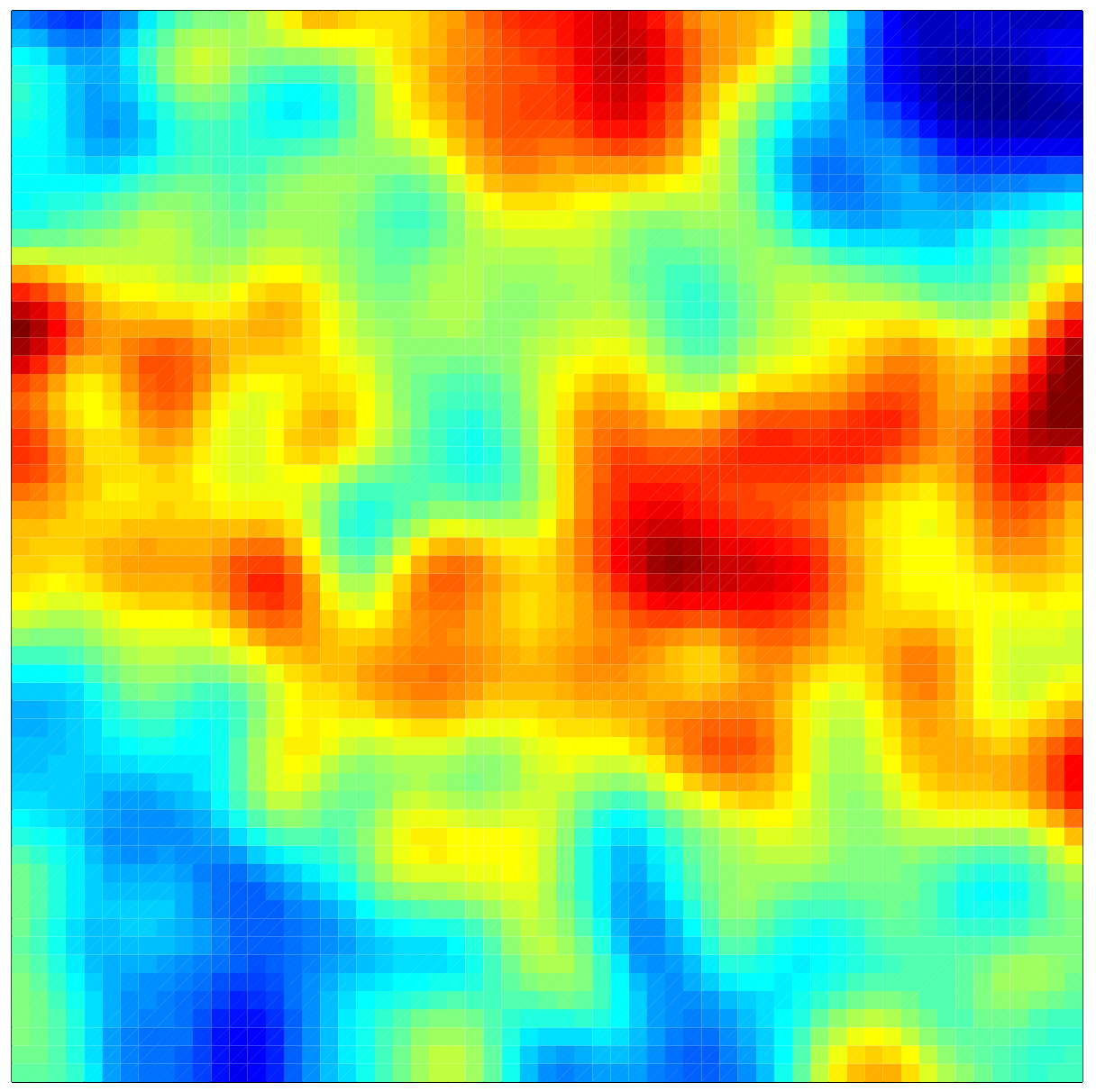}\\

\includegraphics[scale=0.23]{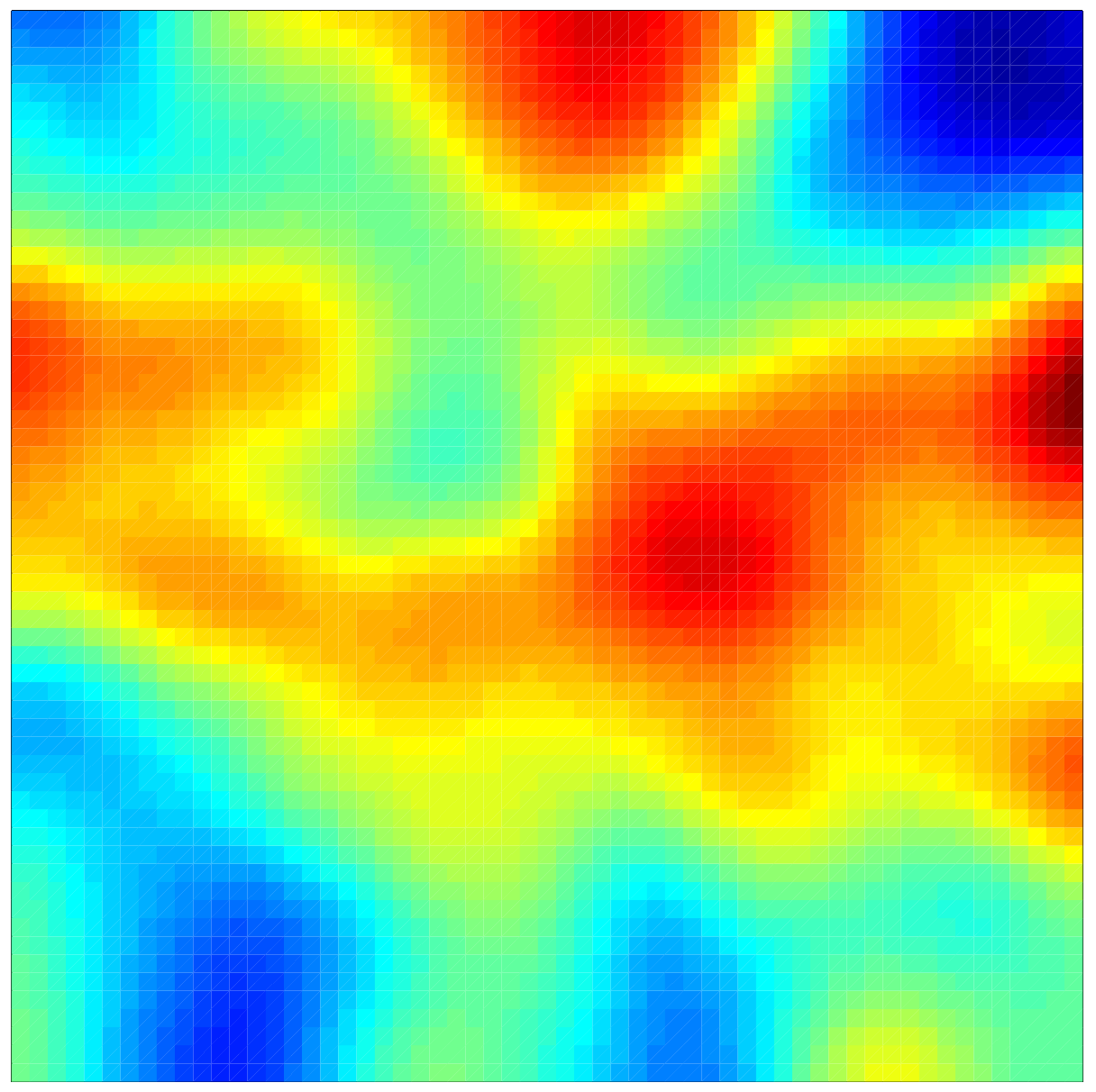}
\includegraphics[scale=0.23]{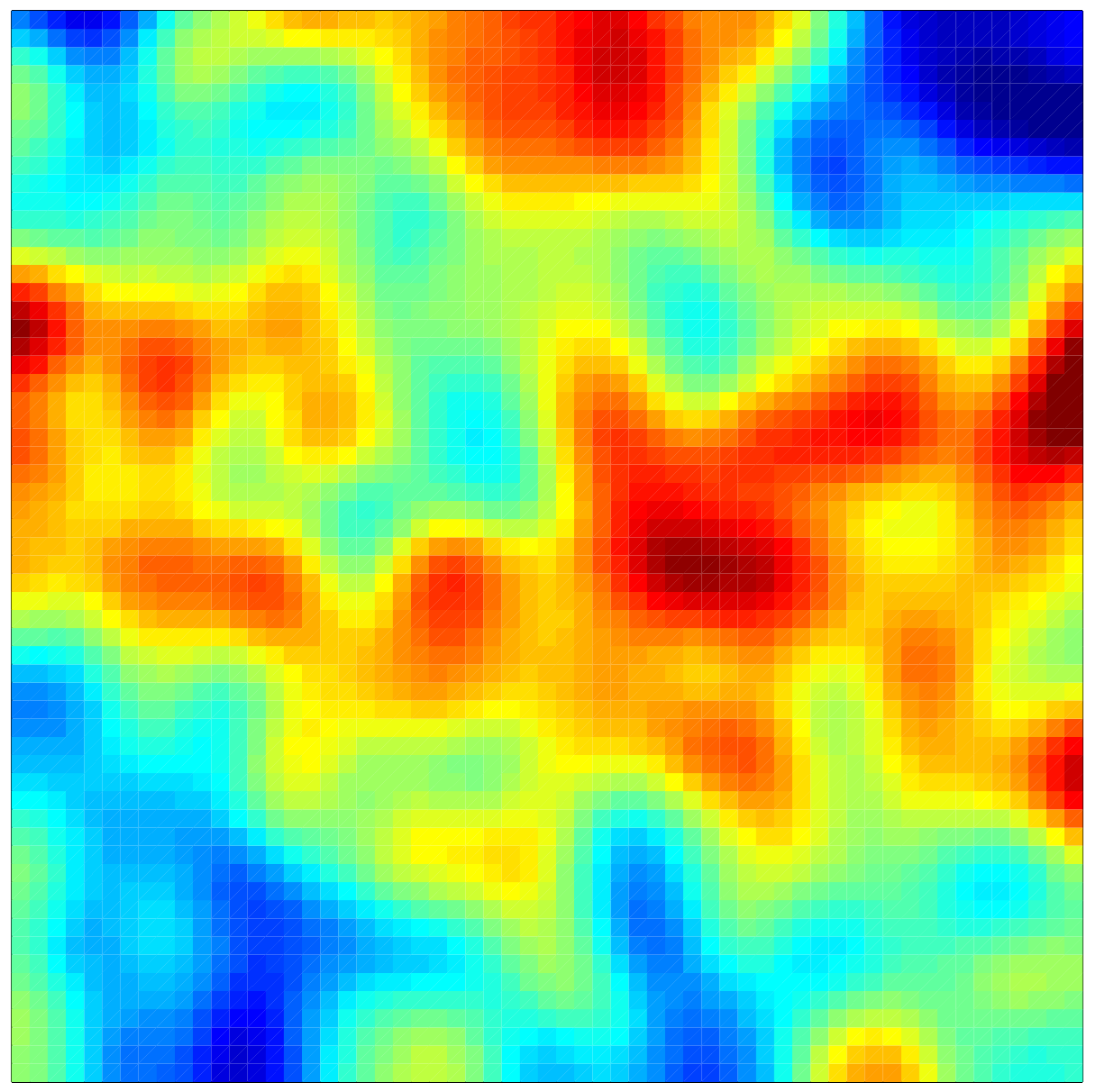}
\includegraphics[scale=0.23]{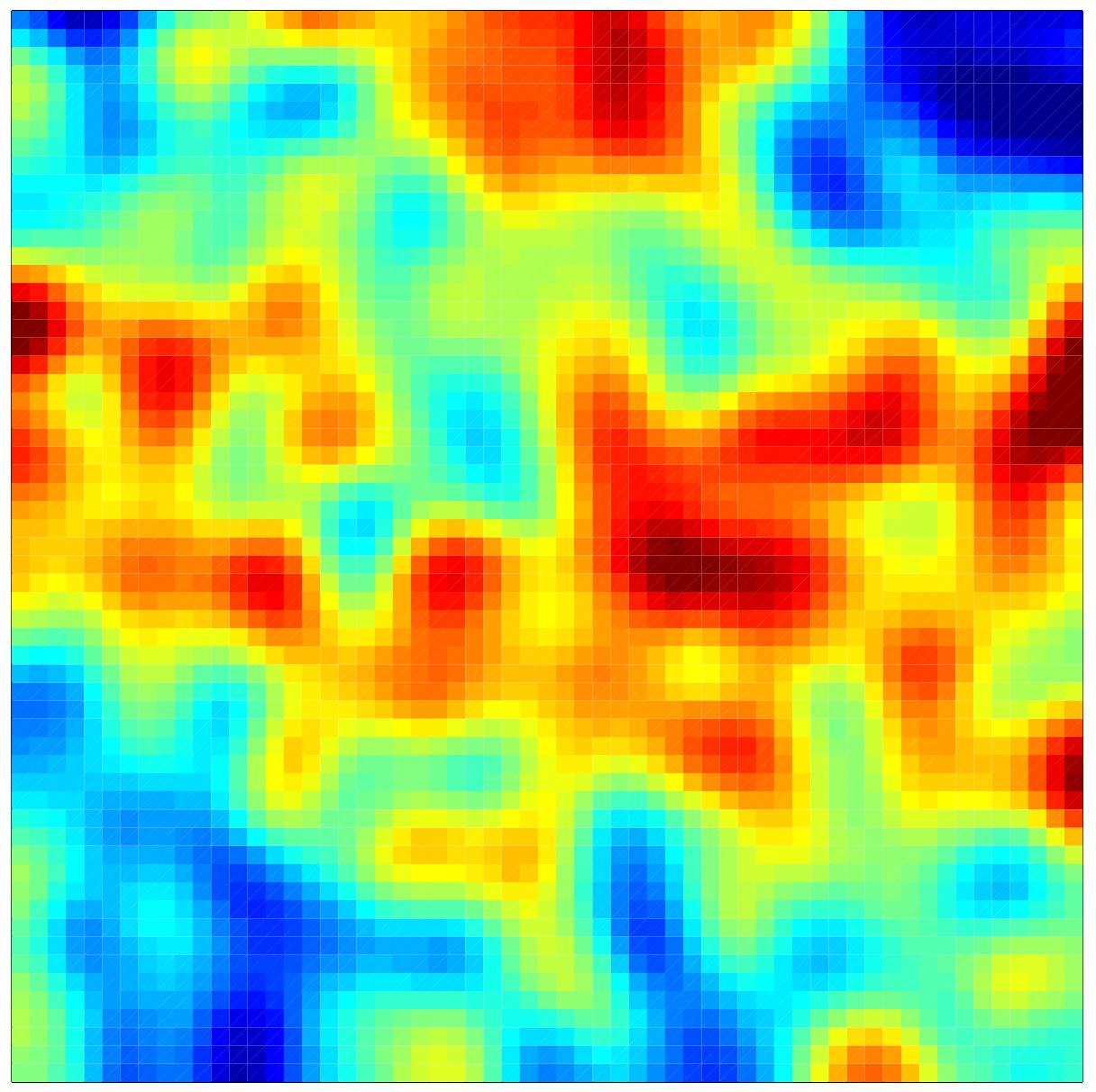}
\includegraphics[scale=0.23]{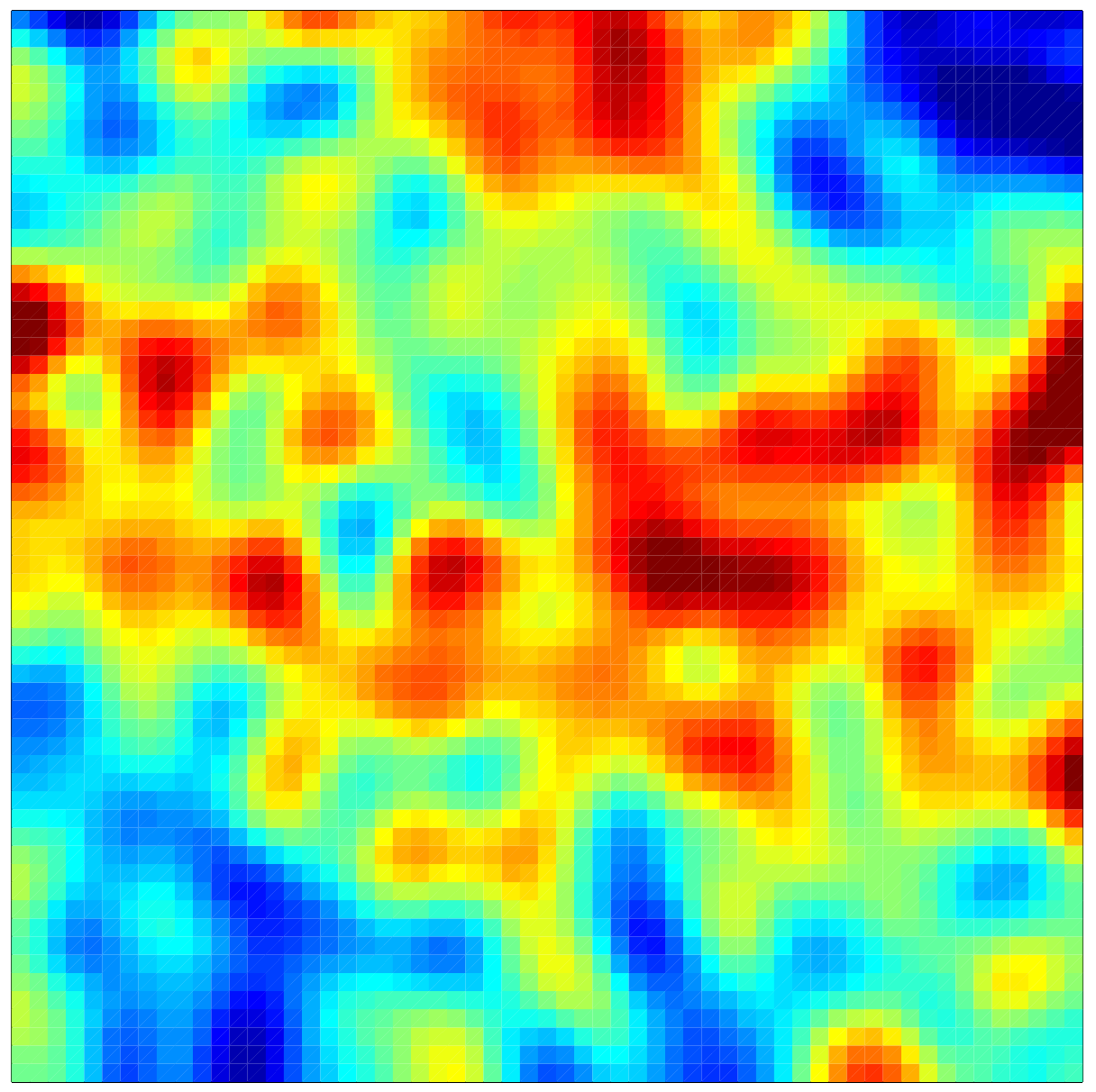}

 \caption{Top: truth $u^{\dagger}$. Top-middle, bottom-middle and bottom: Estimates obtained with 3DVAR and Data Model 2 at iterations (from left to right) 1, 10, 20, 30) for noise levels of $1\%$ (top-middle), $2.5\%$ (bottom-middle) and $5\%$ (bottom)  }. \label{Fig3}
\end{center}
\end{figure}

For Data Model $2$, the selection of $\gamma$ corresponding to noise levels of $1\%$,  $2.5\%$, and $5\%$ yields a maximum number of iterations $N=20,6,3$ respectively. Clearly, for Data Model $2$, smaller observational noise results in schemes that can be iterated longer, presumably achieving more accurate estimates. Similarly to Data Model $1$, we are required to stop the algorithm at the iteration $n=N$. However, for the purpose of this study we iterate until $n=30$. In the right column of Figure \ref{Fig2} we display the plots of the error w.r.t the truth of $m_n$. In contrast to Data Model $1$, we observe that the error w.r.t the truth increases for $n>N$. In other words, the Kalman filter and 3DVAR, when applied to Data Model $2$, need to be stopped at $n=N$ in order to stabilize the scheme and obtain accurate estimates of the truth.  Moreover, stopping the scheme at $n=N$ results in estimates whose accuracy increases with smaller noise levels. Clearly, in the small noise limit, both data models tend to exhibit similar behaviour. In Figure \ref{Fig3} we display $m_n$ obtained from 3DVAR
applied to Data Model $2$ at iterations $n=1, 10, 20, 30$ for data with noise levels of $1\%$ (top-middle), $2.5\%$ (middle-bottom) and $5\%$ (bottom). Similar estimates (not shown) were obtained with the Kalman Filter for Data Model $2$.

It is clear indeed, that for Data Model $1$, the application of Kalman filter and 3DVAR results in more accurate and stable estimates when the noise level in the data is not sufficiently small. The results from this subsection show that the reduction in the variance of the noise due to the law of large numbers
and central limit theorem effect in Data Model $1$ produces more accurate algorithms. Although Data Model $1$ requires multiple instances of the data, in some applications such as in EIT \cite{EIT}, the data collection can be repeated multiple times in order to obtain these data.

\subsection{Numerical verification of convergence rates}

In this subsection we test the convergence rates from Theorems \ref{thm_KFDM1}, \ref{thm_3DVARDM1}, \ref{thm_KFDM2}, and \ref{thm_3DVARDM2}. For the verification of each of these rates we let $\Sigma:=\Sigma_s$ denote the covariance from \eqref{num3}, and we consider $20$ experiments corresponding
to different truths $u^{\dagger}$ generated from $N(0,\Sigma_{s})$ with
the four selections of $s=1,2,3,4$. Note that Assumption \ref{assp_main} (iii) is satisfied only for $s\leq 3$.  Again, for Data Model 1 we generate synthetic data for each of the truths and for as many iterations used for the application of both schemes. Inverse crimes are avoided as described in subsection \ref{se_num:setup}.

\begin{figure}[htbp]
\begin{center}

\includegraphics[scale=0.32]{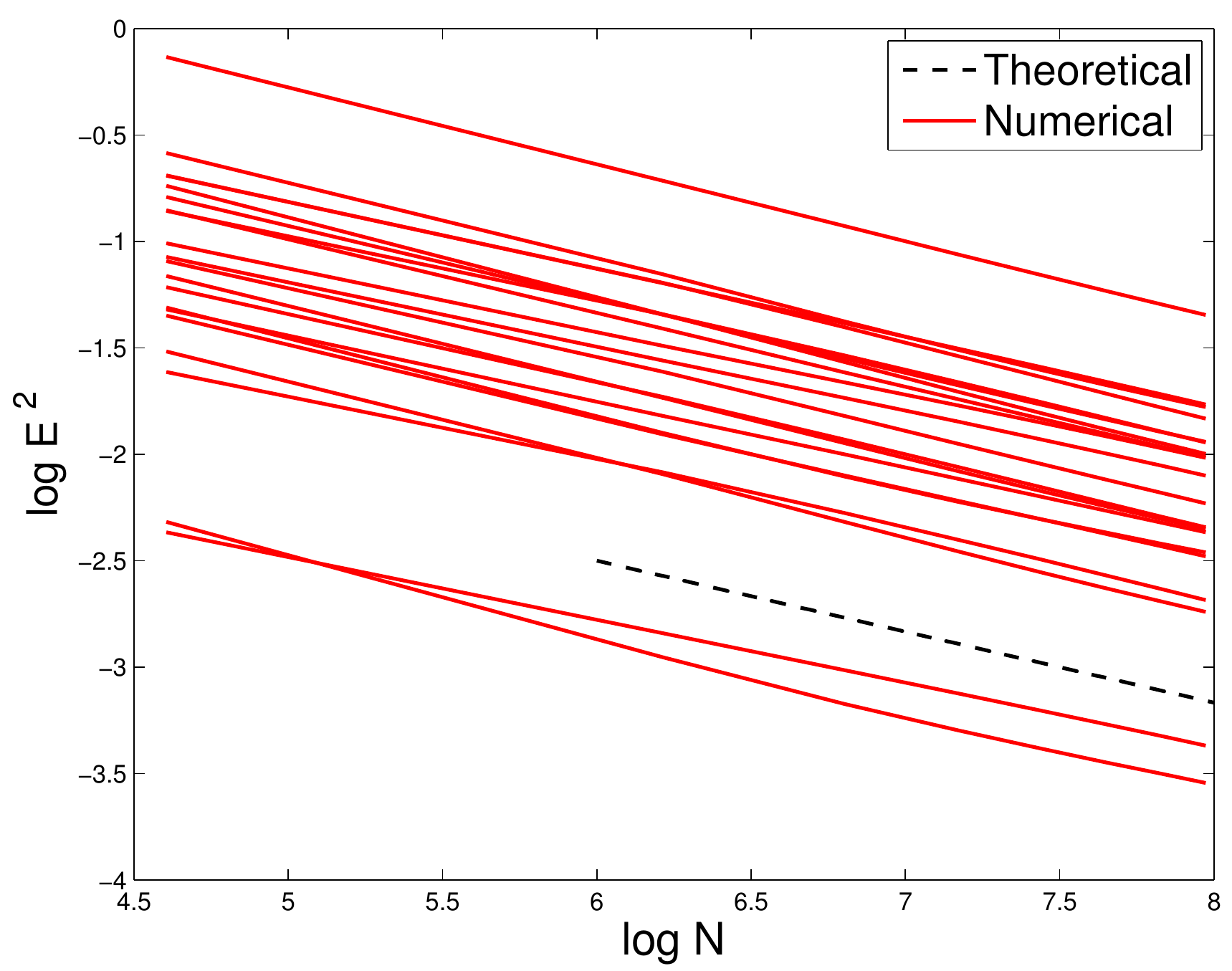}
\includegraphics[scale=0.32]{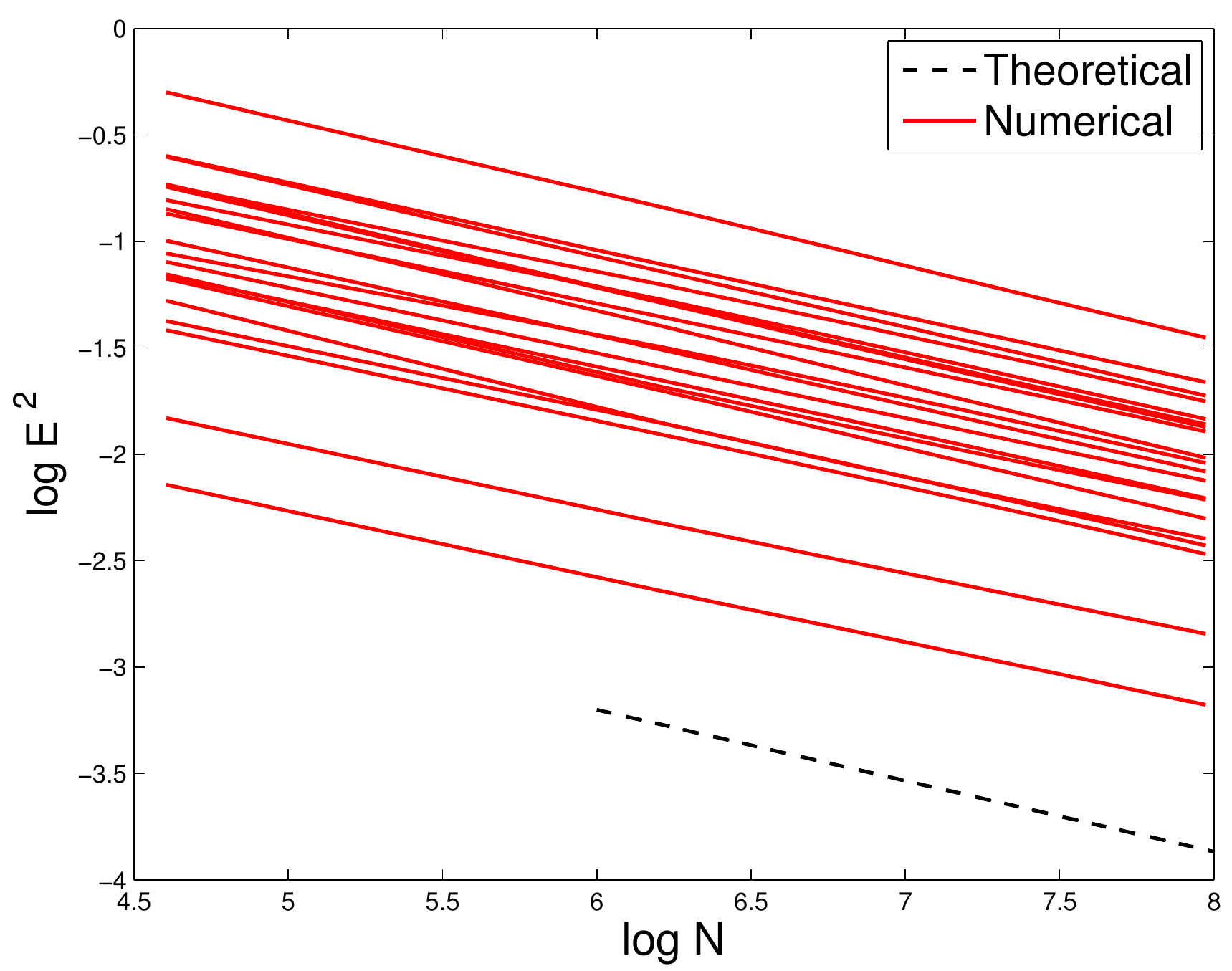}\\
\includegraphics[scale=0.32]{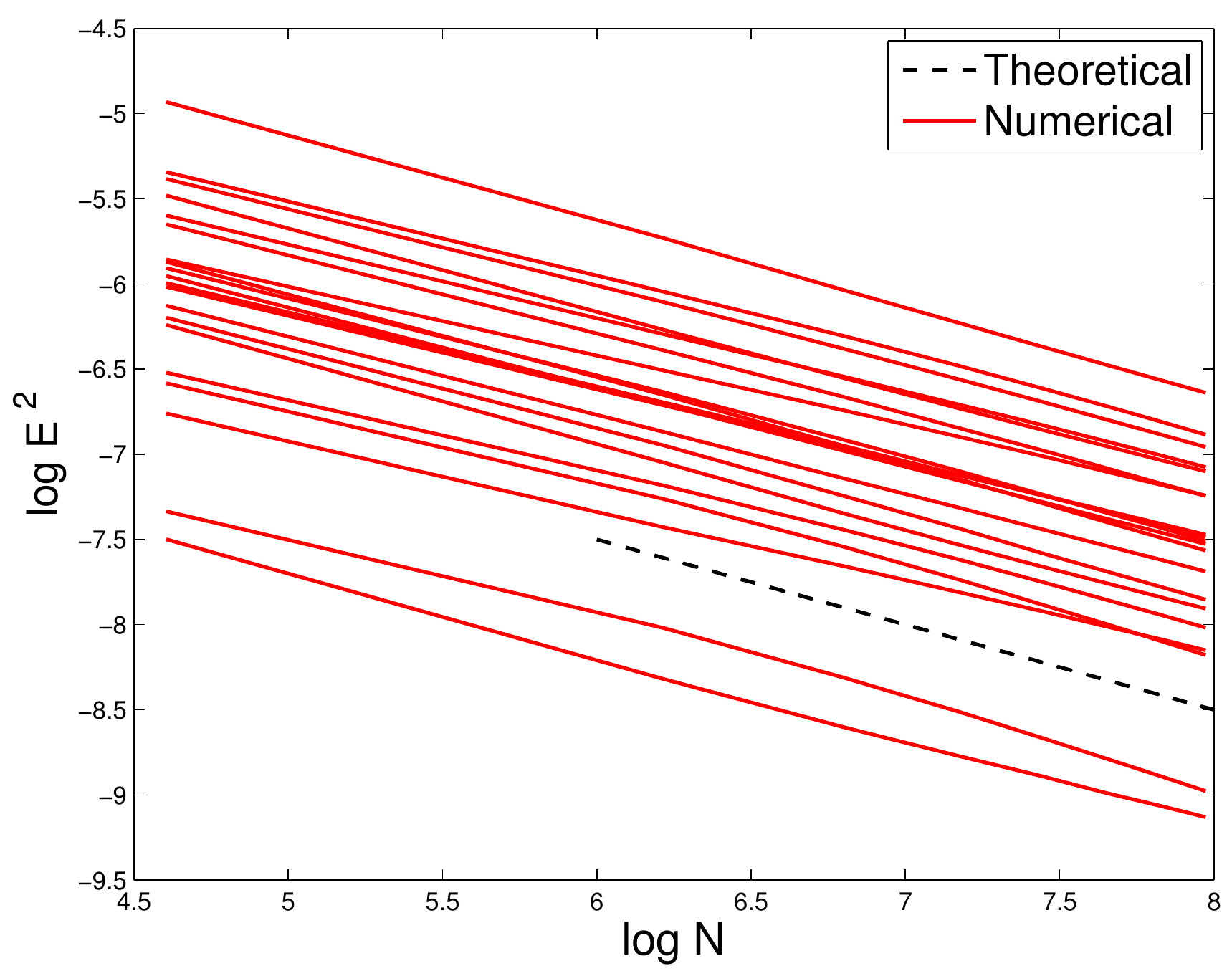}
\includegraphics[scale=0.32]{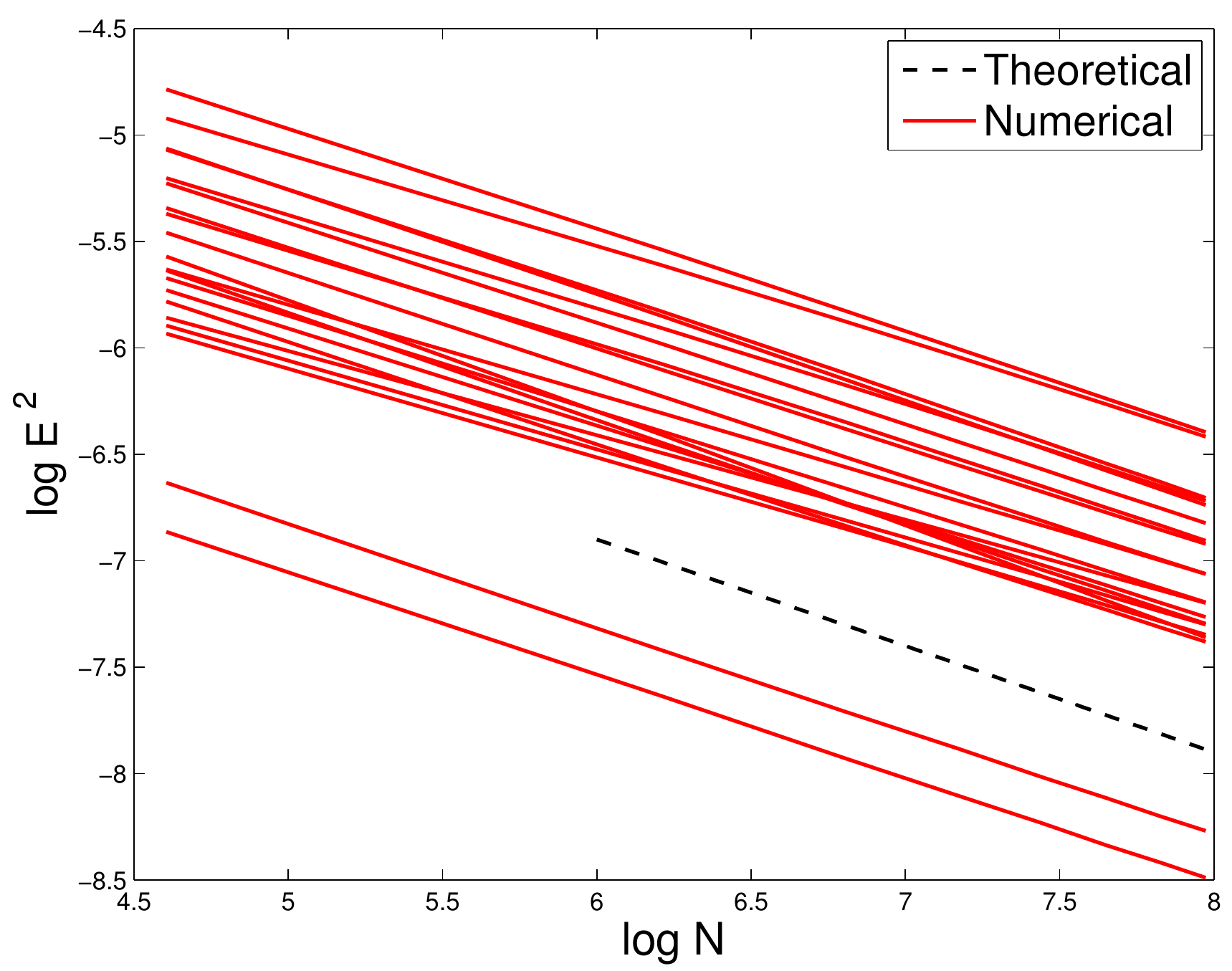}\\
\includegraphics[scale=0.32]{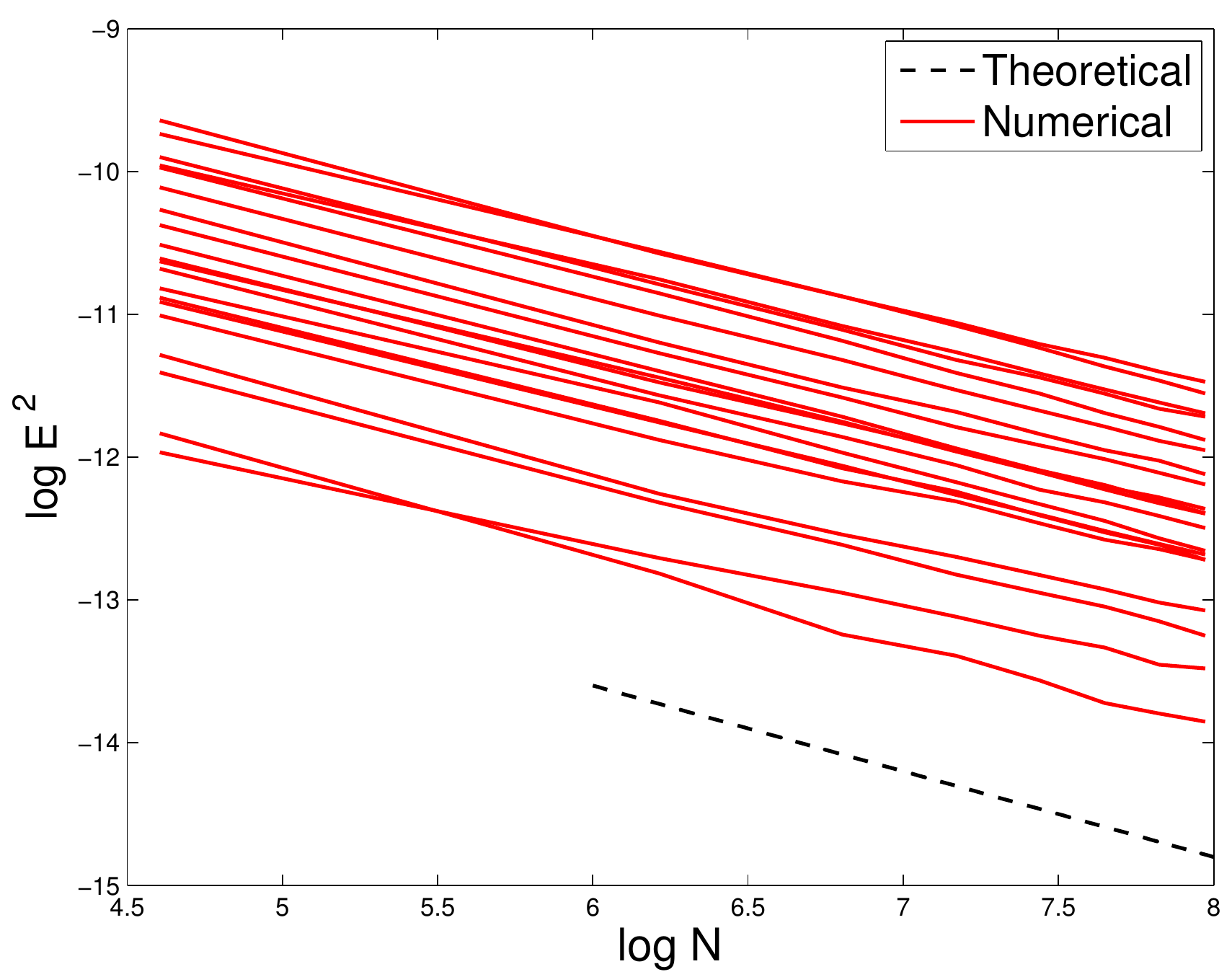}
\includegraphics[scale=0.32]{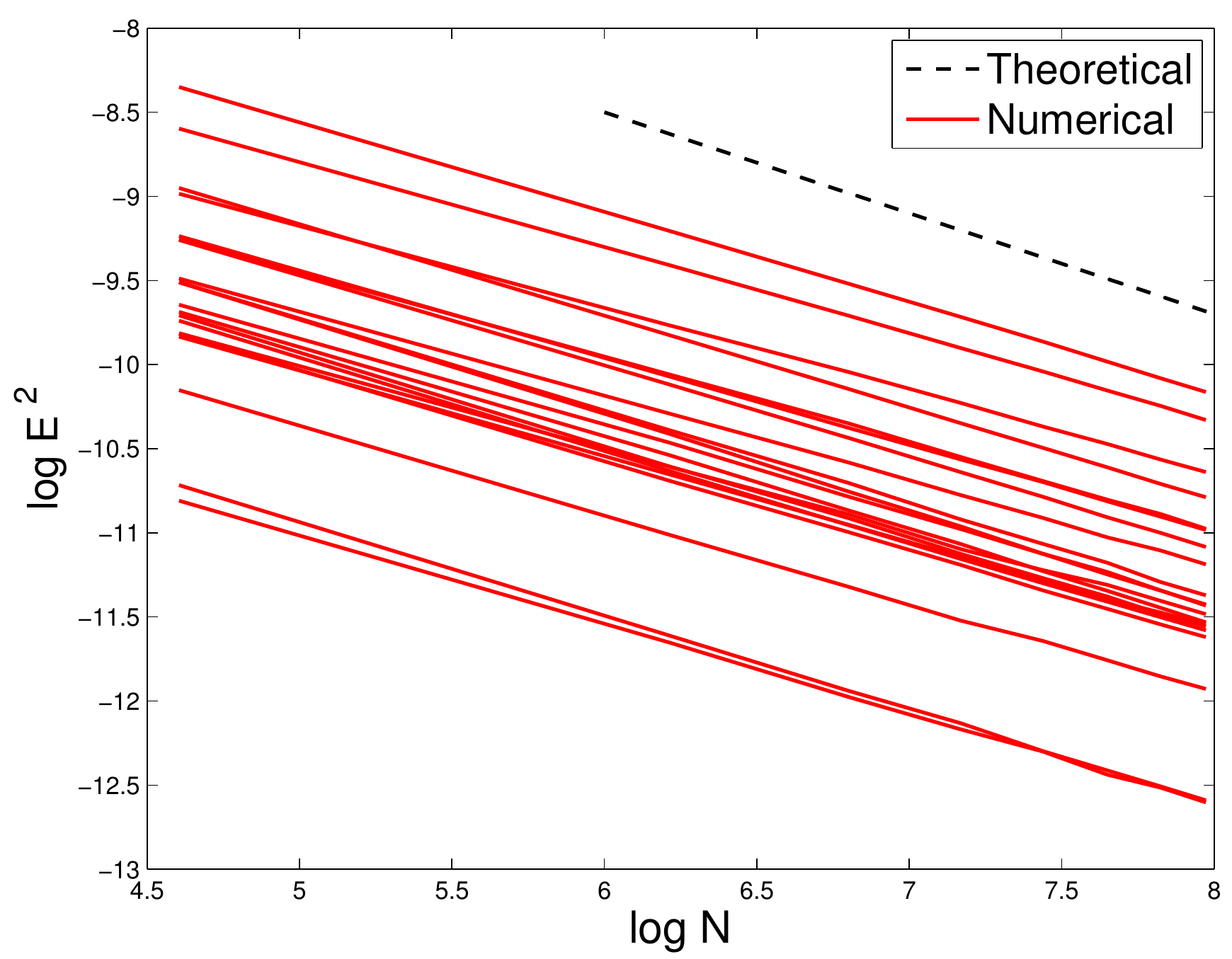}\\
\includegraphics[scale=0.32]{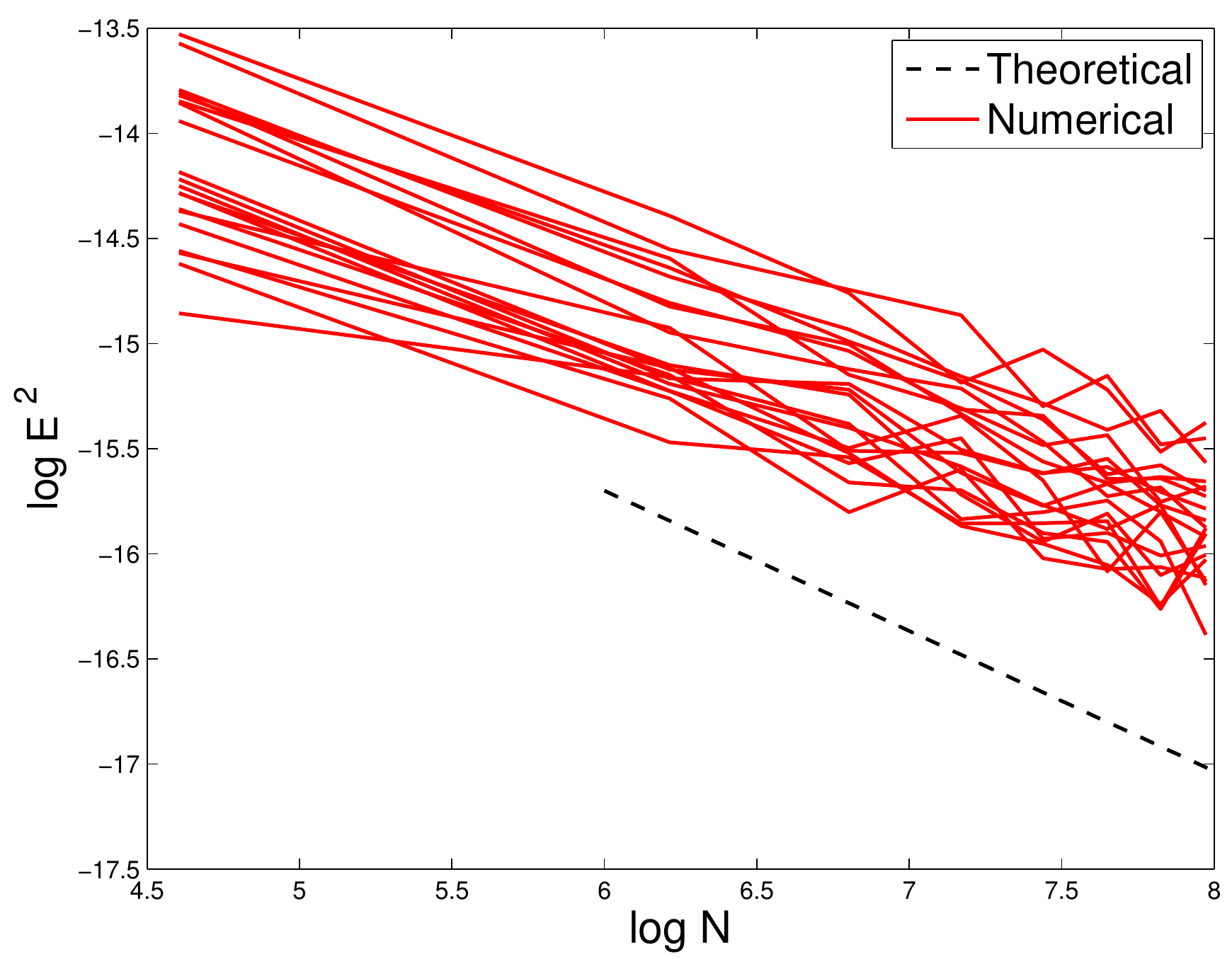}
\includegraphics[scale=0.32]{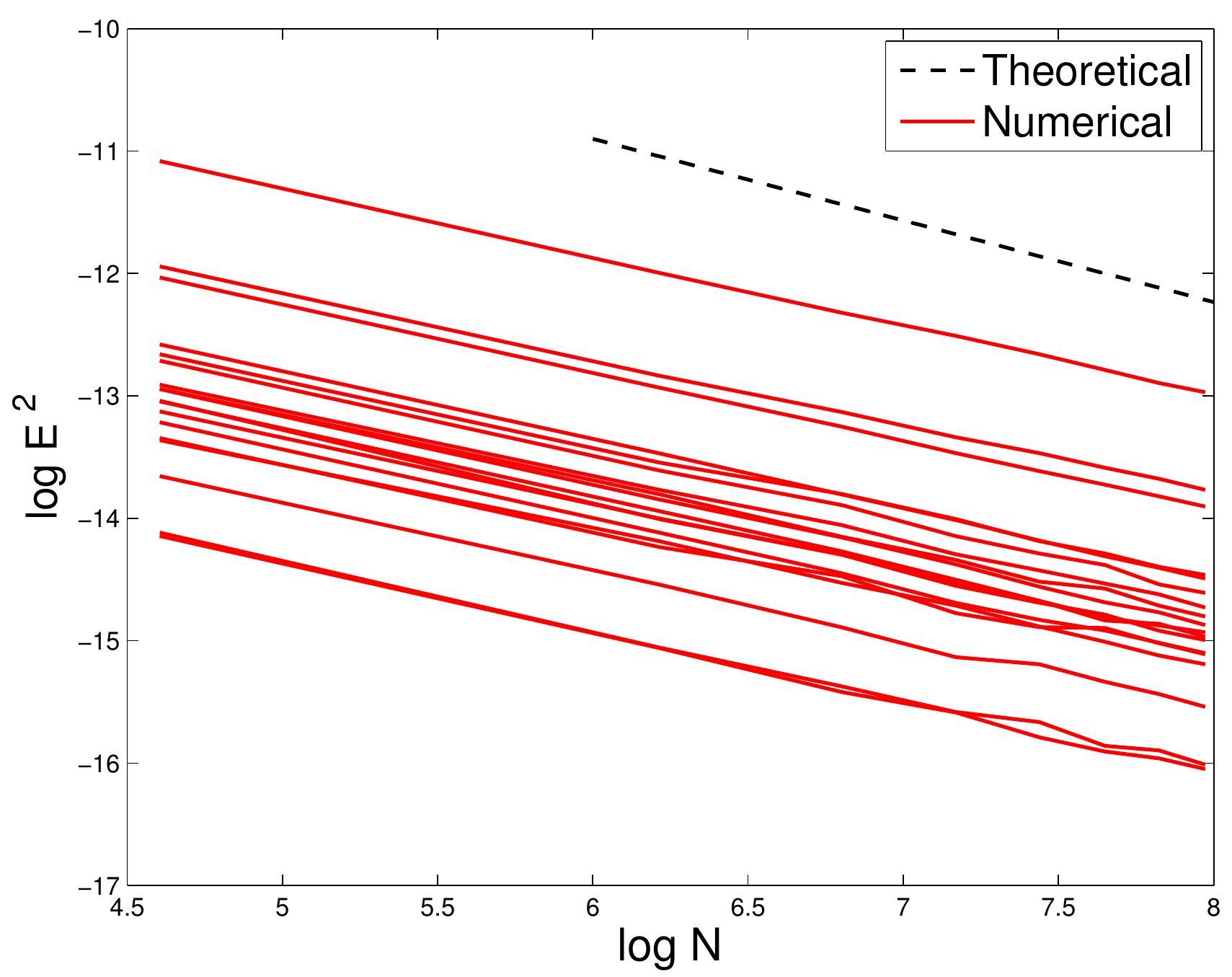}
 \caption{Convergence rates for 3DVAR (Left) and Kalman Filter (Right) with Data Model 1 and synthetic data generated from 20 different truths with regularity $H^{2s}$ with $s$ (from top to bottom) 1,2,3 and 4. } \label{Fig4}
\end{center}
\end{figure}

The verification of  Theorem \ref{thm_KFDM1} and Theorem \ref{thm_3DVARDM1} by means of Algorithm \ref{Al1} in the case of Data Model 1 is straightforward. For each of the set of synthetic data associated to each of the 20 truths $u^{\dagger}$ previously mentioned, we fix $\gamma=5\times 10^{-4}$. For each $N$ (with $N=\{100,\dots, 3000\}$) we run Algorithm \ref{Al1}, stop the schemes at $n=N$ and record the value of $\|\zeta_N-\true\|^2 $. In the right (resp. left) Figure \ref{Fig4} we display a plot of $\|\zeta_N-\true\|^2 $ vs $\log N$ for the Kalman filter (resp. 3DVAR) for each of the set of 20 experiments associated to different truths (red solid lines) generated as described above with (from top to bottom) $s=1,2,3,4$. From Theorem \ref{thm_KFDM1} we note that the corresponding slopes of the convergence rates should be approximately given by $-\frac{s}{s+1+a}$. For Theorem \ref{thm_3DVARDM1} there is an additional term of $\log{N}$, but this is of course negligible compared to the algebraic decay and we
ignore it for the purposes of this discussion.
For comparison, a line (black dotted) with slope $-\frac{s}{s+1+a}$ is displayed in Figure \ref{Fig4}.

\begin{figure}[htbp]
\begin{center}

\includegraphics[scale=0.32]{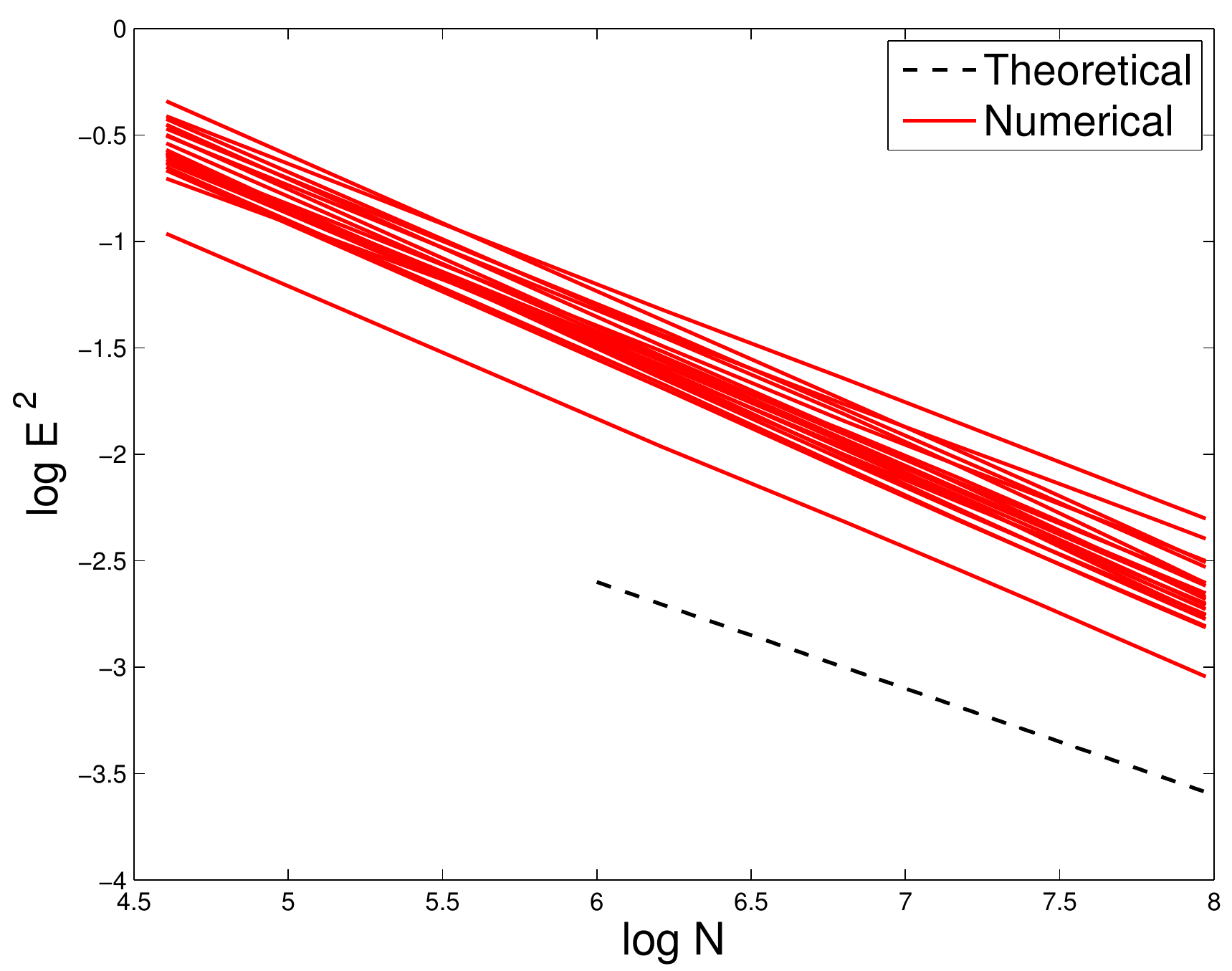}
\includegraphics[scale=0.32]{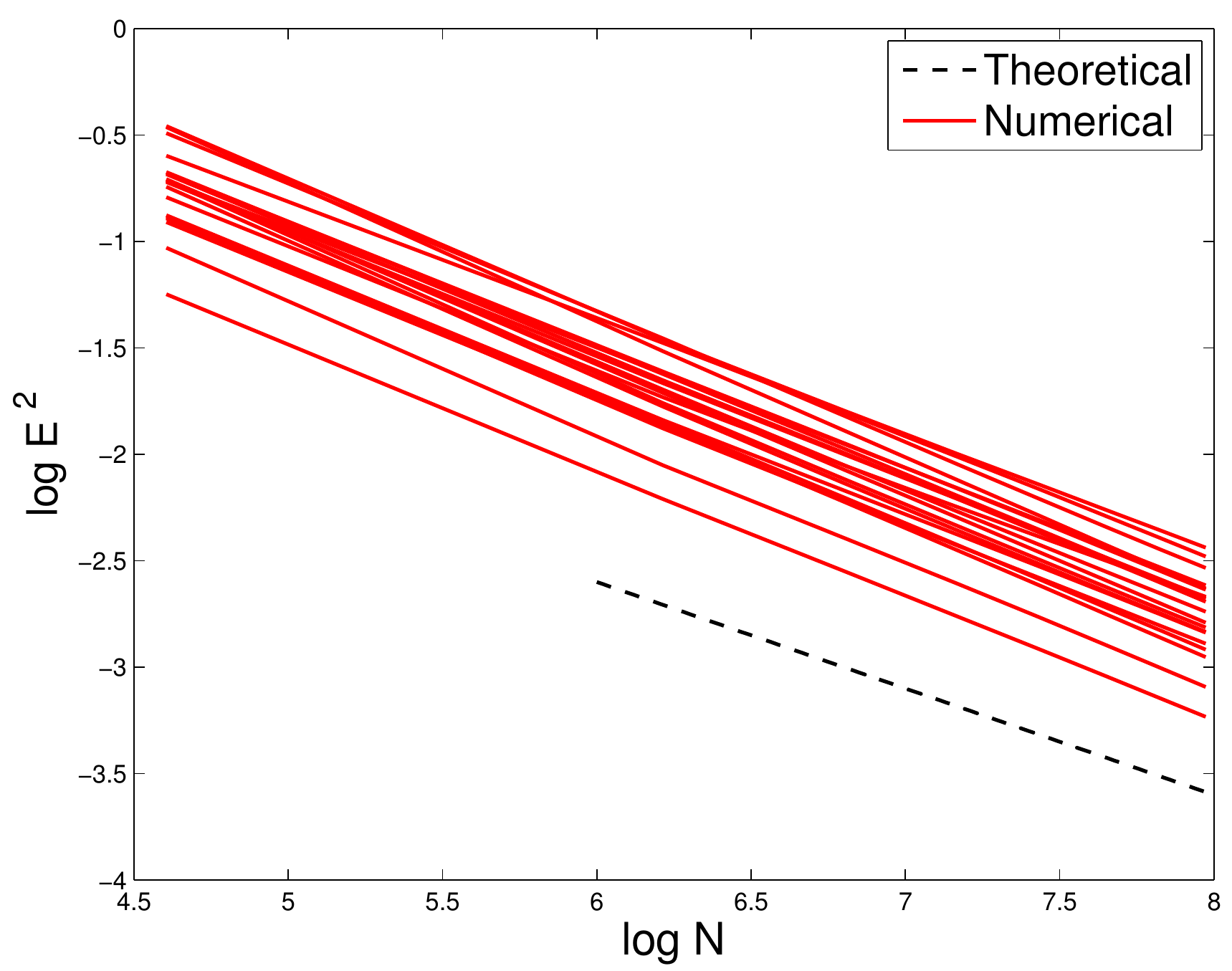}\\
\includegraphics[scale=0.32]{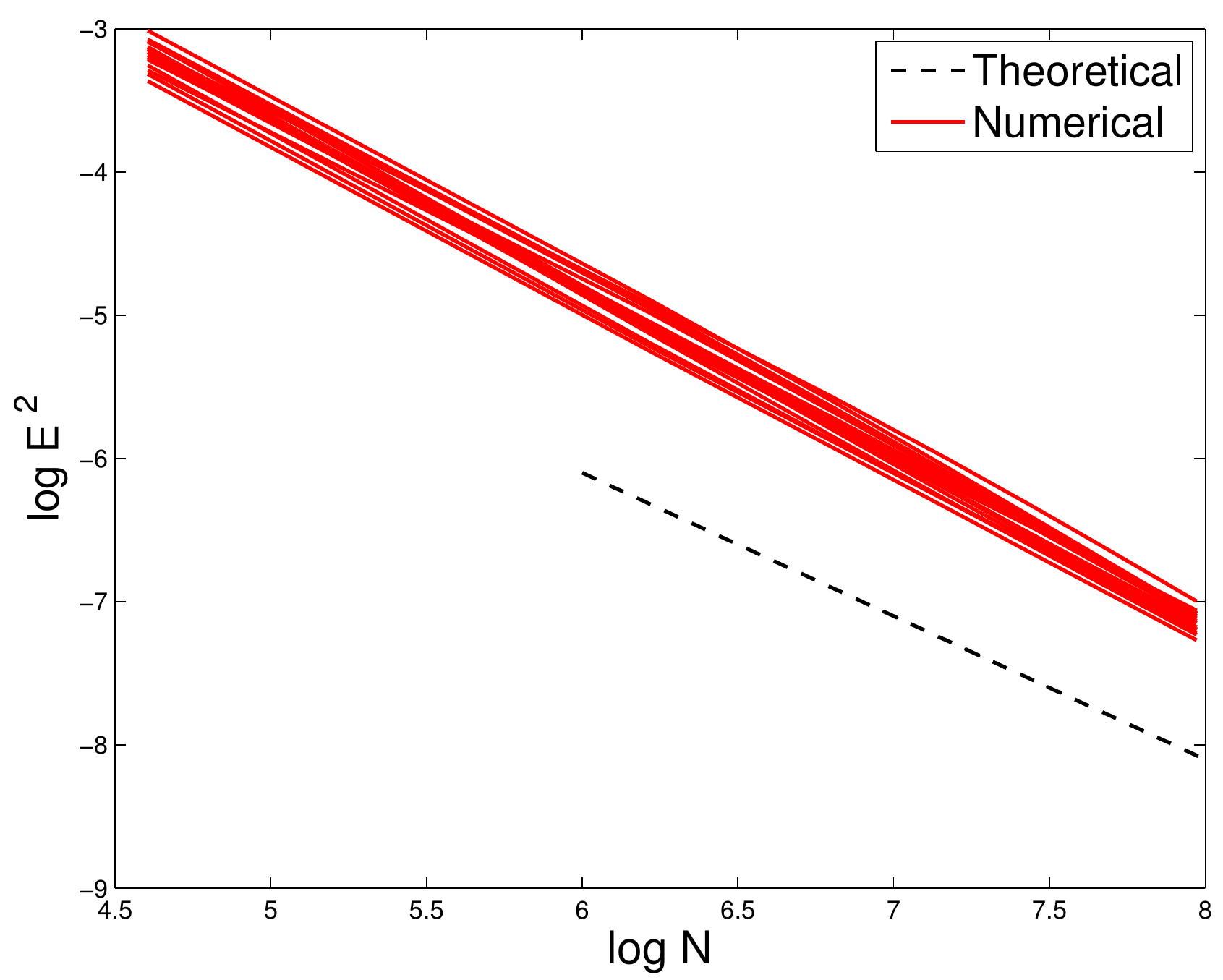}
\includegraphics[scale=0.32]{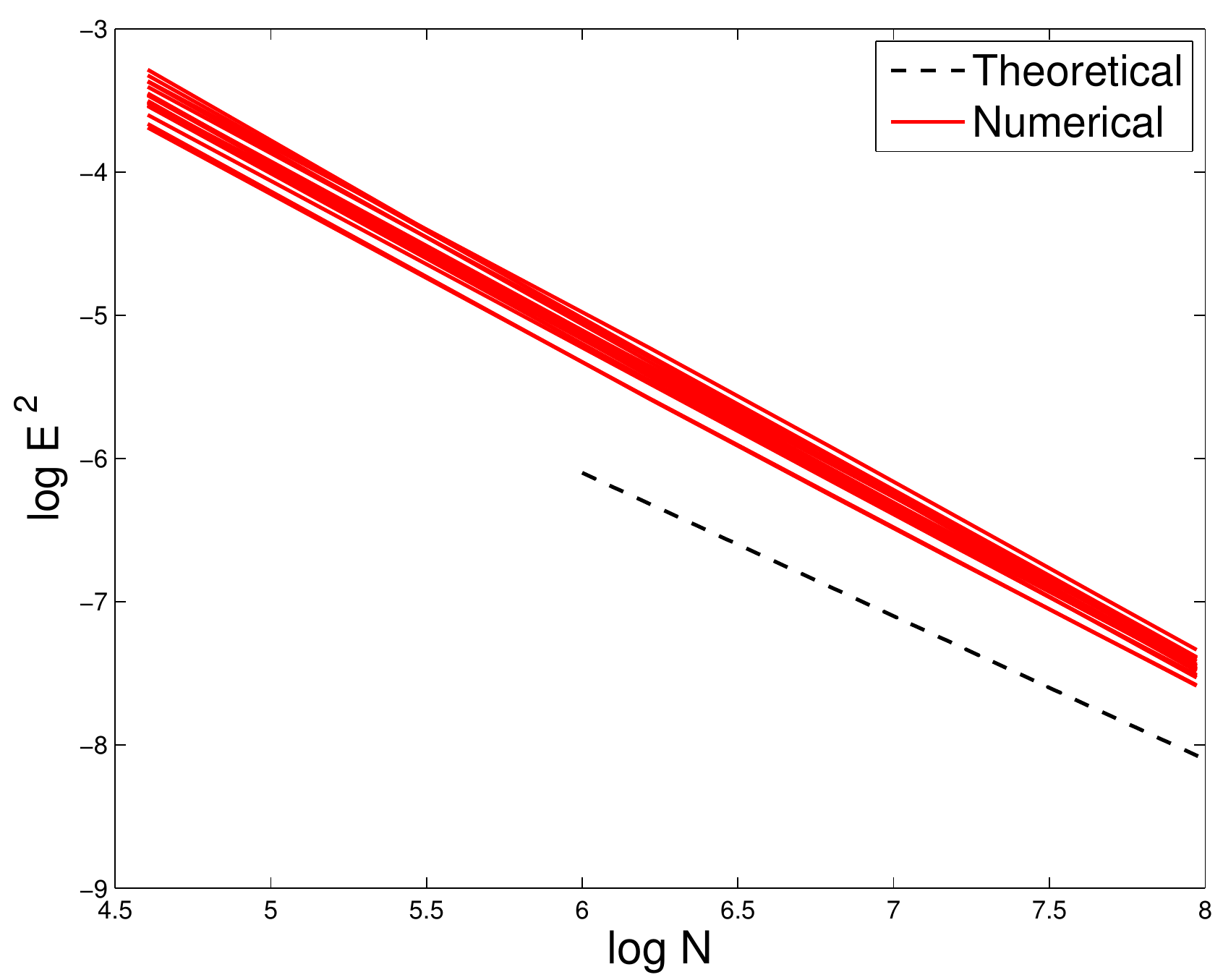}\\
\includegraphics[scale=0.32]{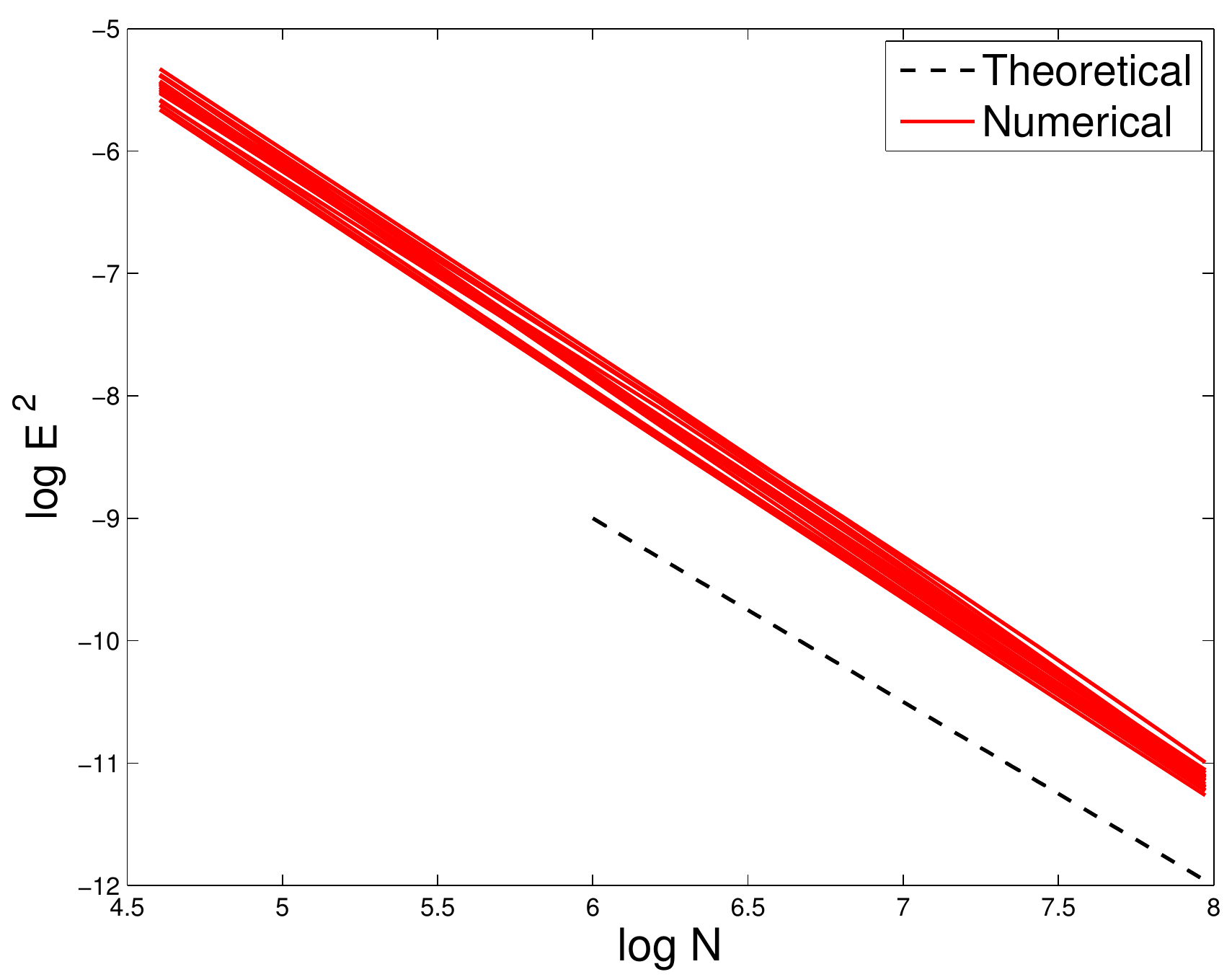}
\includegraphics[scale=0.32]{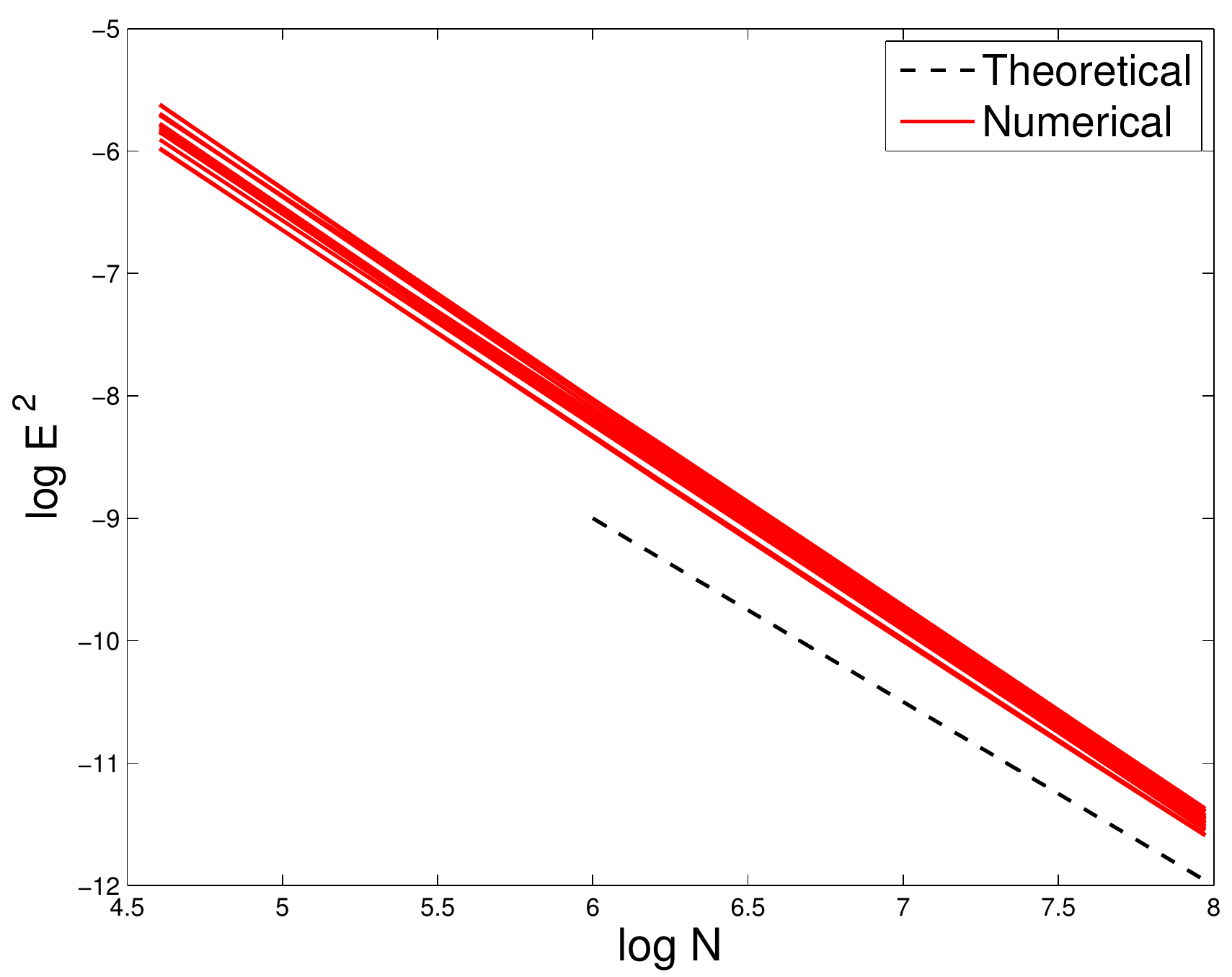}\\
\includegraphics[scale=0.32]{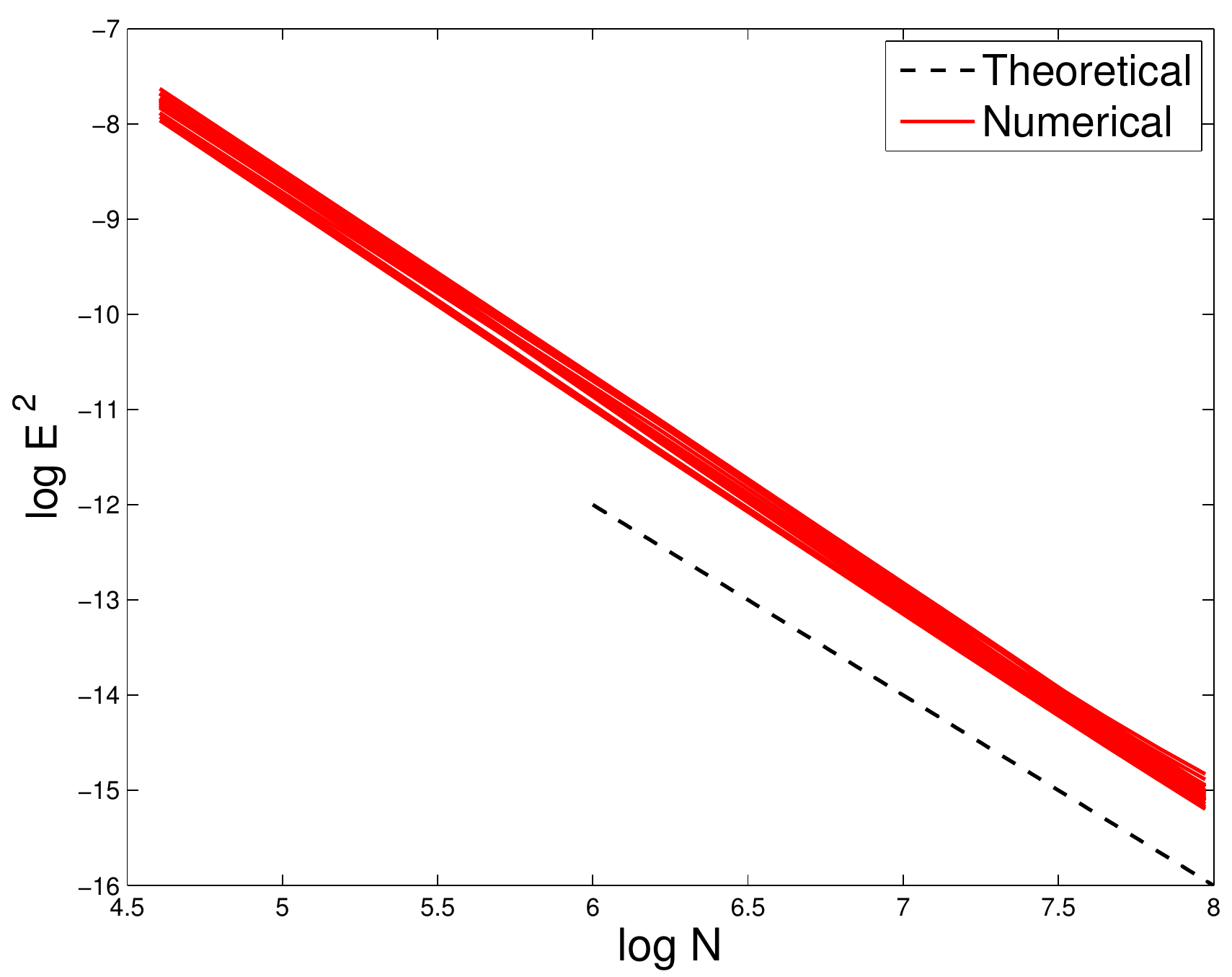}
\includegraphics[scale=0.32]{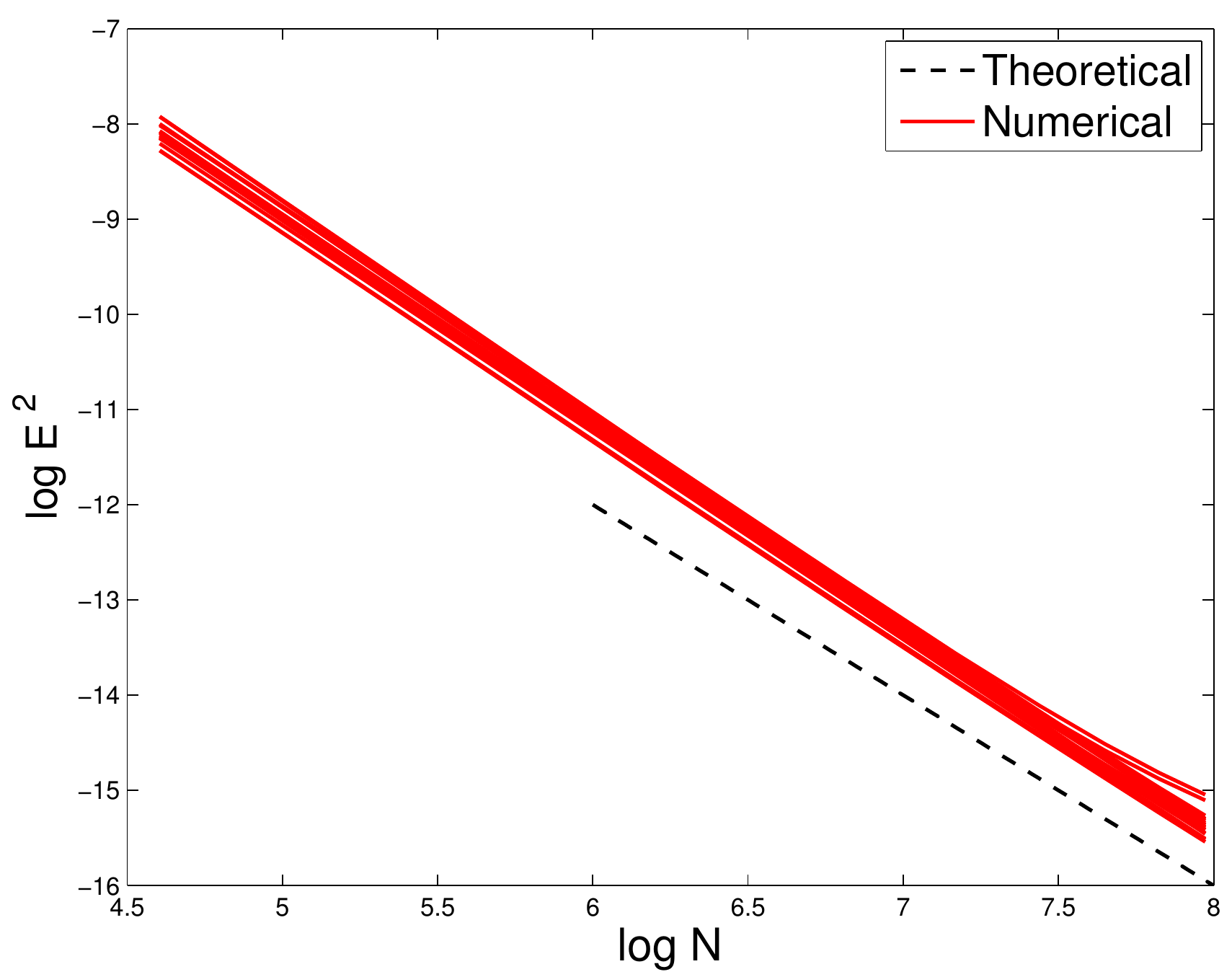}
 \caption{Convergence rates for 3DVAR (left) and Kalman Filter (right) with Data Model 2 and synthetic data generated from 20 different truths with regularity $H^{2s}$ with $s$ (from top to bottom) 1,2,3 and 4. } \label{Fig5}
\end{center}
\end{figure}

We now verify the convergence rates of Theorem \ref{thm_KFDM2} and Theorem \ref{thm_3DVARDM2}. Note first that in Algorithm \ref{Al1} for Data Model $2$ we define $N$ in terms of the given small noise $\gamma$, in order to obtain
convergence. However, for the purpose of the verification of the aforementioned convergence rates we define $\gamma$ in terms of $N$ by means of the same expressions. In other words, for each $N$ ($N=\{100,\dots, 3000\}$) we produce synthetic data (or each of the 20 truths) with $\eta\sim N(0,\gamma^2)$ and $\gamma=N^{-\frac{a+s+1}{2(a+1)}}$. We then run Algorithm \ref{Al1} and stop the schemes at $n=N$.  In the right (resp. left) Figure \ref{Fig5} we display a plot of $\|\zeta_N-\true\|^2 $ vs $\log N$ for the Kalman filter (resp. 3DVAR) for each of the set of 20 experiments associated to different truths (red solid lines) generated as before with (from top to bottom) $s=1,2,3,4$. We again include a line (black dotted) with slope of $-\frac{s}{a+1}$ which is the asymptotic behavior predicted by Theorems \ref{thm_KFDM2} and \ref{thm_3DVARDM2}.

We can clearly appreciate that, for $s$ satisfying Assumption \ref{assp_main} (iii) (i.e. $0<s\leq 3$), the numerical convergence rates fit very well the ones predicted by the theory. Note that the higher the regularity of the truth (i.e. the larger the $s$), the smaller the error. w.r.t the truth in the estimates. We note that for $s=4$, the aforementioned assumption is violated and, in the case of Data Model 1, the slopes of the numerical convergence rates are slightly smaller than the theoretical ones. In this case ($s=4$) there are
also fluctuations of the error w.r.t. the truth obtained with 3DVAR. These fluctuations may be associated with the fact that since for the error w.r.t. the truth is very small for sufficiently large iterations and for Data Model $1$ the noise level is fixed a priori (recall $\gamma=5\times 10^{-4}$). However, for
the Kalman filter these fluctuations are not so evident; presumably updating the covariance has a stabilizing effect.  For Data Model $2$, as $N$ increases, the corresponding $\gamma$ decreases and so these fluctuations in the error are non existent.

\section{Conclusions}

\begin{itemize}
\item We have presented filter based algorithms for the linear inverse problem,
based on introduction of an artificial dynamic. This results in methods
which are closely related to iterative Tikhonov-type regularization.
Two data scenarios are considered, one (Data Model $1$)
involving multiple realizations
of the data, with independent noise; the other (Data Model $1$)
involving a single realization of the data; both are relevant in applications.

\item We present theoretical results demonstrating convergence of the
algorithms in the two data scenarios. For multiple realizations of the
noisy the convergence is induced by the inherent averaging present
in the iterative method, and the link to the law of large numbers
and central limit theorem. For the single instance of data
the small observational noise limit must be considered.

\item For both Data Model $1$ and Data Model $2$ the Kalman Filter and 3DVAR
produced very similar results for relatively small $N$ ($N<100$). In practice it is clear that 3DVAR is preferable as the Kalman filter requires covariance
updates which may be impractical for large scale models. However, updating the covariance in the Kalman filter seems to have an stabilizing effect in the error w.r.t the truth.

\item For Data  Model $1$ the level of accuracy of the estimator is independent of the noise level. Moreover, the stability of the scheme is not conditioned to the early termination of the scheme. In contrast, for Data Model $2$ we need to stop at $n=N$ to avoid an increase in the error w.r.t the truth. Again, this illustrates that, whenever multiple instances of the data are available,
Data Model $1$ offers a more stable and accurate framework for solving the inverse problems under consideration.

\item The theoretical results from this work are verified numerically whenever the assumptions of the theory are satisfied.

\end{itemize}

\section*{Acknowledgements}

Andrew M. Stuart is funded by EPSRC (under the Programme Grant EQUIP), DARPA (under
EQUiPS) and ONR. Shuai Lu is supported by Special Funds for Major State Basic Research Projects of China (2015CB856003), NSFC (91130004, 11522108), Shanghai Science and
Technology Commission Grant (14QA1400400) and the Programme of Introducing Talents of Discipline
to Universities (number B08018), China.

\bibliographystyle{amsplain}

\appendix
\section*{Appendix:}

{\it Proof of Lemma \ref{lemma_Kalmanbias}.}
The proof follows the classic arguments on Tikhonov regularization with Hilbert scales, c.f. \cite[Ch8.4]{EHN1996}.
We recall (\ref{eq:CA}) in Assumption \ref{assp_main} $(ii)$ and rewrite it
in the form
\begin{equation*}
      c_1\|\Sigma_{0}^{\frac{a+1}{2}}x\| \leq \|A\Sigma_{0}^{\frac 1 2}x\|\leq c_2\|\Sigma_{0}^{\frac{a+1}{2}}x\|.
\end{equation*}
Notice that the definition of $B_0$ in \eqref{eq:B0} gives, since $X$
is a Hilbert space, and using Assumption \ref{assp_main} $(i)$
to ensure that $B_0$ and $B_0^*$ are well-defined bounded
linear operators,

$$\|A\Sigma_{0}^{\frac 1 2}x\|=\|B_0x\|=\|(B_0^*B_0)^{\frac 1 2}x\|.$$
Combining the two preceding displays we obtain
$$c_1\|\Sigma_{0}^{\frac{a+1}{2}}x\| \leq \|(B_0^*B_0)^{\frac 1 2}x\| \leq c_2\|\Sigma_{0}^{\frac{a+1}{2}}x\|$$
and a duality argument yield
$$c_2^{-1}\|\Sigma_{0}^{-\frac{a+1}{2}}x\| \leq \|(B_0^*B_0)^{-\frac 1 2}x\| \leq c_1^{-1}\|\Sigma_{0}^{-\frac{a+1}{2}}x\|$$
for any $x\in \mathcal{R}(\Sigma_{0}^{\frac{a+1}{2}})$.
Let $\theta\in[-1,1]$. Then the inequality of
Heinz \cite[Ch.8.4, pp. 213]{EHN1996} and an additional duality argument gives
\begin{equation}
\label{eq:CB}
c_1\|\Sigma_{0}^{\frac{\theta(a+1)}{2}}x\| \leq \|(B_0^*B_0)^{\frac \theta 2}x\| \leq c_2\|\Sigma_{0}^{\frac{\theta(a+1)}{2}}x\|,
\end{equation} which yields $\mathcal{R}\left((B_0^*B_0)^{\frac {\theta} {2}}\right)=\D{\Sigma_0^{-\frac{\theta(a+1)}{2}}}$.

Let Assumption \ref{assp_main} $(iii)$ be valid with $m_0-\true\in \D{\Sigma_{0}^{-\frac{s}{2}}}$, and define $z^{\dag}:=\Sigma_{0}^{-\frac 1 2}(m_0-\true) \in \D{\Sigma_{0}^{-\frac{s-1}{2}}}$.
Since $s-1 \in[-1,a+1]$, and consequently
$\theta=\frac{s-1}{a+1}\in (-1,1],$ we obtain from (\ref{eq:CB})
that $z^{\dag} \in R((B_0^*B_0)^{\frac {s-1} {2(a+1)}}).$ Furthermore
there exists a $\upsilon\in X$ such that
\begin{equation*}
 z^\dag=(B_0^*B_0)^{\frac {s-1} {2(a+1)}}\upsilon.
\end{equation*}
Noting that $\frac{1}{a+1}\in (0,1)$,
and employing (\ref{eq:CB}) with $\theta=\frac{1}{a+1}$, together with (\ref{eq:Rapproxi_km2}) and (\ref{eq:pos_rr}), we have
\begin{align*}
\|J_1\|^2 & = \|\Sigma_0^{\frac 1 2}r_{1,\frac{\alpha}{n}}(B_0^*B_0)\Sigma_0^{-\frac 1 2}(m_0-\true)\|^2 \\
& = \|(B_0^* B_0)^{\frac{1}{2(a+1)}}  r_{1,\frac{\alpha}{n}}(B_0^*B_0)z^{\dag}\|^2\\
& = \|(B_0^* B_0)^{\frac{1}{2(a+1)}}  r_{1,\frac{\alpha}{n}}(B_0^*B_0)(B_0^*B_0)^{\frac {s-1} {2(a+1)}}\upsilon\|^2\\
& = \|(B_0^* B_0)^{\frac{s}{2(a+1)}}  r_{1,\frac{\alpha}{n}}(B_0^*B_0) \upsilon\|^2 \\
& \leq \left(\frac{\alpha}{n}\right)^{\frac{s}{a+1}} \|\upsilon\|^2.
\end{align*}
In case of Assumption \ref{assp_main2} we insert $a =\frac{2p}{1+2\epsilon}$ and $s =\frac{2\beta}{1+2\epsilon}$.

\hfill$\Box$\newline

{\it Proof of Lemma \ref{lemma_Kalmanvariance}.}
Notice that
\begin{align*}
K_n & = C_{n-1} A^* (A C_{n-1} A^* + \gamma^2 I)^{-1} \\
 & = C_{n-1}^{\frac{1}{2}}  C_{n-1}^{\frac{1}{2}} A^* (A C_{n-1}^{\frac{1}{2}}C_{n-1}^{\frac{1}{2}} A^* + \gamma^2 I)^{-1} \\
 & = C_{n-1}^{\frac{1}{2}}  (C_{n-1}^{\frac{1}{2}} A^* A C_{n-1}^{\frac{1}{2}} + \gamma^2 I)^{-1} C_{n-1}^{\frac{1}{2}} A^*
\end{align*}
and
\begin{align*}
C_n & = (I-K_n A) C_{n-1} = \gamma^2 C_{n-1}^{\frac{1}{2}}  (C_{n-1}^{\frac{1}{2}} A^* A C_{n-1}^{\frac{1}{2}} + \gamma^2 I)^{-1} C_{n-1}^{\frac{1}{2}}.
\end{align*}
Thus we obtain
\begin{align}
C_n A^* = \gamma^2 K_n.  \label{eq:CnKn}
\end{align}
By virtue of (\ref{eq:prod}) and (\ref{eq:CnKn})  we derive
\begin{align*}
J_2 & =  \sum_{j=1}^{n-1} \left(\prod_{i=n-j}^{n-1}(I-K_{i+1}A)\right) K_{n-j} \eta_{n-j} + K_n \eta_n\\
& = \sum_{j=0}^{n-1} \left(C_n C_{n-j}^{-1} K_{n-j}\right)\eta_{n-j} \\
& = \sum_{j=0}^{n-1} \left(C_n A^*/\gamma^2\right)\eta_{n-j} \\
& = \sum_{j=0}^{n-1} \left( \left(C_0^{-1} + n\frac{A^* A}{\gamma^2}\right)^{-1} A^* /\gamma^2  \right) \eta_{n-j}.
\end{align*}
We denote $F:= \left(C_0^{-1} + n\frac{A^* A}{\gamma^2}\right)^{-1} A^* /\gamma^2$ and obtain
\begin{align*}
\E \|J_2\|^2 & = \sum_{j=0}^{n-1} \E \|F\eta_{n-j}\|^2 = n\gamma^2 \tr(F F^*).
\end{align*}
By the definition of $F$ and Assumption \ref{assp_main} $(i)$, $(iv)$ we obtain

\begin{align*}
F & = \left(C_0^{-1} + n\frac{A^* A}{\gamma^2}\right)^{-1} A^* /\gamma^2 \\
& = \Sigma_0^{\frac{1}{2}} (\alpha I  + n B_0^* B_0)^{-1} B_0^*
\end{align*}
and consequently derive
\begin{align*}
 \tr(FF^*) &=  \tr\left(\left( \Sigma_0^{\frac{1}{2}}(\alpha I + n B_0^*B_0)^{-1}B_0^*\right)\left(B_0 (\alpha I + nB_0^*B_0)^{-1}\Sigma_0^{\frac{1}{2}}\right)\right)\nonumber\\
  & =  \tr((\alpha I+ nB_0^*B_0)^{-2}B_0^*B_0\Sigma_0)\nonumber\\
  & \leq  \|(\alpha I+ nB_0^*B_0)^{-2}B_0^*B_0\|\tr(\Sigma_0)\nonumber\\
  & =  \frac{1}{\alpha^2}\|r_{2,\frac{\alpha}{n}}(B_0^*B_0)B_0^*B_0\|\tr(\Sigma_0)\nonumber\\
  & \leq  \frac{1}{\alpha^2}\frac{\alpha}{n}\tr(\Sigma_0)= \frac{1}{\alpha n}\tr(\Sigma_0).
\end{align*}
with the operator-valued function $r_{2,\frac{\alpha}{n}}(\lambda) := \left(\frac{\frac{\alpha}{n}}{\frac{\alpha}{n}+\lambda}\right)^2 = \left(\frac{\alpha}{\alpha+n\lambda}\right)^2$.
Such an observation then yields
\begin{align*}
\E \|J_2\|^2 = n\gamma^2 \tr(F F^*) \leq \frac{\gamma^2}{\alpha} \tr(\Sigma_0).
\end{align*}
Concerning Assumption \ref{assp_main2}, we further estimate, by
exploiting \cite[Lemma 8.2]{KVZ2011},
\begin{align*}
\tr(FF^*) & = \frac{1}{\alpha^2}\sum_{i=1}^{\infty} \frac{i^{-(1+2\epsilon+2p)-(1+2\epsilon)}}{\left(1+\frac{n}{\alpha} i^{-(1+2\epsilon+2p)}\right)^2} \\
& = \frac{1}{n} \sum_{i=1}^{\infty} \frac{\frac{n}{\alpha^2} i^{-4\epsilon-2p-2} }{\left(1+\frac{n}{\alpha} i^{-2\epsilon-2p-1}\right)^2} \\
& \asymp \frac{1}{n \alpha} \left(\frac{n}{\alpha}\right)^{-\frac{2\epsilon}{1+2\epsilon+2p}}
\end{align*}
and
\begin{align*}
\E\|J_2\|^2 \asymp \gamma^2 n^{-\frac{2\epsilon}{1+2\epsilon+2p}} \alpha^{-\frac{1+2p}{1+2\epsilon+2p}}.
\end{align*}

\hfill$\Box$\newline

{\it Proof of Lemma \ref{lemma_KalmanvarianceDM2}.}
Notice that for Data Model $2$, we derive
\begin{align*}
J_2 & =  \sum_{j=1}^{n-1} \left(\prod_{i=n-j}^{n-1}(I-K_{i+1}A)\right) K_{n-j} \eta_{n-j} + K_n \eta_n \\
& = \sum_{j=0}^{n-1}  \left( \left(C_0^{-1} + n\frac{A^* A}{\gamma^2}\right)^{-1} A^* /\gamma^2 \right)\eta_{n-j}\\
& = n F \eta
\end{align*}
which yields
\begin{align*}
\E \|J_2\|^2 & = n^2\gamma^2 \tr(F F^*).
\end{align*}
The remainder of the proof follows Lemma \ref{lemma_Kalmanvariance}.

\hfill$\Box$\newline

{\it Proof of Theorem \ref{thm_3DVARDM1Ass2}.}
As for the other theorems, the proof rests, of course, on the bias variance
decomposition, and then use of
Lemmas \ref{lemma_3DVARbias} and \ref{lemma_3DVARvarianceDM1}. This yields
\begin{align*}
\E\|\zeta_n - \true\|^2 \leq C \left(\frac{\alpha}{n}\right)^{\frac{2\beta}{1+2\epsilon+2p}} + C \gamma^2 \alpha^{-\frac{1+2p}{1+2\epsilon+2p}}
\end{align*}
and simultaneously
\begin{align*}
\E\|\zeta_n - \true\|^2 \leq C \left(\frac{\alpha}{n}\right)^{\frac{2\beta}{1+2\epsilon+2p}} + C \frac{\gamma^2 \ln n}{\alpha}.
\end{align*}
Choosing $\alpha = N^{\frac{2\beta}{1+2\beta+2p}}$ for the former inequality and $\alpha = N^{\frac{2\beta}{1+2\epsilon+2\beta+2p}}$ for the latter inequality, we conclude that, by stopping the iteration when $n=N$,
\begin{align*}
\E\|\zeta_N - \true\|^2 \leq C N^{-\frac{2p}{1+2\epsilon+2\beta+2p}}\ln N.
\end{align*}
\hfill$\Box$\newline

{\it Proof of Lemma \ref{lemma_3DVARbias}.}
Analogously to the proof of Lemma \ref{lemma_Kalmanbias}, it may
be shown that
  \begin{eqnarray*}
    \|I_1\|^2&=&   \|\Sigma_{0}^{\frac {1} {2}}\rr{B_0^*B_0}\Sigma_{0}^{-\frac 1 2}(\true-m_0)\|^2\\
       &\leq& \|(B_0^*B_0)^{\frac{1}{2(a+1)}}\rr{B_0^*B_0}z^{\dag}\|^2 \\
       &=& \|(B_0^*B_0)^{\frac{1}{2(a+1)}}\rr{B_0^*B_0}(B_0^*B_0)^{\frac {s-1} {2(a+1)}}\upsilon\|^2\\
        &= & \|(B_0^*B_0)^{\frac{s}{2(a+1)}}\rr{B_0^*B_0}\upsilon\|^2\\
        & \leq & \left(\frac \alpha n\right)^\frac{s}{a+1}\|\upsilon\|^2.
  \end{eqnarray*}
The final inequality follows from the asymptotic behavior of
$\rr{\lambda}$, established, for example. in  \cite[Ch.2, pp. 63]{LP2013}.
In the case of Assumption \ref{assp_main2}
we insert $a =\frac{2p}{1+2\epsilon}$ and $s =\frac{2\beta}{1+2\epsilon}$.

\hfill$\Box$\newline

{\it Proof of Lemma \ref{lemma_3DVARvarianceDM1}.}
We denote $F_j = (I-\mathcal{K}A)^j \mathcal{K}$ and obtain
\begin{align*}
I_2 = \sum_{j=0}^{n-1} F_j \eta_{n-j}.
\end{align*}
Furthermore, we derive
\begin{align*}
\E\|I_2\|^2 = \sum_{j=0}^{n-1}\E \|F_j \eta_{n-j}\|^2 =\gamma^2 \sum_{j=0}^{n-1} \tr(F_j F_j^*)
\end{align*}
and
\begin{align*}
\sum_{j=0}^{n-1} \tr \left(F_j F_j^*\right) & = \sum_{j=0}^{n-1} \tr\left((I-\mathcal{K}A)^j \mathcal{K} \mathcal{K}^* \big((I-\mathcal{K}A)^*\big)^j\right) \\
& = \frac{1}{\alpha^2}\sum_{j=0}^{n-1} \tr\left( \Sigma_0^{\frac{1}{2}} r_{2j+2,\alpha}(B_0^* B_0) B_0^* B_0 \Sigma_0^{\frac{1}{2}} \right) \\
& \leq  \frac{1}{\alpha^2}\sum_{j=0}^{n-1} \|r_{2j+2,\alpha}(B_0^* B_0) B_0^* B_0\| \tr(\Sigma_0) \\
& \leq  \frac{\tr(\Sigma_0)}{\alpha^2}\sum_{j=0}^{n-1} \frac{\alpha}{2j+2} \\
& \asymp C \ln(n)\frac{\tr(\Sigma_0)}{ \alpha}.
\end{align*}
Thus,
\begin{align*}
\E\|I_2\|^2 \leq C\frac{\ln(n)\gamma^2}{\alpha} \tr(\Sigma_0).
\end{align*}
For Assumption \ref{assp_main2} we need to estimate $\tr(F_jF_j^*)$ carefully.
We substitute the given decay rate of the different eigenvalues, to obtain
\begin{align*}
\tr(F_jF_j^*) & = \frac{1}{\alpha^2}\tr\left( \Sigma_0^{\frac{1}{2}} r_{2j+2,\alpha}(B_0^* B_0) B_0^* B_0 \Sigma_0^{\frac{1}{2}} \right) \\
&  = \frac{1}{\alpha^2}\sum_{i=1}^{\infty} \frac{i^{-2-4\epsilon-2p}}{\left(1+\frac{1}{\alpha} i^{-1-2\epsilon-2p}\right)^{2j+2}} \\
& \leq \frac{1}{\alpha^2}\sum_{i=1}^{\infty} \frac{i^{-2-4\epsilon-2p}}{\left(1+\frac{j+1}{\alpha} i^{-1-2\epsilon-2p}\right)^2}.
\end{align*}
By arguments similar to those used
in the proof of Lemma \ref{lemma_Kalmanvariance} (or \cite[Lemma 8.2]{KVZ2011}), we further estimate
\begin{align*}
 \frac{1}{\alpha^2}\sum_{i=1}^{\infty} \frac{i^{-2-4\epsilon-2p}}{\left(1+\frac{j+1}{\alpha} i^{-1-2\epsilon-2p}\right)^2} & \asymp \left(\frac{1}{j+1 } \right)^{1+\frac{2\epsilon}{1+2\epsilon+2p}} \alpha^{-\frac{1+2p}{1+2\epsilon+2p}}
\end{align*}
and
\begin{align*}
\E\|I_2\| \leq C \gamma^2 \alpha^{-\frac{1+2p}{1+2\epsilon+2p}} \sum_{j=0}^{n-1} \left(\frac{1}{j+1 } \right)^{1+\frac{2\epsilon}{1+2\epsilon+2p}}
\end{align*}
where the summation term in the right-hand side is bounded.

On the other hand, we can also estimate
\begin{align*}
\tr(F_jF_j^*) & = \frac{1}{\alpha^2}\tr\left( \Sigma_0^{\frac{1}{2}} r_{2j+2,\alpha}(B_0^* B_0) B_0^* B_0 \Sigma_0^{\frac{1}{2}} \right) \\
&  = \frac{1}{\alpha^2}\sum_{i=1}^{\infty} \frac{i^{-2-4\epsilon-2p}}{\left(1+\frac{1}{\alpha} i^{-1-2\epsilon-2p}\right)^{2j+2}} \\
& \leq \frac{1}{\alpha^2}\sum_{i=1}^{\infty} \frac{i^{-2-4\epsilon-2p}}{\left(1+\frac{2(j+1)}{\alpha} i^{-1-2\epsilon-2p}\right)} \\
& < \frac{1}{2(j+1) \alpha} \sum_{i=1}^{\infty} i^{-1-2\epsilon} \\
& \leq C \frac{1}{(j+1) \alpha}
\end{align*}
and
\begin{align*}
\E\|I_2\| \leq C \frac{\gamma^2 \ln n}{\alpha}.
\end{align*}

\hfill$\Box$\newline

{\it Proof of Lemma \ref{lemma_3DVARvarianceDM2}.}
Since $\eta_n = \eta$ in this case, we derive
\begin{align*}
I_2 = \sum_{j=0}^{n-1} (I-\mathcal{K}A)^j \mathcal{K} \eta
\end{align*}
and by operator-valued calculation we obtain
\begin{align*}
\sum_{j=0}^{n-1} (I-\mathcal{K}A)^j \mathcal{K}  & = \frac{1}{\alpha} \Sigma_0^{\frac{1}{2}} \sum_{j=0}^{n-1} \left(\alpha(\alpha I + B_0^* B_0)^{-1}\right)^{j+1} B_0^* \\
& = \Sigma_0^{\frac{1}{2}} q_{n,\alpha}(B_0^* B_0) B_0^*
\end{align*}
where $q_{n,\alpha}(\lambda) := \frac{1}{\lambda}\left(1-\frac{\alpha^n}{(\alpha+\lambda)^n}\right)$.
Thus we obtain, by the asymptotic behavior of $q_{n,\alpha}(\lambda)$
derived in \cite[Ch.2, pp. 64]{LP2013},
\begin{align*}
\E\|I_2\|^2 & = \gamma^2 \tr\big(\Sigma_0^{\frac{1}{2}} q_{n,\alpha}(B_0^* B_0) B_0^* B_0 q_{n,\alpha}(B_0^* B_0) \Sigma_0^{\frac{1}{2}} \big ) \\
& \leq \gamma^2 \|q_{n,\alpha} (B_0^* B_0) (B_0^* B_0)^{\frac{1}{2}}\|^2 \tr(\Sigma_0) \\
& \leq \frac{n\gamma^2}{\alpha} \tr(\Sigma_0).
\end{align*}
\hfill$\Box$\newline

{\it Proof of Theorem \ref{thm_3DVARvariant}.}
Similar to the Kalman filter method and 3DVAR, by Assumption \ref{assp_main} $(i)$-$(iii)$, we obtain the bias error estimate
\begin{align*}
\|\mathcal{I}_1\| & \leq \left\|\Sigma_0^{1/2} \prod_{j=1}^n \frac{\alpha_j I}{B_0^* B_0 +\alpha_j I} \Sigma_0^{1/2} \epsilon_0 \right\|^2 \\
& \leq \left\|\prod_{j=1}^n \left(\frac{\alpha_j I}{B_0^*B_0 +\alpha_j I}\right) \left(B_0^*B_0 \right)^{\frac{s}{2(a+1)}} \upsilon\right\|^2.
\end{align*}
Now we need upper bounds of the following operator-valued function
\begin{align}
f_{n,v}(\lambda)  = \lambda^v \prod_{j=1}^n \frac{\alpha_j}{\lambda+\alpha_j}, \quad \lambda\in (0,\infty), \quad n> v>0. \label{eq:f_nv}
\end{align}
Define $\sigma_n :=\sum_{j=1}^n \frac{1}{\alpha_j}$ and assume the sequence $\{\alpha_j\}_{j=1}^n$ satisfying
\begin{align}
\frac{1}{\alpha_n} \leq \tilde{c} \sigma_{n-1} \label{eq:sequence}
\end{align}
with a constant $\tilde{c}$. Then the results in \cite{HG1998} yield
\begin{align*}
\|\mathcal{I}_1\| \leq C \sigma_n^{-\frac{s}{a+1}}, \quad n>1.
\end{align*}

It remains to estimate the $\mathcal{I}_2$ term. Define
\begin{align*}
\mathcal{F}_j : & = \prod_{i=n-j}^{n-1} \left(I-\mathcal{K}_{i+1} A\right) \mathcal{K}_{n-j} \\
& = \Sigma_0^{1/2} \prod_{i=n-j}^{n-1} \left(\frac{\alpha_{i+1} I}{B_0^* B_0 + \alpha_{i+1} I}\right) \frac{1}{B_0^* B_0 +\alpha_{n-j} I} B_0^*.
\end{align*}
Then we obtain
\begin{align*}
\mathcal{I}_2 & = \sum_{j=1}^{n-1} \mathcal{F}_j \eta_{n-j} + K_n \eta_n
\end{align*}
which yields
\begin{align*}
\mathbb{E} \|\mathcal{I}_2\|^2 &  = \gamma^2 \sum_{j=1}^{n-1} \tr(\mathcal{F}_j \mathcal{F}_j^*) + \gamma^2 \tr(K_n K_n^*).
\end{align*}

Notice that for any $j\geq 1$
\begin{align*}
\mathcal{F}_j \mathcal{F}^*_j & = \frac{1}{\alpha_{n-j}^2} \Sigma_0^{1/2} \prod_{i=n-j}^{n-1} \left(\frac{\alpha_{i+1} }{B_0^* B_0+\alpha_{i+1} I}\right)^2 \left(\frac{\alpha_{n-j} }{B_0^* B_0 +\alpha_{n-j} I}\right)^2 B_0^* B_0 \Sigma_0^{1/2}
\end{align*}
we derive,
\begin{align*}
\tr(\mathcal{F}_j \mathcal{F}_j^*) & \leq \frac{\tr(\Sigma_0)}{\alpha_{n-j}^2} \left\|\prod_{i=n-j}^{n-1} \left(\frac{\alpha_{i+1} }{B_0^* B_0+\alpha_{i+1} I}\right)^2 \left(\frac{\alpha_{n-j} }{B_0^* B_0 +\alpha_{n-j} I}\right)^2 B_0^* B_0\right\| \nonumber \\
& \leq \frac{\tr(\Sigma_0)}{\alpha_{n-j}^2} \left\|\left(\frac{\alpha_{n-j} }{B_0^* B_0 +\alpha_{n-j} I}\right)^2 B_0^* B_0\right\| \nonumber \\
& \leq  \frac{1}{2 \alpha_{n-j} } \tr(\Sigma_0)
\end{align*}
and
\begin{align*}
\tr(K_n K_n^*) & \leq \frac{1}{2 \alpha_n } \tr(\Sigma_0).
\end{align*}
A rough variance estimate for the variant method is
\begin{align*}
\mathbb{E}\|v_n - u^{\dag}\|^2 & = \|\mathcal{I}_1\|^2 + \mathbb{E}\|\mathcal{I}_2\|^2  \\
& \leq C \left(\sigma_{n}^{-\frac{s}{a+1}} + \gamma^2\tr(\Sigma_0)  \sigma_{n} \right).
\end{align*}
The first term vanishes but the second term blows up when $n\rightarrow \infty$ and $\sigma_n \rightarrow \infty$.

To further investigate the blow up, we consider a special geometric sequence $\alpha_n = \alpha q^{n-1}$ with $0<q<1$. Thus, we have
\begin{align*}
\sigma_n = \alpha^{-1} q^{1-n} \frac{1-q^n}{1-q} \geq \alpha^{-1} q^{1-n} = q/\alpha_{n+1}
\end{align*}
and (\ref{eq:sequence}) is satisfied with $\tilde{c} = 1/q$. Actually, we derive
\begin{align*}
\frac{1}{\alpha_{n-j}} + \sigma_{n}-\sigma_{n-j} =\frac{1}{\alpha_{n-j}} + \frac{1}{\alpha_{n-j+1}} +\ldots+ \frac{1}{\alpha_{n}} = \alpha^{-1} q^{1-n} \frac{1-q^{j+1}}{1-q} \geq  \alpha^{-1} q^{1-n} = q/\alpha_{n+1}.
\end{align*}
Thus the results in \cite{HG1998} refine, by the asymptotic behavior of (\ref{eq:f_nv}),
\begin{align*}
\tr(\mathcal{F}_j \mathcal{F}_j^*) & \leq \frac{\tr(\Sigma_0)}{\alpha_{n-j}^2} \left\|\prod_{i=n-j}^{n-1} \left(\frac{\alpha_{i+1} }{B_0^* B_0+\alpha_{i+1} I}\right)^2 \left(\frac{\alpha_{n-j} }{B_0^* B_0 +\alpha_{n-j} I}\right)^2 B_0^* B_0\right\| \nonumber \\
& \leq \frac{1}{\alpha^2_{n-j}} \Big(\frac{1}{\alpha_{n-j}}+ \sigma_n-\sigma_{n-j}\Big)^{-1} \tr(\Sigma_0) \\
& = \alpha^{-1} q^{1+j-n} \left(1+q^{-j} \frac{1-q^{j}}{1-q}\right)^{-1} \tr(\Sigma_0)\\
& \leq \alpha^{-1} q^{1+2j - n} \tr(\Sigma_0)
\end{align*}
and we derive
\begin{align*}
\mathbb{E}\|\mathcal{I}_2\|^2 & \leq \frac{q \gamma^2 \tr(\Sigma_0)}{ \alpha } q^{-n} \left(\sum_{j=1}^{n-1}q^{2j} + 1\right).
\end{align*}
Summing up, for the geometric sequence, we obtain
\begin{align*}
\mathbb{E}\|v_n-u^{\dag}\|^2 \leq C\left(q^{\frac{s}{a+1}n } + \gamma^2  \tr(\Sigma_0) q^{-n }\right).
\end{align*}
The second term blows up exponentially when $n$ goes to infinity.

\hfill$\Box$\newline

\end{document}